\documentclass[10pt]{article}
\usepackage[usenames]{color} 
\usepackage{amssymb} 
\usepackage{amsmath,amsthm,graphicx} 
\usepackage[utf8]{inputenc} 
\usepackage{hyperref}
\usepackage{wrapfig}
\usepackage{tikz}
\usepackage{bbm}
\usepackage{algorithm,algorithmic}
\usepackage{subfigure}
\usepackage{comment}
\usepackage{caption}
\usepackage[capitalise,noabbrev]{cleveref}
\usepackage{todonotes}
\usepackage{tikz}

\usepackage{pgfplots}
\usepackage{pgfplotstable}
\usetikzlibrary{fillbetween}
\usepgfplotslibrary{fillbetween}

\usetikzlibrary{shapes.misc}
\usetikzlibrary{math} 

\usetikzlibrary{arrows}

\def \R {\mathbb{R}}

\definecolor{lightgreen}{rgb}{0,1,0}
\def \E {\mathcal{E}}
\newcommand{\est}[2]{\left\langle #1,#2 \right\rangle_{\E}}
\def \NFE {{N_h}}
\def \NRB {{N_{RB}}}
\def \NEIM {{N_{EIM}}}
\def \CFL {\texttt{CFL}}

\newcommand{\mass}{\mathbb M}
\newcommand{\boldeta}{\boldsymbol{\eta}}
\newcommand{\bq}{\boldsymbol{q}}

\newcommand{\diff}{\mathcal{J}}
\newcommand{\dispersion}{\mathbb{T}^t}
\newcommand{\deriv}{\mathbb{D}}
\newcommand{\iden}{\mathbb{I}}
\newcommand{\sponge}{\mathbb{S}}

\newcommand{\iwg}{{\textrm{iwg}}}

\newcommand{\NLRHS}{\mathcal{N}}

\graphicspath{{figures/}}

\newtheorem{remark}{Remark}
\newcommand{\footremember}[2]{%
	\footnote{#2}
	\newcounter{#1}
	\setcounter{#1}{\value{footnote}}%
}

\title{Model order reduction strategies for   weakly dispersive waves} 
\author{Davide Torlo\footremember{sissa}{Mathematics Area, mathLab, SISSA, via Bonomea 265, I-34136 Trieste, Italy}\thanks{Corresponding author: \href{mailto:davide.torlo@sissa.it}{davide.torlo@sissa.it}}\, and Mario Ricchiuto\footremember{inria}{Team CARDAMOM, Inria Bordeaux Sud-Ouest, 200 Avenue de la Vieille Tour, 33405 Talence, France}}

\date{\today}

\begin{document}

\maketitle
\abstract{We focus on the numerical modelling of  water waves by means of depth averaged models. 
We consider in particular  PDE systems which consist in a nonlinear hyperbolic model plus 
a  linear dispersive perturbation involving an elliptic operator.
We propose  two strategies to construct reduced order models for these problems, with  the main focus being the control of the overhead 
related to the inversion of the elliptic operators, as well as the robustness with respect to variations of the flow  parameters.
In a first approach, only a linear  reduction strategies is applied only to the elliptic component, while the computations of the nonlinear fluxes are still performed explicitly. 
This hybrid approach, referred to as pdROM, is compared to a hyper-reduction strategy based on  the empirical interpolation method  to  reduce also the nonlinear fluxes. 
We evaluate the two approaches on a variety of benchmarks involving a generalized variant of the BBM-KdV model with a variable bottom, and a
one-dimensional enhanced  weakly dispersive shallow water system. The results show the potential of both approaches in terms of  cost reduction,
with a clear advantage for the  pdROM in terms of robustness, and for the  EIMROM in terms of cost reduction.
}

\section{Introduction}
Water waves equations can be modeled with various strategies \cite{lannes2020modeling}. Different models provide reasonable approximations in different contexts, according to which type of solution we are interested in. The most complete models consider Euler or Navier--Stokes equations and describe the motion of the water in each point in the 3 dimensional space. A first approximation level is given by averaging the water speed on the vertical direction, reducing by one the dimensions and allowing to display only the water surface as a function of the horizontal variables. The well--known shallow water equations, Boussinesq equations, Green-Naghdi equations, and other models fall in this category. In this work we are interested in studying dispersive models, which are able of reproducing dispersion phenomena. In contrast with other models, as shallow water equations, in these dispersive models shocks are avoided and instead the dispersion fragments the waves into different waves with different wavelength before they break.
We will consider two dispersive models which approximate water waves at different levels: the Benjamin--Bona--Mahony (BBM) equation \cite{benjamin1972model}, in a more general form that it is linked to the Korteweg--De Vries (KdV) equation \cite{dingemans1997water,lannes2020modeling}, and an enhanced Boussinesq (EB) system of equations \cite{madsen1992newForm,ricchiuto2014upwind}.

The peculiarity of both models is the combination of hyperbolic systems and dispersive terms.
Numerical solutions of such problems are often obtained using explicit solver for the hyperbolic part, while the dispersive terms are obtained first solving an elliptic problem, which requires greater computational costs, and then using this solution inside the hyperbolic system \cite{ricchiuto2014upwind,filippini2016flexible}. In this work, we will try to reduce the computational costs of such problems, focusing on the compression of the elliptic operators, which are responsible of the largest part of the computational time of such methods. In case of parametric problems, where a fast response is necessary or we have a multi-query task, the reduction could lead to strong advantages for the computational costs, without degrading the quality of the solution.

Hyperbolic problems are known for developing shocks and advection dominated solutions. Classically, they are badly reducible as the Kolmogorov $n$--width decays very slowly for such problems \cite{RB_freezing,taddei2020registration,Cagniart2019,torlo2020model,mojgani17aleRB,sPOD,peherstorfer18adaptiveBases}. Nevertheless, for dispersive problems there is no shock formations and the wave traveling is often transformed into a oscillation of the whole domain. Hence, also the advection character of the solution is less pronounced. Hence, there is room to attempt a reduction with classical model order reduction (MOR) algorithms coming from the parabolic and elliptic community. In particular, we will use the proper orthogonal decomposition (POD) \cite{POD} to reduce the solution manifold and then we will apply a Galerkin projection to obtain a reduced problem. A further step will consist of interpolating the nonlinear fluxes with the empirical interpolation method (EIM) \cite{barrault04}. This will provide a second level of reduction that allows more reduction in the computational costs. Few works have already performed a model order reduction directly on shallow water equations, \textit{inter alia} \cite{strazzullo2018model,strazzullo2021pod,cstefuanescu2013pod,cstefuanescu2014comparison}. Up to our knowledge, there are no other works trying to reduce the computational costs of dispersive wave equations by the means of model order reduction techniques.

The paper is structured as follows.
In \cref{sec:BBMKdV} we introduce the BBM--KdV model, some energy properties and its classical discretization, i.e., the full order model (FOM), and we introduce few tests that we will study along the paper.  In \cref{sec:projection} we develop the reduction algorithms and the hyper-reduction steps. In \cref{sec:simulations} we test the presented algorithms on all the presented tests, showing the power and limits of these reduction techniques. 
In \cref{sec:EB} we introduce an enhanced Boussinesq system of equations, its FOM discretization and some reference tests that we will use. We introduce then its ROM and its hyper-reduction algorithm. Then we apply the reduced algorithms for the presented tests in \cref{sec:simulationsEB}, where we compare the two reduced methods and the FOM one.
In \cref{sec:conclusion} we highlight the major steps of this work and we suggest possible developments and extensions.

\section{Generalities: dispersive wave models and model reduction}\label{sec:general_model}
We discuss the general form of the dispersive wave model used in this work, which also allows to explain the main underlying  idea. 
Weakly dispersive wave models can be often seen as perturbations of some hyperbolic partial differential system \cite{lannes2020modeling}, and  can be written in general form 
\begin{equation}\label{eq:generalModel0}
	\partial_t u + \partial_x F(u) + S(u)  - \mu^2 \mathcal{X}^{t}\left( \partial_t u + \partial_x G(u)  \right)  +  \mu^2\mathcal{D}^x u   =0
\end{equation}
where $u:\R\to \R^d$ are the unknowns of the system, $F:\R^d \to \R^d$ is a (nonlinear) flux which describes the hyperbolic operator, $S:\R^d \to \R^d$ is a source term which might include bathymetry effects, $\mathcal{X}^t:\R^d \to \R^d$ is a linear operator that contains second derivative terms, while $\mathcal{D}^x:\R^d \to \R^d$ is a linear operator defined by some dispersion terms with third order derivatives.
%
%
%
The first three terms in the above equation define a hyperbolic balance law, as e.g. the shallow water equations with bathymetry,
and   $\mu^2$ is a small parameter multiplying the weakly dispersive regularization \cite{lannes2020modeling}. 
Two specific examples of this general family are used in the paper: the BBM-KdV model \cite{benjamin1972model,dingemans1997water,lannes2020modeling} and an enhanced Boussinesq model \cite{madsen1992newForm,ricchiuto2014upwind}. \\

The main idea of this work is to exploit the smallness of $\mu^2$ to defined a hybrid model, which enhances the hyperbolic equation  with a reduced order approximation of the small terms.
To this end, we follow \cite{filippini2016flexible,Cauquis2022} and recast the system as
\begin{equation}\label{eq:generalModel}
	( 1 - \mu^2 \mathcal{X}^{t}) (\partial_t u + \partial_x F(u)  + S(u)) + \mu^2 \mathcal{X}^{x} u =0 ,
\end{equation}
having set
\begin{equation}\label{eq:generalModel2}
	\mathcal{X}^{x} u:=   \mathcal{X}^{t}(  \partial_x F(u) -  \partial_x G(u) + S(u) ) + \mathcal{D}^x u  .
\end{equation}
%
%
%
The models share the property of being composed of two parts: an elliptic (linear) operator and a hyperbolic (nonlinear) operator. 
The two operators are treated very differently in their discretization, where the hyperbolic part is often discretized in an explicit way, while the elliptic operator is discretized implicitly. 
Moreover, the combination of the two operators allows to obtain stable schemes with little or zero extra numerical viscosity, which is not the case for pure hyperbolic problems.
In particular, we can recast the system introducing an auxiliary variable $\Phi \in \R^d$:
\begin{equation}
	\begin{cases}
		\left(1-\mu^2\mathcal{X}^{t}\right) \Phi    + \mu^2 \mathcal{X}^x u  = 0,\\[5pt]
		\partial_t u + \partial_x F(u)  + S(u) = \Phi.
	\end{cases}
\end{equation}
The operator $(1-\mathcal{X}^t)$ is   elliptic  (and invertible), and by classical formal arguments we can show the smallness of $\Phi$ simply
by noting that  (see also  \cite{lannes2020modeling,book_lanes,LB09})
$$
\Phi = -  \mu^2 \left(1-\mu^2\mathcal{X}^{t}\right)^{-1}\mathcal{X}^x u  = \mathcal{O}(\mu^2).
$$
Our main idea is to devise an approach where only the term $\Phi$ is reduced, possibly benefiting from the combination of the modelling error (controlled by $\mu^2$)
and model reduction error. The motivation for this is discussed in the next section.

\subsection{Time and space discretization of the model}
As in \cite{filippini2016flexible,Cauquis2022} for all models we use a splitting of the elliptic and hyperbolic operators, as in an IMEX time discretization. This  leads to the semi discrete prototype 
\begin{equation}
	\begin{cases}
		\left(1-\mathcal{X}^{t}\right) \Phi^{n+1}   + \mathcal{X}^x u^n = 0,\\[5pt]
		\dfrac{u^{n+1}-u^n}{\Delta t} + \partial_x F(u^n) + S(u^n, u^{n+1}) = \Phi^{n+1}.
	\end{cases}
\end{equation}
Here, we denote with $n$ the time discretization index, and the $\mu^2$ factor has been included in the operators $\mathcal{X}^t$ and $\mathcal{X}^x$ to lighten the expressions. 
As we can see, all the operations are vectorial except the solution of $\Phi^{n+1}$, which requires the solution of a linear system.
After a further spatial discretization the system can be written in the general form
\begin{equation}\label{eq:generalModelFOM}
	\begin{cases}
		\left(\mathbb{M}^\Phi-\mathbb{X}^{t}\right) \Phi^{n+1}   + \mathbb{X}^x u^n = 0,\\[5pt]
		\mathbb{M}^u\dfrac{u^{n+1}-u^n}{\Delta t} + \mathbb F(u^n) + \mathbb S^{ex}(u^n) + \mathbb{S}^{im} u^{n+1} = \Phi^{n+1},
	\end{cases}
\end{equation}
where, with an abuse of notation, we use the same symbols $u\in \R^{d\times \NFE}$ and $\Phi\in \R^{\NFE}$ to refer to the discretized variables, where $\NFE$ is the number of degrees of freedom for each variable. The discretization can be performed by finite differences (FD) or finite element methods (FEM). Specific choices used for our simulations are discussed later in the text.
$\mathbb{M}^u$ and $\mathbb{M}^\Phi$ are mass matrices, $\mathbb{X}^t$ and $\mathbb{X}^x$ are the discretizations of $\mathcal{X}^t$ and $\mathcal{X}^x$ respectively. $\mathbb S^{ex}$ is the discretization of the explicit part of the source function $S$ and $\mathbb S^{im}$ is the matrix that discretizes the linear implicit source terms and $\mathbb{F}(u)$ is a discretization of the (nonlinear) flux $\partial_x F(u)$, which might include some stabilization terms as well.
In the following, we will denote \eqref{eq:generalModelFOM} as full order model (FOM). The time discretization in \eqref{eq:generalModelFOM} is done through an implicit--explicit (IMEX) scheme, where some of the terms are discretized with implicit Euler ($\Phi^{n+1}$) and others with explicit Euler. In practice, we will will use a IMEX Runge--Kutta (RK) schemes with order matching the spatial discretization orders, but all the discussion can be easily derived from the first order discretization of \eqref{eq:generalModelFOM}. Hence, we will keep discussing it, and only for the specific problems we will go in detail with the IMEX RK discretization. To have stability of the hyperbolic operator, we need to impose some CFL conditions on the timestep of the type $\Delta t \rho(JF(u)) \leq  \texttt{CFL} \Delta x$, where $\rho(JF(u))$ is the spectral radius of the Jacobian of $F$ in $u$ and $\texttt{CFL}$ is a constant smaller than 1.

Typically, the most computationally intensive part is the
solution of the linear systems, 
derived from the elliptic operator. 
While in one dimensional cases the linear systems are tridiagonal for simple discretizations and very efficient algorithms like Thomas algorithm can render very fast the solution of such systems, this is not the case in general.
Even in this very favourable case,  we find that 
the cost of the solution of these linear systems is still the largest of the problems (between $60\%$  and $90\%$).
This is also due to the fact that the computation of the explicit flux and source terms can be efficiently done in parallel. 
The main goal of this work is to investigate hybrid methods in which we reduce the cost of the linear system, while keeping the evaluation of the fully nonlinear terms. 

\subsection{Reduction of the general model}
In this section, we describe how to perform some reduction techniques on the operators presented in the FOM. As described above, there are two types of operations in the previous discretized model. The first ones are the vector based operations constituted of the fluxes, sources and diffusion terms, of derivatives and of sums of all terms and (sparse) matrix vector multiplications. The second one are the matrix based operations which include essentially the solution of linear systems. Even if in some particular cases efficient \textit{ad hoc} algorithms can be used (e.g. Thomas algorithm for the tridiagonal matrices),
we have in mind more challenging applications, e.g. two--dimensional problems or high order methods, where such algorithms cannot be used. In those cases, we need to recast to sparse linear solver methods, which in general allow at best to obtain costs of the order of $\mathcal{O}(\NFE \log(\NFE))$ at each time step. These systems require greater computational costs to be solved, with respect to the vector part, and are the ones that mostly weigh on the FOM algorithm.

Our objective here is to  perform some model reduction to gain control  on the computational costs.
To this end,  we introduce different approximations of the model, via reduced basis  expansions  \cite{RB_book_rozza}. 
We want to find an appropriate  set of  $\NRB$ basis of vectors in $\R^\NFE$,  with $\NRB \ll \NFE$, allowing to   express the approximate solutions as 
%
\begin{equation}\label{eq:rbApporx}
	u^n \approx  V \hat{u}^n, \qquad \hat{u}^n \in \R^\NRB
\end{equation}
with the columns of $V \in \R^{(d\times \NFE ) \times \NRB}$ containing the basis vectors, and $ \hat{u}^n$ the vector of the reduced basis coefficients.

A classical way of obtaining a reduced basis method is by projecting onto the reduced space.
We introduce a test matrix $W \in \R^{(d\times \NFE )\times \NRB}$ and then we project the whole equation \eqref{eq:FOMmodel0}, i.e.,
\begin{subequations}\label{eq:ROM_general}
	\begin{align}
		&W^T\left(\mathbb{M}^\Phi-\mathbb{X}^{t}\right) V \hat{\Phi}^{n+1} + W^T \mathbb{X}^x V \hat u^n = 0,\\
		&W^T\mathbb{M}^u V  \frac{\hat u^{n+1}-\hat u^n}{\Delta t} + W^T\mathbb F(V \hat u^n)  + W^T \mathbb S^{ex}(V \hat u^n) + W^T\mathbb S^{im} V \hat u^{n+1}  = W^T V \hat \Phi^{n+1}.
	\end{align}
\end{subequations}
This open the question on which matrix $W$ we should consider. If we choose $W=V$ we obtain a Galerkin projection, while different choices leads to Petrov--Galerkin projections. There are several choices that one can take to obtain $W$, for example controlling the energy of the system \cite{grimberg2021Stability}. We will comment in each model what we use as test matrix $W$.

It should be noted that with the reduction \eqref{eq:ROM_general}, if $\NRB \ll \NFE$, one can get rid of the costs of the solution of the large systems of equations. If we define, for every matrix $\mathbb A \in \R^{\NFE \times \NFE}$, its reduced version as
\begin{equation}
	\hat{\mathbb{A}}:=W^T\mathbb{A} V \in \R^{\NRB \times \NRB},
\end{equation}
we can reduce the computational costs of \eqref{eq:ROM_general}, obtaining
\begin{subequations}\label{eq:ROM_general_compressed}
	\begin{align}
		&\left(\hat{\mathbb{M}}^\Phi-\hat{\mathbb{X}}^{t}\right)  \hat{\Phi}^{n+1}  + \hat{\mathbb{X}}^x  \hat u^n = 0,\\
		&\hat{\mathbb{M}}^u \frac{\hat u^{n+1}-\hat u^n}{\Delta t} + W^T\mathbb F(V \hat u^n) + W^T \mathbb S^{ex}(V \hat u^n) + \hat{\mathbb S}^{im}\hat u^{n+1}   = \hat{\mathbb{I}} \hat\Phi^{n+1},
	\end{align}
\end{subequations}
where $\mathbb{I}\in \R^{\NFE\times \NFE}$ is the identity matrix. In this model, the costs of the linear systems are independent of $\NFE$, while still some computations (the computations of fluxes and source terms) are still dependent of $\NFE$. Nevertheless, their implementation can be easily parallelized and their computational costs do not impact as much as the solver of the full linear system. Still, the computational cost reduction that one can obtain in this context is limited by these terms.

\subsubsection{Choice of the reduced space}
The underlying technique used in this work is 
proper orthogonal decomposition (POD) \cite{POD,RB_book_rozza} or principal component analysis (PCA). 
It consists in a singular value decomposition (SVD) of the snapshot matrix (a simple juxtaposition of the FOM solutions for different parameters and timesteps) of which we retain the most energetic \textit{eigenvectors}, ordering the singular values. The retained $\NRB$ eigenvectors form 
the  basis $V \in \R^{\NFE\times \NRB}$ for the 
reduced order model (ROM).

The number of retained modes/eigenvectors can be chosen either directly, or according to a tolerance on the discarded energy, i.e., 
\begin{equation}
	\NRB:= \arg\max_{N} \left\lbrace \frac{\sum_{j=1}^N \sigma_i }{\sum_{j=1}^{N_{\max}} \sigma_i } \leq 1- tol_{POD} \right \rbrace,
\end{equation}
where $N_{\max}$ is the minimum between $\NFE$ and the number of snapshots considered. 

The POD (and the linear ROMs) are based on the superposition ansatz, for which the solutions can be reasonably represented by a linear combination of few modes $u(x)\approx \sum_{j=1}^{\NRB} v_j(x) \hat{u}_j$. This is not always true, in particular for hyperbolic problems.
More sophisticated algorithms exist  to obtain nonlinear ROMs \cite{Cagniart2019,torlo2020model,RB_freezing,taddei2020registration,sPOD,mojgani17aleRB}, but
these will not be considered in this work. 

\subsubsection{Hyper-reduction of nonlinear operators}\label{sec:EIM}
To further alleviate the computational cost one can  introduce some linearization of the nonlinear operator. 
There are several algorithms for this. 
The approach selected in this work is the  so-called 
empirical interpolation method (EIM) \cite{barrault04}.  The EIM   is a very simple   approximate-then-project approach. 
It is perhaps not the most advanced technique to handle nonlinearities, however it will serve as a reference in terms of cost. 
Using it with different tolerance will allow to have a 
fair comparison. The hyper-reduced model will be referred to as EIMROM in the following. 
We want to mention that several other techniques have been developed in recent years and more stable algorithms are available \cite{yano2019discontinuous,peherstorfer2020stability,drmac2016new,chen2021eim,farhat2015structure}. Still, the focus of this work is not on the hyper-reduction, hence we will stick to a simple EIM algorithm.

Given a nonlinear flux, the goal of the EIM procedure is to approximate it through an interpolation into few \textit{magic points} $\lbrace z_j \rbrace_{j=1}^{\NEIM}$,
and some basis functions $\lbrace\psi_j \rbrace_{j=1}^{\NEIM} \subset \R^\NFE$.  We denote by $Z\in \R^{\NEIM}$ the vector of the magic points, and by  $\Psi_{ij}:=\psi_j(x_i)\in \R^{\NFE \times \NEIM}$ the  matrix of the reduced basis functions evaluated in the magic points, i.e., the slice of $V$ in the magic points. Given   a nonlinear
flux $\varphi\lbrace\boldsymbol{u}\rbrace:(\R^{\NFE})^d \to \R^\NFE$, 
we approximate it as 
\begin{align}
	&\varphi^{EIM}\lbrace\boldsymbol{u}\rbrace (x):= \sum_{j=1}^{\NEIM} \psi_j(x) \varphi\lbrace\boldsymbol{u}\rbrace(z_j),
\end{align}
which, after the projection, reads
\begin{align}
	\hat \varphi^{EIM} (\boldsymbol{u}) = W^T &\varphi^{EIM}(\boldsymbol{u}) = \sum_{j=1}^{\NEIM} W^T \psi_j(x) \varphi\lbrace\boldsymbol{u}\rbrace(z_j) = W^T \Psi \varphi\lbrace\boldsymbol{u}\rbrace(Z).
\end{align}
The reduction lies in the computation of the flux $\varphi\lbrace\boldsymbol{u}\rbrace$ only in few magic points $\lbrace z_j\rbrace_{j=1}^\NEIM$, while we can pre-compute and store the matrix $W^T\Psi \in \R^{\NRB \times \NEIM}$ in the offline phase.

The choice of the points and of the number of the EIM basis functions is crucial in order to maintain a good accuracy of the whole algorithm. Typically, they are chosen in a greedy fashion that minimizes the $\infty$-norm of the worst approximated fluxes. All the details can be found in \cite{barrault04,Drohmann2012}. In practice we will stop such algorithm when the error gets lower than a certain tolerance set proportionally with the error of the RB approximation.

\begin{remark}[Implementation details]
	To assemble and solve the systems at each stage we will use a Python implementation based on \texttt{numpy} \cite{numpy} and on \texttt{Numba} \cite{lam2015numba}. These packages allow to perform a multithread computation of the vector assembly (all the nonlinear terms), while the solution of the sparse FOM systems will be performed with Thomas algorithm, since the matrices are tridiagonal. This leads to have the largest cost in the solution of the system $\mathcal{O}(\NFE)$, while the nonlinear operators can be assembled at a lower cost, depending on the architecture of the machine. We will compare also the sparse solver of \texttt{scipy} to emulate the case of a general sparse solver without the \textit{a priori} knowledge of the tridiagonality. In that version of the algorithm, the whole computational time, except the fluxes, will be affected by a non precompiled code. For ROM systems, which are not sparse, we will use the \texttt{solve} function of the linear algebra package of \texttt{numpy}. The computations all along this work will be performed on an 8 cores Intel(R) Core(TM) i7-10875H CPU @ 2.30GHz. 
\end{remark}

\section{BBM-KdV model}\label{sec:BBMKdV}
\begin{figure}
	\centering
	\subfigure[Dimensional variables]{	
		\begin{tikzpicture}[node distance=2cm] 
			
			\coordinate (y) at (0,2.3);
			\coordinate (x) at (6.,0);
			\draw[<->] (y) node[anchor=east] {} -- (0,0)node[anchor=east]{} --  (x);
			
			\draw[domain=0.:6.,color=brown,samples=100]
			plot (\x, {0.13 - 0.7 *exp(-(6.-\x-1.)^2) - 0.3 *exp(-(6.-\x-4.)^2) + 0.2 *exp(-(6.-\x-3.)^2)}) ;
			\draw[domain=0.:6.,color=blue,samples=100]
			plot (\x, {1. + 0.6 *exp(-(\x-1.)^2) + 0.45 *exp(-(\x-4.1)^2) - 0.3 *exp(-(\x-2.4)^2)}) ;
			\draw[domain=0.:6.,color=blue,samples=100, dotted]
			plot (\x, {1.}) ;
			
			\draw[stealth-stealth] (1.2,1.) -- (1.2,1.5) node[midway,anchor=east] {$\eta(x,t)$};
			\draw[stealth-stealth] (1.4,-0.08) -- (1.4,1.4) node[midway,anchor=west] {$h(x,t)$};
			\draw[stealth-stealth] (5.,-0.57) -- (5.,0.) node[midway,anchor=east] {$b(x)$};
			\draw[stealth-stealth] (0.5,0.) -- (0.5,1.) node[midway, anchor=east] {$h_0$};
			\draw[stealth-stealth] (0.9,1.6) -- (4, 1.6) node[midway, anchor=south] {$L$};
			\draw[stealth-stealth] (6.1, 0.7) -- (6.1, 1.5) node[midway, anchor=west] {$a_0$};
			\draw[stealth-stealth] (6.1, -0.57) -- (6.1, 0.3) node[midway, anchor=west] {$b_0$};
			\draw[stealth-stealth] (3.5, 0.2) -- (3.5, 1.) node[midway, anchor=west] {$\bar{h}(x)$};
			
	\end{tikzpicture}}
	\subfigure[Dimensionless variables]{	
		\begin{tikzpicture}[node distance=2cm] 
			
			\coordinate (y) at (0,2.3);
			\coordinate (x) at (6.,0);
			\draw[<->] (y) node[anchor=east] {} -- (0,0)node[anchor=east]{} --  (x);
			
			\draw[domain=0.:6.,color=brown,samples=100]
			plot (\x, {0.13 - 0.7 *exp(-(6.-\x-1.)^2) - 0.3 *exp(-(6.-\x-4.)^2) + 0.2 *exp(-(6.-\x-3.)^2)}) ;
			\draw[domain=0.:6.,color=blue,samples=100]
			plot (\x, {1. + 0.6 *exp(-(\x-1.)^2) + 0.45 *exp(-(\x-4.1)^2) - 0.3 *exp(-(\x-2.4)^2)}) ;
			\draw[domain=0.:6.,color=blue,samples=100, dotted]
			plot (\x, {1.}) ;
			
			\draw[stealth-stealth] (4.2,1.) -- (4.2,1.45) node[midway,anchor=west] {$a_0\eta(x,t)$};
			\draw[stealth-stealth] (1.4,-0.08) -- (1.4,1.4) node[midway,anchor=west] {$h_0h(x,t)$};
			\draw[stealth-stealth] (5.,-0.57) -- (5.,0.) node[midway,anchor=east] {$b_0b(x)$};
			\draw[stealth-stealth] (0.5,0.) -- (0.5,1.) node[midway, anchor=east] {$h_0$};
			\draw[stealth-stealth] (0.9,1.6) -- (4, 1.6) node[midway, anchor=south] {$L$};
			\draw[stealth-stealth] (6.1, 0.7) -- (6.1, 1.5) node[midway, anchor=west] {$a_0$};
			\draw[stealth-stealth] (6.1, -0.57) -- (6.1, 0.3) node[midway, anchor=west] {$b_0$};
			\draw[stealth-stealth] (3.5, 0.2) -- (3.5, 1.) node[midway, anchor=west] {$h_0\bar{h}(x)$};
			
	\end{tikzpicture}}
	\caption{Water waves: notation\label{fig:sw-notation}}
\end{figure}
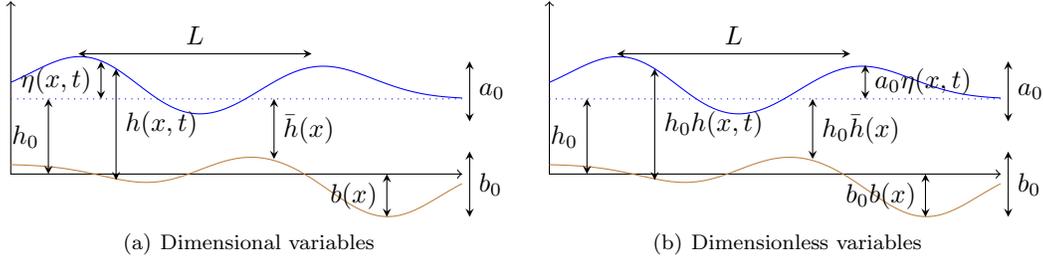

As a first exercise, we study  the  BBM-KdV equations, which provide a scalar one way approximation   of the water wave equations   \cite{benjamin1972model,dingemans1997water,lannes2020modeling}.
It is useful to write the scalar partial differential equation in terms of the dimensionless parameters
$$\mu = \frac{h_0}{L},\quad  \epsilon = \frac{a_0}{h_0}, \quad \beta=\frac{b_0}{h_0},$$
with $h_0$  a characteristic water depth, $L$  the wave length, $a_0$  the  wave amplitude, and with $b_0$ a  characteristic bathymetry amplitude, see \cref{fig:sw-notation} for both dimensional and non-dimensional notations. 
These parameters are classically used  to define the different propagation regimes. More precisely, the dispersive
character of the waves is measured by $\mu$, while $\varepsilon$ is a measure of their nonlinearity.

For $\varepsilon=\mathcal{O}(\mu^2)$, the BBM--KdV equation provides  a 1-way  $\mathcal{O}(\mu^4+\epsilon\mu^2)$ approximation of the water equations.
Its classical form for constant  bathymetry reads

\begin{equation}\label{eq:BBMKdV}
	\left(1-\alpha_p\mu^2\partial_{xx}\right)\partial_t \eta + \partial_x \left( \eta + \frac{3}{2}\varepsilon \frac{\eta^2}{2}  \right) +\mu^2 p \partial_{xxx}\eta =0 ,\qquad p\leq \frac{1}{6},
\end{equation}
having set $\alpha_p=1/6 -p$.
For different values of the parameter $p$ we obtain different approximations. In particular, the KdV equation is obtained for 
$p=\frac{1}{6}$, while  for $p=0$ the equation reduces  to the BBM model. 

For the purpose of its numerical approximation,  we manipulate the model   equation following \cite{filippini2016flexible}, and recast it
as 
\begin{equation}\label{eq:BBMKdVPhi-a}
	\left(1-\alpha_p\mu^2\partial_{xx}\right)\left( \partial_t \eta + \partial_x F(\eta) \right)  
	+ \dfrac{\mu^2}{6}\partial_{xxx}  \eta  
	= - \varepsilon \mu^2\alpha_p\partial_{xxx} \left( \frac{3}{2} \frac{\eta^2}{2} \right), 
\end{equation}
having introduced the short notation for nonlinear flux
\begin{equation}\label{eq:flux}
	F(\eta):=\left( \eta + \frac{3}{2}\varepsilon \frac{\eta^2}{2}  \right).
\end{equation}
Within the same approximation hypotheses, the BBM-KdV equation  is equivalent to
\begin{equation}\label{eq:BBMKdVPhi}
	\left(1-\alpha_p\mu^2\partial_{xx}\right)\left( \partial_t \eta + \partial_x F(\eta) \right)  
	+ \dfrac{\mu^2}{6}\partial_{xxx}  \eta  
	=  0. 
\end{equation}
which we refer to as the modified BBM-KdV equation.
Note that both for \eqref{eq:BBMKdVPhi-a} and \eqref{eq:BBMKdVPhi}, the limit $\mu\rightarrow 0$ corresponds to the  hyperbolic conservation law
\begin{equation}\label{eq:hyp}
	\partial_t \eta + \partial_xF(\eta)=0.
\end{equation}

Equation \eqref{eq:BBMKdVPhi} can be recast into the form of \eqref{eq:generalModel} by setting $u:=\eta$, $\mathcal{X}^{x}:= \frac{\mu^2}{6}\partial_{xxx}$, $\mathcal{X}^t:=\alpha_p \mu^2 \partial_{xx}$ and $S(u):=0$.

\begin{remark}[Dimensional equations]
	One can easily recover the dimensional form of the equation, and more importantly of the variables involved, from the scaling: 
	\begin{equation}\label{eq:dimensionalization}
		\tilde{x} = L x, \quad \tilde{t} = \frac{L}{c_0} t, \quad \tilde{\eta} = \epsilon h_0 \eta, \quad \tilde{b} = \beta h_0 b,\quad \tilde{\bar{h}}=h_0 \bar{h}, \quad \tilde{h}=h_0 h
	\end{equation}
	with $c_0=\sqrt{gh_0}$.
\end{remark}


\subsection{Energy conservation}
A  property of \eqref{eq:BBMKdV}   which we will use in the following is the existence of a certain number of additional derived conservation
laws (see also \cite{ali,ali1}). In particular,   \eqref{eq:BBMKdV}  is endowed with an energy-energy flux pair also verifying a conservation law. 
This can be shown quite classically by   premultiplying  by $\eta$ and manipulating the resulting expression. Setting 
$\Xi(u)=\eta^2/2 + \varepsilon\eta^3/2$, this leads to
\begin{equation}\label{eq:nrg}
	\begin{split}
		&\partial_t\mathcal{E} + \partial_x\mathcal{F}=0,\\
		& \mathcal{E}(\eta):= \frac{\eta^2}{2} + \mu^2\alpha_p \frac{(\partial_x \eta)^2}{2},\\
		& \mathcal{F}(\eta):= \Xi(\eta)
		-\mu^2\alpha_p\eta \partial_{xt} \eta    +\mu^2 p\eta \partial_{xx} \eta   - \mu^2 p \frac{(\partial_x \eta)^2}{2}.
\end{split}\end{equation}

From this relation, we know that the energy $\mathcal{E}(\eta)$ 
is conserved for   \eqref{eq:BBMKdV}.  For the simplified model \eqref{eq:BBMKdVPhi}, this balance is only valid within the $\mathcal{O}(\epsilon\mu^2, \mu^4)$ asymptotic error of the model, 
namely
\begin{equation}\label{eq:nrg-a}
	\partial_t\mathcal{E} + \partial_x\mathcal{F}=\mathcal{O}(\epsilon\mu^2, \mu^4).
\end{equation}
Other balance laws may be derived in a similar spirit, see e.g.   \cite{ali}.  In both cases, it is  convenient  to consider the  scalar product associated to $\mathcal{E}$
\begin{equation}\label{eq:energyAdim}
	\est{\eta}{\rho} := \frac{1}{2}\int_{\Omega} \eta \rho +  \mu^2 \alpha_p \partial_x \eta \partial_x \rho
\end{equation}
which is essentially a modified $H^1$ norm scalar product. 

\subsection{Including the bathymetric effects}
Generalizations  of the KdV equation and of other dispersive one-wave models are discussed in some detail in \cite{israwi2010variable,durufle2012numerical,karczewska2020can}. 
For our purposes,  these generalizations  are useful as they allow to broaden the spectrum of benchmarks.
We use  here the   $\mathcal{O}(\varepsilon \mu^2, \mu^2 \beta, \mu^4)$ approximation
obtained from the variable depth model of    \cite{israwi2010variable},  still referred to here as to the modified BBM-KdV model, and reading 
\begin{equation}\label{eq:bathymetryModel}
	(1-\alpha\partial_{xx})(\partial_t \eta+ (\gamma+\delta \eta) \partial_x \eta) + \omega \partial_{xxx}\eta + \nu \eta=0,
\end{equation}
where  the coefficient depend now on the bathymetry profile $b(x)$ as (in  dimensionless form) 
\begin{equation}\label{eq:parametersKdV}
	\begin{split}
		c(x)= \sqrt{ 1 -\beta b(x) },\;\;\;
		\alpha(x) =&\alpha_p  \mu^2 ,\;\;\;
		\omega(x) = \frac{1}{6} \mu^2 c(x)^5,\;\;\;\\
		\gamma(x) =c(x),\;\;\;
		\delta(x) =& \frac{3}{2} \varepsilon,\;\;\;
		\nu(x) = \frac{3}{2} c'(x).
		\end{split}
	\end{equation}
	The equation can be rewritten up to an $\mathcal{O}(\beta \mu^2)$ as 
	\begin{equation}\label{eq:bathymetryModel2}
		(1-\alpha\partial_{xx})(\partial_t \eta+ (\gamma+\delta \eta) \partial_x \eta+ \nu \eta) + \omega \partial_{xxx}\eta =0,
	\end{equation}
	to better match \eqref{eq:generalModel}. We can indeed define $S(u)=\nu \eta$ to the previous formulation.
	
	As before, a conservation law for the energy  \eqref{eq:nrg-a} can only be shown within   an $\mathcal{O}(\beta, \epsilon\mu^2, \mu^4)$.
	
	\subsection{Full order model discretization}
	We discretize \eqref{eq:bathymetryModel}  using a finite difference  (FD) approach. The resulting discrete model will be denoted by full order model (FOM).
	Classically, we  consider a  discretization of the spatial domain $\Omega = \cup_{i=1}^{\NFE} [x_i,x_{i+1}]$ into  $\NFE$ equi-spaced intervals, and define $\Delta x:= x_{i+1}-x_i$. 
	The temporal domain $[0,T]$ is subdivided into  $K$ time steps $[0,T]=\cup_{k=1}^K [t^{k-1},t^k],$.
	We  denote by  $u_i^k$ the numerical approximation of a quantity $u$ in a point $x_i$ at time $t^k$, and, with an abuse of notation, we use the same symbol for the continuous quantity $u:\R^+\times \R \to \R$, and for the discrete one: $u\in \R^{\NFE\times K}$.
	For convenience, we use the following notation for the classical centered finite difference   operators  $\deriv$, $\deriv_2$ and $\deriv_3:\R^\NFE \to \R^\NFE$:
	\begin{equation}
		(\deriv u)_i = \frac{u_{i+1}-u_{i-1}}{2\Delta x}, \quad (\deriv_2 u)_i = \frac{u_{i+1}-2u_i+u_{i-1}}{\Delta x^2},\quad (\deriv_3 u)_i = \frac{u_{i+2}-2u_{i+1}+2u_{i-1}-u_{i-2}}{2\Delta x^3}.
	\end{equation}
	We start by defining the  IMEX first order Euler discrete  approximation
	\begin{equation}\label{eq:FOMmodel0}
		\begin{cases}
			(\iden- \alpha \deriv_2)\Phi^{n+1}+ \omega \deriv_3 \eta^n =0,\\
			\frac{\eta^{n+1}-\eta^n}{\Delta t} + (\gamma + \delta \eta^n)\deriv\eta^n  + \nu \eta^n +\diff(\eta^n,\lambda)  = \Phi^{n+1}.
		\end{cases}
	\end{equation}
	Here, the vector-by-vector multiplication is meant component-wise in $\eta^n \deriv\eta^n$.
	The term $\diff(\eta^n,\lambda)$ is a high order numerical diffusion term, which depends also on the spectral radius $\lambda$, where $\lambda_i=\gamma +\delta |\eta_i|$. 
	In this work we have chosen to use to this purpose 
	an operator originated from the continuous interior penalty (CIP) approach \cite{burman2004edge} which can be written as 
	\begin{equation}\label{eq:CIP}
		\begin{split}
			\diff(v,\lambda)_j:=d \Delta x^3  \left\{ \lambda_{j+1}(\deriv_2 v)_{j+1} - 2  \lambda_{j}(\deriv_2 v)_{j} + \lambda_{j-1}(\deriv_2 v)_{j-1}
			\right\}
		\end{split}
	\end{equation}
	where $d\in \R^+$ is a  stabilization parameter, which in practice we have  set to  $1$  in all numerical simulations.
	Also \eqref{eq:FOMmodel0} can be reshaped into \eqref{eq:generalModelFOM} simply defining $\mathbb{M}^\Phi = \mathbb{M}^u := \iden$, $\mathbb{X}^t := \alpha \deriv_2$, $\mathbb{X}^x := \omega \deriv_3$, $\mathbb{F}(\eta):= (\gamma + \delta \eta)\deriv\eta + \diff(\eta,\lambda)$ and $\mathbb{S}^{ex}(\eta):= \nu \eta$.

	Starting from \eqref{eq:FOMmodel0}, we  build a  second order SSPRK(2,2) IMEX approximation. The time-step is computed from the CFL condition:  $\Delta t\leq \CFL \Delta x/(\max_i \lambda_i^n )$, with $\CFL <1$.\\
	%

	\subsection{Benchmarks for the  modified BBM-KdV equation}
	In this section we describe some problems that we will use to test the reduced algorithms.
	
	\begin{figure}
		\begin{center}
			\subfigure[Propagation of a periodic monochromatic wave \label{fig:FOMfinalTest1}]{
				\includegraphics[width=0.48\textwidth]{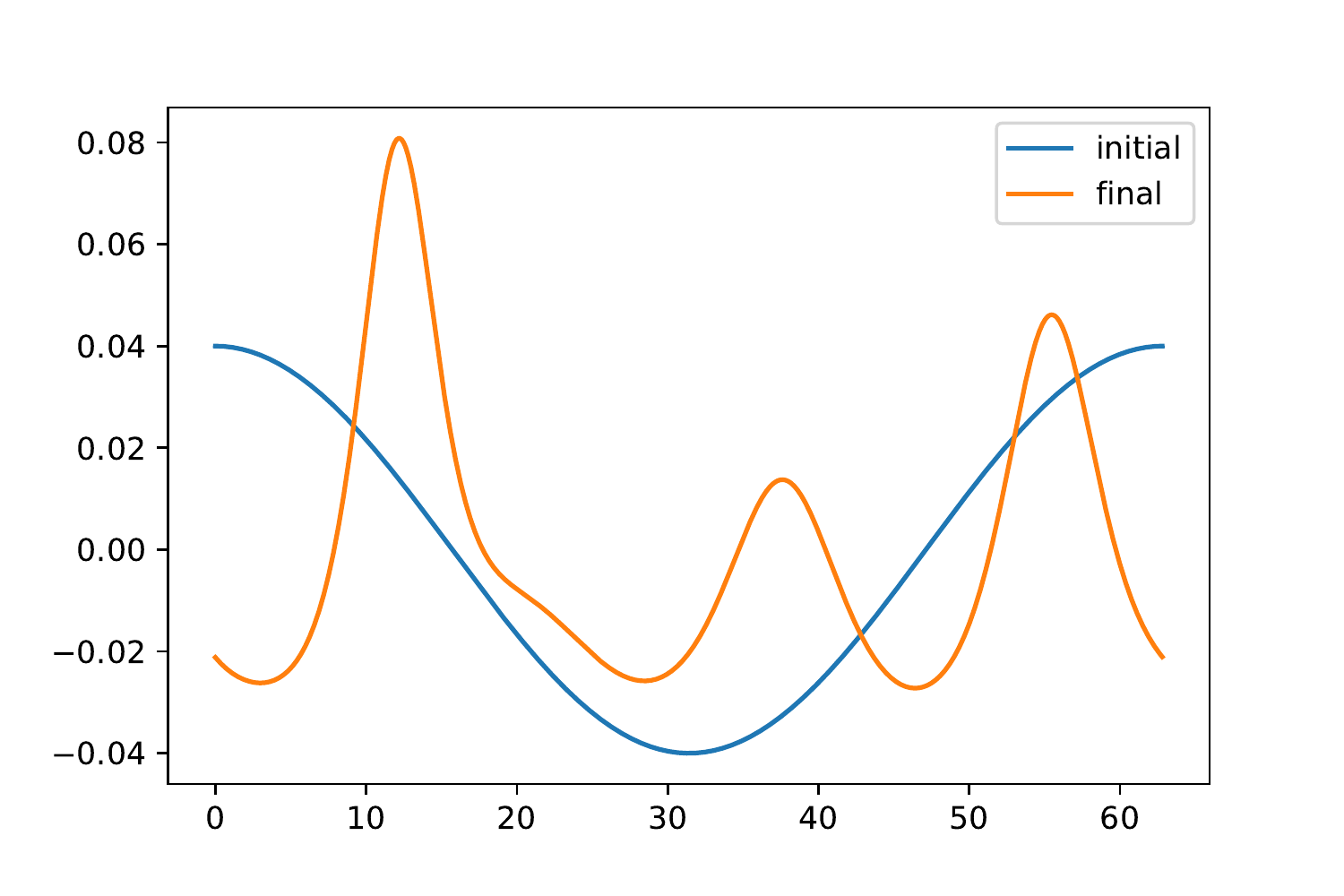}}
			\subfigure[Undular bore propagation\label{fig:FOMfinalTest2}]{
				\includegraphics[width=0.48\textwidth]{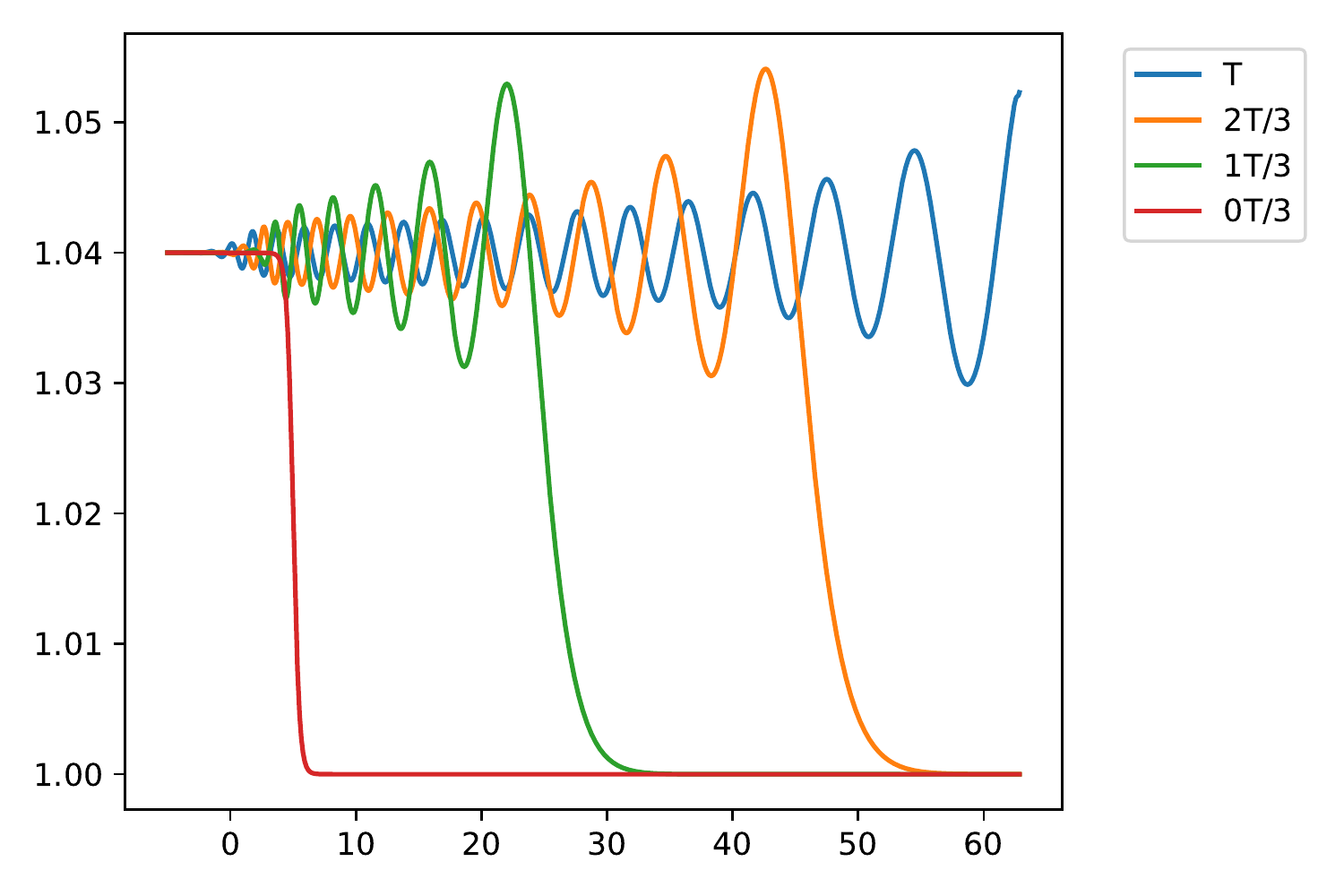}} \\
			\subfigure[Solitary waves interacting with a submerged bar until $T/2$ \label{fig:FOMfinalTest31}]{
				\includegraphics[width=0.48\textwidth]{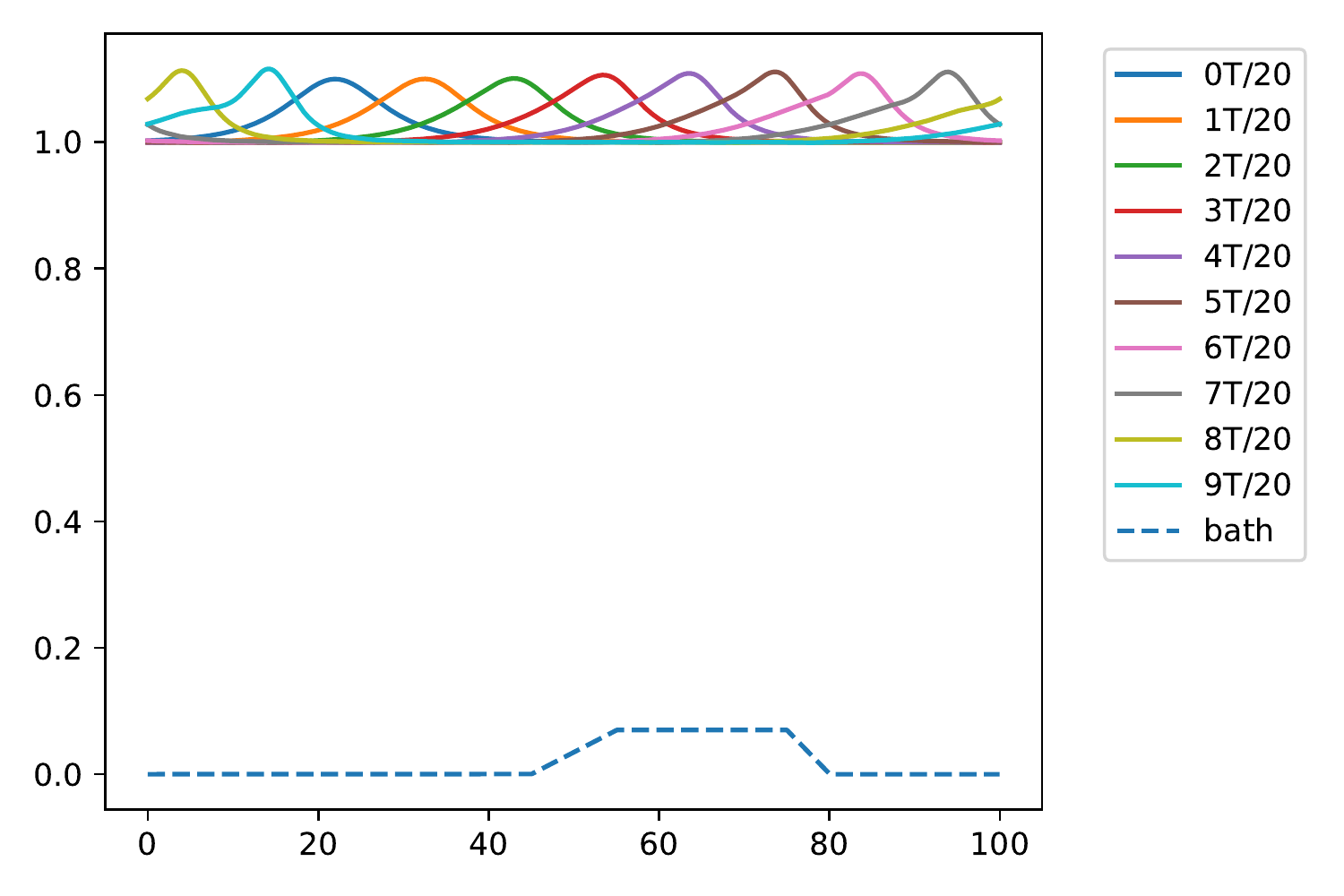}}\hfill
			\subfigure[Solitary waves interacting with a submerged bar from $T/2$ to $T$\label{fig:FOMfinalTest32}]{
				\includegraphics[width=0.48\textwidth]{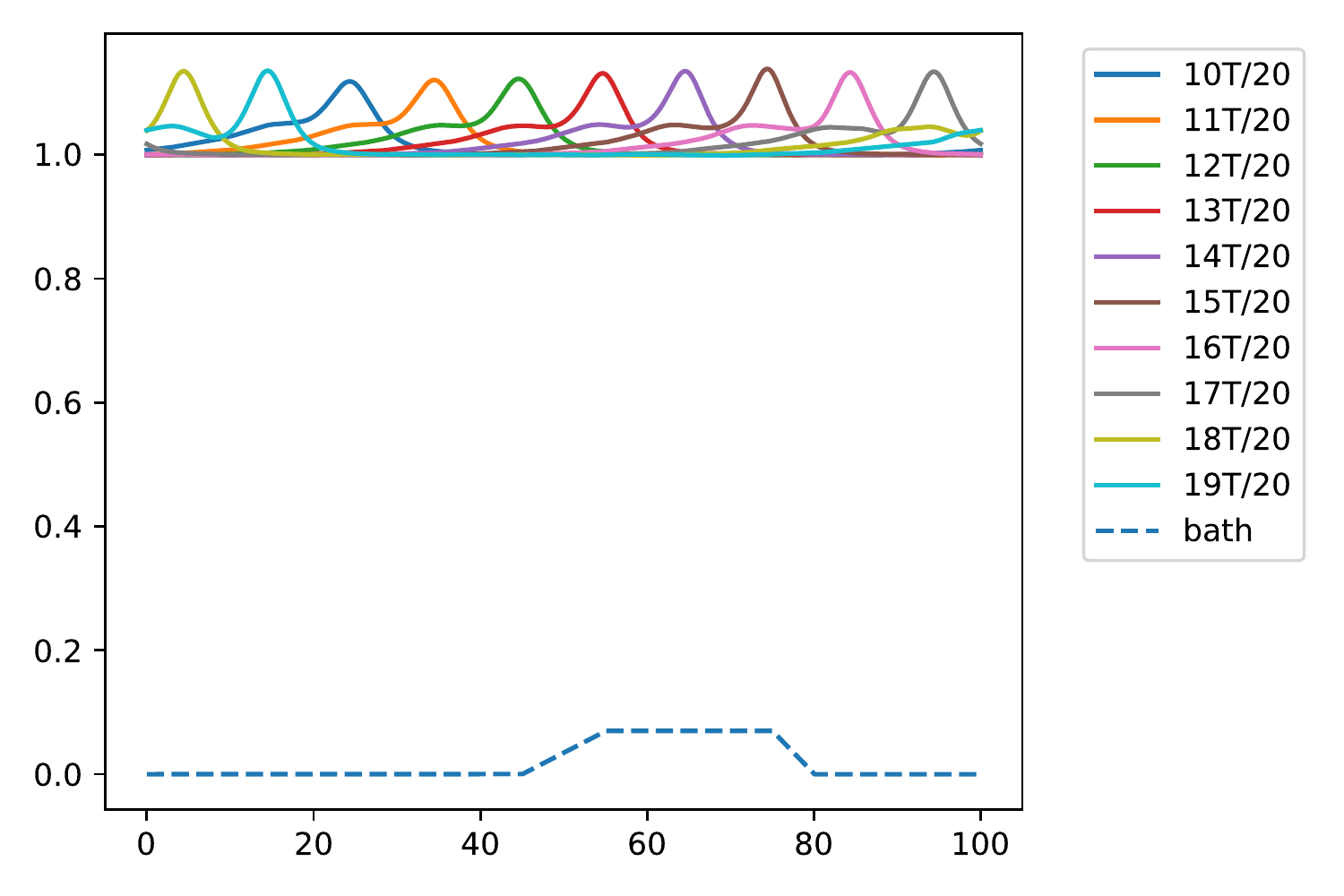}}
			\caption{Initial and final solution for tests.}
		\end{center}
	\end{figure}
	\subsubsection{Propagation of a periodic monochromatic wave}\label{test1}
	The first one starts with a long wave on the whole domain, i.e., a cosine on a periodic domain. This cosine deforms and splits into different waves that travel at different speeds.
	We use in \eqref{eq:BBMKdVPhi} with the dimensionalized coefficients passed through \eqref{eq:dimensionalization}: $h_0=1$, $g=9.81$, $a_0=0.04$, $u_0=a_0 \cos(x /10)$ on the domain $[0,20\pi]$ and final time is $T=200$. The FOM will use $\NFE = 2000$ degrees of freedom and time step will be set with a $\CFL$ number of 0.2. In \cref{fig:FOMfinalTest1} the initial and the final solution are depicted. 
	
	\subsubsection{Undular bore propagation}\label{sec:undular} 
	This test emulates a dam break which develops into several waves. The initial condition is $$u_0 = a_0 \left(1- \frac{1}{1+e^{-4(x-5)}}\right) $$ and the parameters in \eqref{eq:BBMKdV} are $h_0=1$, $g=9.81$, $a_0=0.04$, domain $[0,20\pi]$, $\CFL = 0.1$ and final time $T=20$. Initial and final simulations are shown in \cref{fig:FOMfinalTest2}.
	
	\subsubsection{Solitary waves interacting with a submerged bar}\label{sec:benchSolitaryTrapezoid}
	This test consists of a soliton evolving on a ramp bathymetry defined by 
	$$
	b(x) = 0.07 \begin{cases}
		\frac{x-45}{10} & \text{if }45<x\leq 55,\\
		1 & \text{if } 55<x\leq 75, \\
		-\frac{x-80}{5} &\text{if }75<x\leq 80,\\
		0 & \text{else.}
	\end{cases}
	$$
	The initial condition is defined by $u_0= a_0 \cosh^{-1}(\frac{x-22}{5})$, final time is $T=60$, the parameters are $h_0=1$, $g=9.81$, $a_0=0.1$, $\CFL = 0.1$  and the domain is $[0,100]$. In \cref{fig:FOMfinalTest31,fig:FOMfinalTest32} the evolution of the solution in time is shown.
	
	\begin{table}
		\centering
		\caption{Estimate percentage of computational time for solving the linear systems for the three benchmark problems: comparison of the implementation of Thomas algorithm all in \texttt{numba} framework and sparse solver by \texttt{scipy} with fluxes computed in \texttt{numba}  \label{tab:kdv_system}}
		\begin{tabular}{|c||c|c|c||c|c|c|}
			\hline 
			& \multicolumn{3}{|c||}{\texttt{numba} and Thomas} &\multicolumn{3}{|c|}{\texttt{scipy} and \texttt{numba}}  \\ \hline
			$N_h$ & \ref{test1} & \ref{sec:undular} & \ref{sec:benchSolitaryTrapezoid} & \ref{test1} & \ref{sec:undular} & \ref{sec:benchSolitaryTrapezoid} \\ \hline 
			2000 &86.4\% & 63.4\% & 85.2\%& 86.6\% & 84.4\% & 86.6\% \\
			4000 &84.1\% & 62.4\% & 83.9\%&86.2\% & 83.5\% & 88.8\%  \\
			8000 &86.1\%&61.5\% & 85.0\%& 94.1\% & 85.3\% & 91.0\% \\
			\hline
		\end{tabular}\\ \vspace{2mm}
		\caption{Estimate of ratio of computational time of the linear system over computational time of the nonlinear fluxes for the three benchmark problems: comparison of the implementation of Thomas algorithm all in \texttt{numba} framework and sparse solver by \texttt{scipy} with fluxes computed in \texttt{numba} \label{tab:kdv_flux_system}}
		\begin{tabular}{|c||c|c|c||c|c|c|}
			\hline 
			& \multicolumn{3}{|c||}{\texttt{numba} and Thomas} &\multicolumn{3}{|c|}{\texttt{scipy} and \texttt{numba}}  \\ \hline
			$N_h$ & \ref{test1} & \ref{sec:undular} & \ref{sec:benchSolitaryTrapezoid} & \ref{test1} & \ref{sec:undular} & \ref{sec:benchSolitaryTrapezoid} \\ \hline 
			2000 &11.2 &3.5 & 11.3 &72.4 &52.1 & 71.8\\
			4000 &11.8 &3.7 & 11.6 &73.4 &54.9 & 74.3 \\
			8000 &11.2 &3.7 & 12.8 &78.4 &56.1 & 75.8\\
			\hline
		\end{tabular}
	\end{table}
	Now, we compare the computational times of the solution of the linear systems for the three benchmark problems with respect to other costs in the problem solution. In \cref{tab:kdv_system} we compare the costs of the linear systems with respect to the whole simulation, while in \cref{tab:kdv_flux_system} we see the ratio of the cost of the solution of the linear systems over the cost of the nonlinear fluxes. Let us focus on Thomas algorithm first, we notice that boundary conditions make a huge difference. Indeed, periodic boundary conditions imply the solution of Thomas algorithm three times instead of just one and this increases the costs of the system solver.
	
	On the other side, using a general solver for sparse system, like the one provided by \texttt{scipy}, increases the computational costs of the solver even by a factor of around 10. In particular, the distinction between different boundary conditions is less evident in this case. 
	This algorithm is more interesting in perspective of a more general discretization or in case of multidimensional problem, but less efficient.

	\subsection{Projection based reduction}\label{sec:projection}
	As said before, the computational operation in \eqref{eq:FOMmodel0} that we mainly aim at reducing is the solution of the linear system, while the flux and source terms are easily parallelized. 
	
	As before, we introduce the reduced variable of dimension $\NRB \ll \NFE$ so that we approximate 
	%
	\begin{equation}\label{eq:rbApprox_bbm}
		\eta^n \approx  V \hat{\eta}^n, \qquad \hat{\eta}^n \in \R^\NRB
	\end{equation}
	with  the columns of $V \in \R^{\NFE \times \NRB}$ containing the basis vectors, and $ \hat{\eta}^n$ the vector of the reduced basis coefficients.
	
	Two projection based  model order reduction (MOR)   strategies exist \cite{grimberg2021Stability}, equivalent under certain conditions.
	The first one is the projection with a test matrix $W \in \R^{\NFE \times \NRB}$ and it is obtained by projecting the whole equation \eqref{eq:FOMmodel0}, i.e.,
	\begin{subequations}
		\begin{align}\label{eq:projection_bbm}
			W^T \left(V\hat{\eta}^{n+1}-r(V\hat{\eta}^n,\Phi^n(V\hat{\eta}^n))\right)=0,\\
			r(\eta,\Phi)=\eta  -\Delta t \left( (\gamma+\delta \eta ) \deriv \eta +\diff(\eta,\lambda) -  \Phi \right).
		\end{align}
	\end{subequations}
	This opens the question on which matrix $W$ we should consider. If we choose $W=V$ we obtain a Galerkin projection, while different choices lead to Petrov-Galerkin projections.
	
	An alternative approach, allowing to answer    this question, is   to construct the reduced space using   a   minimization problem. Using in particular the scalar product \eqref{eq:energyAdim}, we can for example define
	\begin{equation}\label{eq:minimization}
		\hat{\eta}^{n+1} = \arg\min _{x} \est{Vx-r(V\hat{\eta}^n,\Phi^{n+1}(V\hat{\eta}^n)}{Vx-r(V\hat{\eta}^n,\Phi^n(V\hat{\eta}^n)}.
	\end{equation} 
	Using  the discrete form of \eqref{eq:energyAdim} (in dimensional variables), we can introduce the scalar product matrix
	$$\Theta:=\Delta x (\iden + \alpha \deriv^T\deriv ),
	$$ 
	where $\iden\in \R^{\NFE \times \NFE}$ is the identity matrix and $\deriv\in \R^{\NFE \times \NFE}$ is the  (first) derivative matrix, and 
	such that at the discrete level
	$$\est{u}{v}= u^T\Theta v.$$ 
	Then, we can write the minimization \eqref{eq:minimization} as the solution of
	\begin{subequations}
		\begin{align}
			2V^T\Theta V \hat{\eta}^{n+1} -2V^T \Theta r(V\hat{\eta}^n,\Phi^{n+1}(V\hat{\eta}^n))=0 ,
		\end{align}
	\end{subequations}
	suggesting  that minimizing the energy in the projected space requires  testing with  
	$W=\Theta^T V$ in \eqref{eq:projection_bbm}.
	
	\subsubsection{Approximation accuracy of the problems}
\begin{figure}
		\begin{center}
			\subfigure[Propagation of a periodic monochromatic wave \label{fig:eigenDecay1}]{\includegraphics[width=0.32\textwidth]{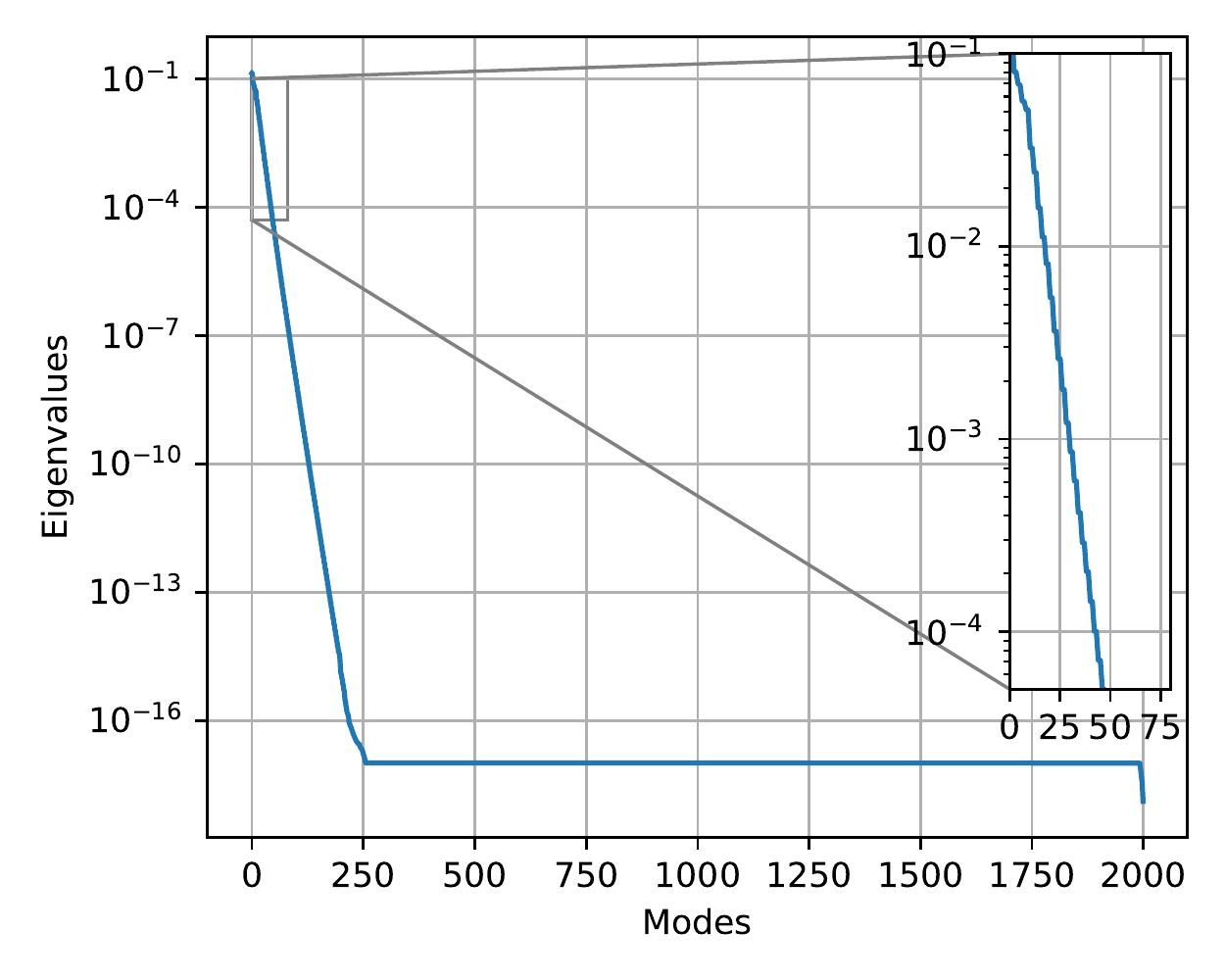}}
			\subfigure[Undular bore propagation \label{fig:eigenDecay2}]{\includegraphics[width=0.32\textwidth]{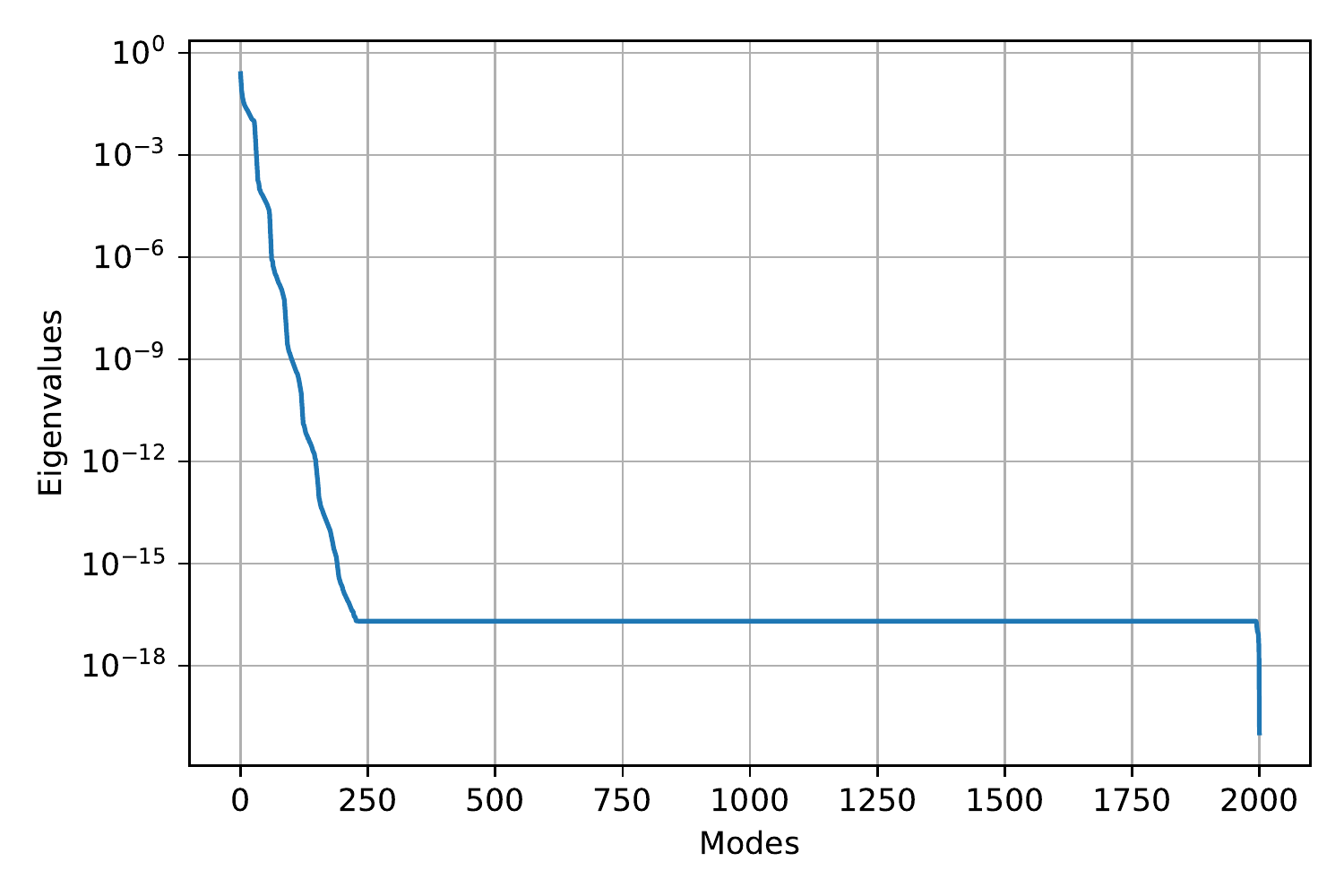}}
			\subfigure[Solitary waves interacting with a submerged bar \label{fig:eigenDecay3}]{\includegraphics[width=0.32\textwidth]{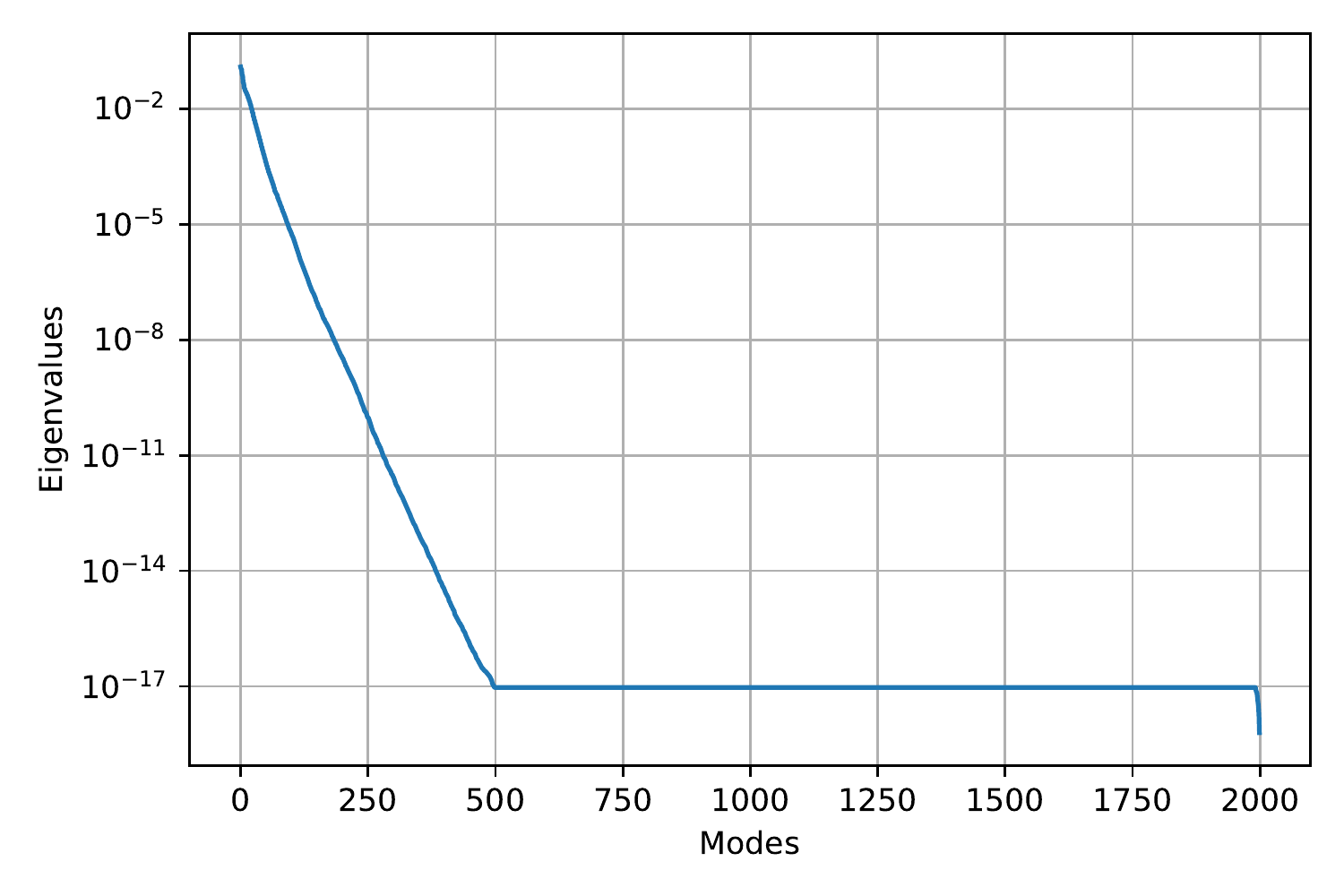}}
			\caption{Singular values decay}
		\end{center}
	\end{figure}

	For the benchmarks involving the modified BBM-KdV equation we have applied the POD
	on 1000 snapshots in time, with fixed parameters.   The FOM dimension is here $\NFE\approx 10^3$, which is a representative order of magnitude. In practice, we have used 2000 points.
	Concerning the potential for reduction, we need to look at the error decay in terms of singular values. 
	For the tests considered here, the decay has been plotted  in \cref{fig:eigenDecay1,fig:eigenDecay2,fig:eigenDecay3}.
	As we can see, the energy drop is not exceedingly fast, still being exponential from the very beginning.
	In particular, for  compact travelling waves   as in the  last benchmark, the advection dominated regime slows down the convergence of the eigenvalues. Nevertheless,
	if we want to keep the error of the order of $10^{-2}$  or $10^{-3}$ the number of eigenmodes needed is in between 50 and 70, which is much below the typical mesh size used in the FOM and can already reduce the computational cost of the linear system (even if full). We remark that the computational reduction is not proportional to the dimension as the original FOM has a sparse linear system to be solved, while the ROM has a full linear system.
	
	\subsubsection{Reduction of the linear operator}
	We consider here the reduction of 
	the linear system which gives $\Phi$ in \eqref{eq:FOMmodel0}. In general this operation can scale  as an $\mathcal{O}(\NFE^2)$, unless particularly optimized
	algorithms are implemented. Indeed, 
	in our FOM simulations we will use Thomas algorithm which only scales as  $\mathcal{O}(\NFE)$. 
	When the linear operator is the only one reduced,  we still need to assemble the discretization of the nonlinear  part of the PDE.
	For this reason we will refer to this approach as to a partial discrete reduced order model or pdROM.\\
	
	The pdROM discrete  evolution equation  can be written as 
	\begin{subequations}
		\begin{align}
			&W^T V \hat{\eta}^{n+1}=W^T V\hat{\eta}^n -W^T \Delta t \left( (\gamma + \delta V\hat{\eta}^n)\deriv V\hat{\eta}^n + \diff ( V\hat \eta^n, \lambda )\right)+\Delta t W^T\Phi^{n+1},\\
			&W^T\Phi^{n+1}= -  \omega W^T(\iden - \alpha \deriv_2)^{-1}  \deriv_3 V(\hat{\eta}^n+ \nu \hat{\eta}^n).
		\end{align}
	\end{subequations}
	
	Also here, we can see the connection with \eqref{eq:ROM_general}, where all the terms can be matched.
	In practice, for $\NRB$ sufficiently small, the  $\NRB\times \NRB$ matrix $A_\NRB:= W^T(\iden - \alpha \deriv_2)^{-1}  \deriv_3 V$ can be pre-computed and stored.
	Denoting with $M_\NRB:=W^T V\in \R^{\NRB\times \NRB}$, the pdROM update can be effectively implemented as 
	\begin{subequations}\label{eq:pdROM_bbm}
		\begin{align}
			&\hat{\Phi}^{n+1}= - \omega A_\NRB^{-1} (\hat{\eta}^n + \nu \hat{\eta}^n),\\
			&\hat{\mathbb{F}}(\eta) = W^T\left( (\gamma + \delta \eta)\deriv \eta+ \diff ( \eta, \lambda ) \right),\\
			&\hat{\eta}^{k+1} = \hat{\eta}^n +\Delta t M_\NRB^{-1} \left( -\hat{\mathbb{F}}(V\hat{\eta}^n) + \hat{\Phi}^{n+1} \right).
		\end{align}
	\end{subequations}
	
	It is clear that, with this approach, the parameters $\gamma, \delta, \omega$ can be easily modified in the online computations. Also $\alpha$ can be decoupled from the reduction process and used to assemble the linear reduced system, but in this context $\alpha$ represent the topography parameters and we do not aim to change it. In the system case we will see how to modify a parameter in the LHS matrix.
	
	\begin{remark}[Modification to nondispersive code]\label{rem:modification_sw}
		It is very important to stress that  the computational results, as well as all conclusions concerning  cost reduction 
		can be obtained keeping a classical solver for nondispersive nonlinear equations, to which we add the reduced dispersive term, i.e.,
		\begin{subequations}
			\begin{align}
				&{\eta}^{k+1}={\eta}^n - \Delta t \left( (\gamma + \delta {\eta}^n)\deriv {\eta}^n+\nu {\eta}^n + \diff (  \eta^n, \lambda )\right)+\Delta t \beta V \hat{\Phi}^{n+1},\\
				&\hat{\Phi}^{n+1}:= - (W^T(\iden - \alpha \deriv_2)V)^{-1} W^T \deriv_3 {\eta}^n,
			\end{align}
		\end{subequations}
		where $B_{\NRB}:=W^T(\iden - \alpha \deriv_2)V\in \R^{\NRB\times \NRB}$ can be stored and inverted at  low computational cost. Nevertheless, working on two different spaces contemporary does not guarantee the same accuracy as in the previously presented algorithm. We have observed that reducing and projecting only the $\Phi$ term leads to much more inaccurate methods. We believe that this is due to spurious modes that the reconstruction of $\Phi$ generates and that are not appropriately damped in the full dimensional space, while they are in the pdROM approach. In the test section we will show that both algorithms have comparable computational costs, but \eqref{eq:pdROM_bbm} provides more stability.
	\end{remark}

	\begin{remark}[EIMROM stability]
		As shown in \cref{sec:EIM}, we can reduce the computational costs of the source and flux terms with a hyper-reduction technique. The extra reduction of the EIM algorithm does not guarantee anymore the energy minimization \eqref{eq:BBMKdV}. This is noticeable when the dimension of the EIM space $\NEIM$ is smaller or close to the dimension of the RB space $\NRB$. In that regime, it is easy to see instabilities and Runge phenomena. In the simulation section, we will observe in which regime we should stay with the dimension $\NEIM$. Other works already shown that the dimension of the EIM space should be larger that the reduced space one to obtain stability of the reduced model \cite{peherstorfer2020stability,zimmermann2018geometric,ghavamian2017pod,argaud2017stabilization,chen2021eim}. Some novel algorithms have been proposed in recent papers \cite{chen2021eim,peherstorfer2020stability} based on different techniques as the over-collocation method and variants of the Gappy-POD. We will not investigate in deep this behavior in this work.
	\end{remark}
	
	\begin{remark}[Fully implicit]
		We want to remark that even for the pdROM simulations, even if the solution is energy stable up to an $\mathcal{O}(\varepsilon\mu^2, \mu^4)$, oscillations may arise due to the under-resolution in the reduced space. This would not lead to blow ups of the solution, but it might have some unphysical behaviors. In the case of the EIM algorithm the oscillations might also lead to instabilities and blow ups when the EIM space does not properly approximate the nonlinear flux manifold.
		In order to alleviate this issue, it might be helpful to consider a fully implicit discretization or, even better, an implicit Euler time discretization. On the other side, this would lead to extra computational costs due to nonlinear solvers that might heavily slow down the reduced formulation, both in the pdROM and EIMROM contexts. In particular, the ROM computational costs with respect of the FOM ones, where the benchmark algorithms are explicit for the hyperbolic equation, might be disadvantageous, at least for one dimensional simulations.
	\end{remark}
	
	\section{Simulations for KdV-BBM}\label{sec:simulations}
	\subsection{Modus operandi}\label{sec:modusOperandi}
	In this section, we will test the pdROM and EIMROM in different regimes, showing their error and the computational costs with respect to the FOM solutions. 
	\begin{enumerate}
		\item We will proceed studying a time-only manifold of solutions, in order to understand which errors our algorithms can reach with different dimensions of reduced basis spaces.
		In this first phase, we will fix the parameters of the problem, we will compute the FOM solution on a mesh with $\NFE=2000$ points and we will collect 1000 snapshots in time $u(t^k)$ and the respective nonlinear fluxes $F(u(t^k))$. These snapshots will form the training set used to choose the reduced basis space and the EIM space, setting either the dimension of the space or the tolerance.
		Then, different simulations of the original problem are run with pdROM and EIMROM, recording their errors and their computational times. We will picture also one significant simulation of the reduced algorithms for a visual comparison for each test. 
		\item The second stage will consider a solution manifold generated by time and parameters. In practice, we will vary $h_0$, on which all the parameters depend in the BBM-KdV problem, see \eqref{eq:parametersKdV}. In particular, we build the training sample using 10 randomly sampled $h_0\in [h_{min},h_{max}]=[0.7 \bar{h}_0, 1.3 \bar{h}_0 ]$, and for each of these parameters, we collect 1000 snapshots at uniform timesteps. We populate the ROM and EIM spaces using these 10.000 snapshots and their fluxes. Then we test the reduced algorithms for $\bar{h}_0$ for a larger final time (for every test we will specify which one), checking again how the reduced procedures perform with different tolerances for EIM and POD. Finally, we test the methods for a parameter $h_0 = 0.63\cdot \bar{h}_0 $ outside the training set and, when sensible, with a different initial condition (for every test we will specify which one), to assess the property of the reduced algorithms in the extrapolation regime.
	\end{enumerate}

	\subsection{Propagation of a periodic monochromatic wave}
	\newcommand{\folderTest}{Test_periodicSoliton}
	\begin{figure}
		\begin{center}
			\subfigure[EIM ROM errors time-only reduction \label{fig:ROMEIMerrPerSol}]{\includegraphics[width=0.32\textwidth]{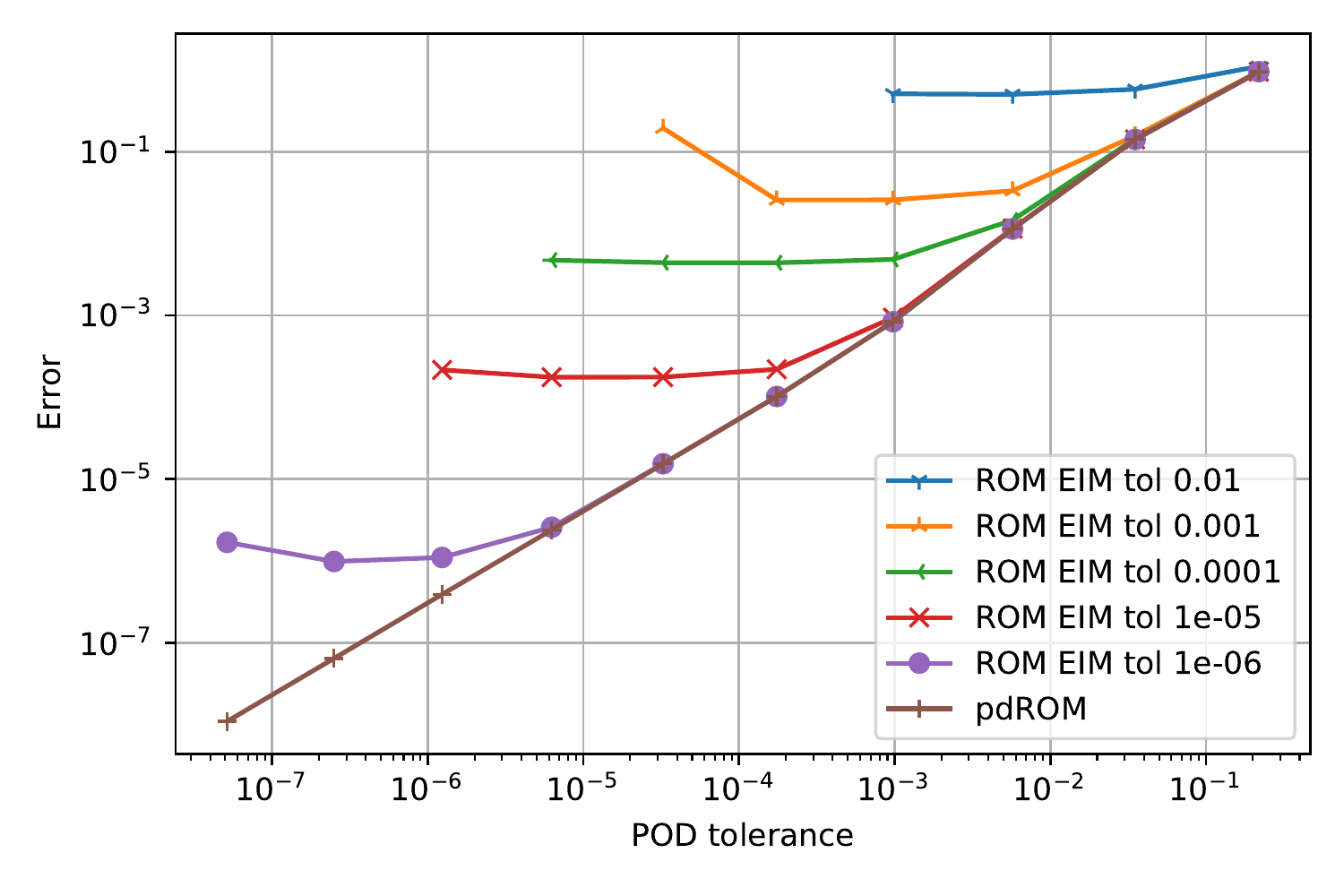}}\;
			\subfigure[ROM errors time-parameter reduction on parameter and $\eta_0$ in the training set \label{fig:ROMEIMparamErrPerSolin}]{\includegraphics[width=0.32\textwidth]{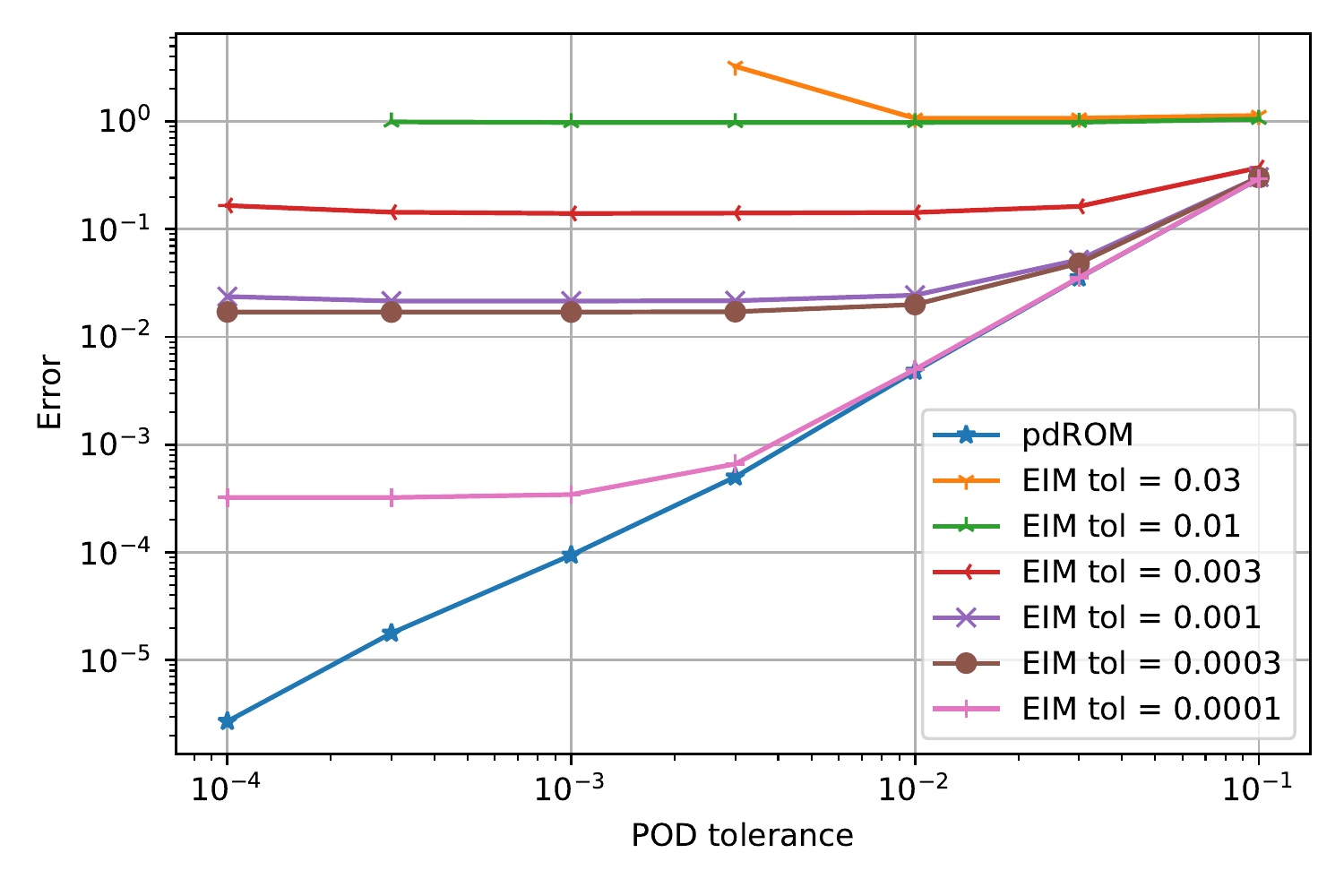}}\;
			\subfigure[ROM errors time-parameter reduction on parameter and $\eta_0$ outside the training set \label{fig:ROMEIMparamErrPerSolout}]{\includegraphics[width=0.32\textwidth]{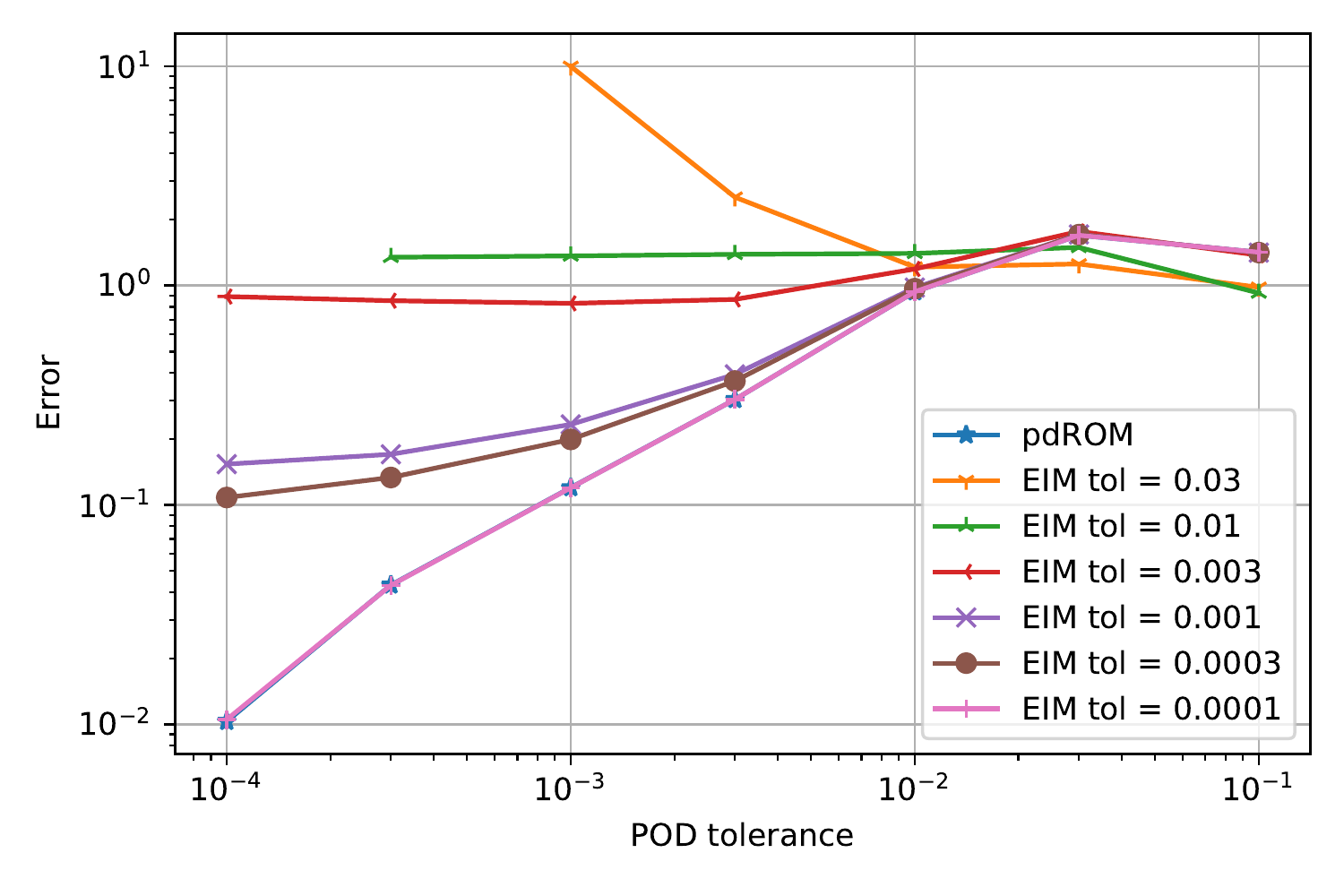}}\;
			\\
			\subfigure[EIM ROM computational times time-only reduction (legend above) \label{fig:ROMEIMtimePerSol}]{\includegraphics[width=0.32\textwidth]{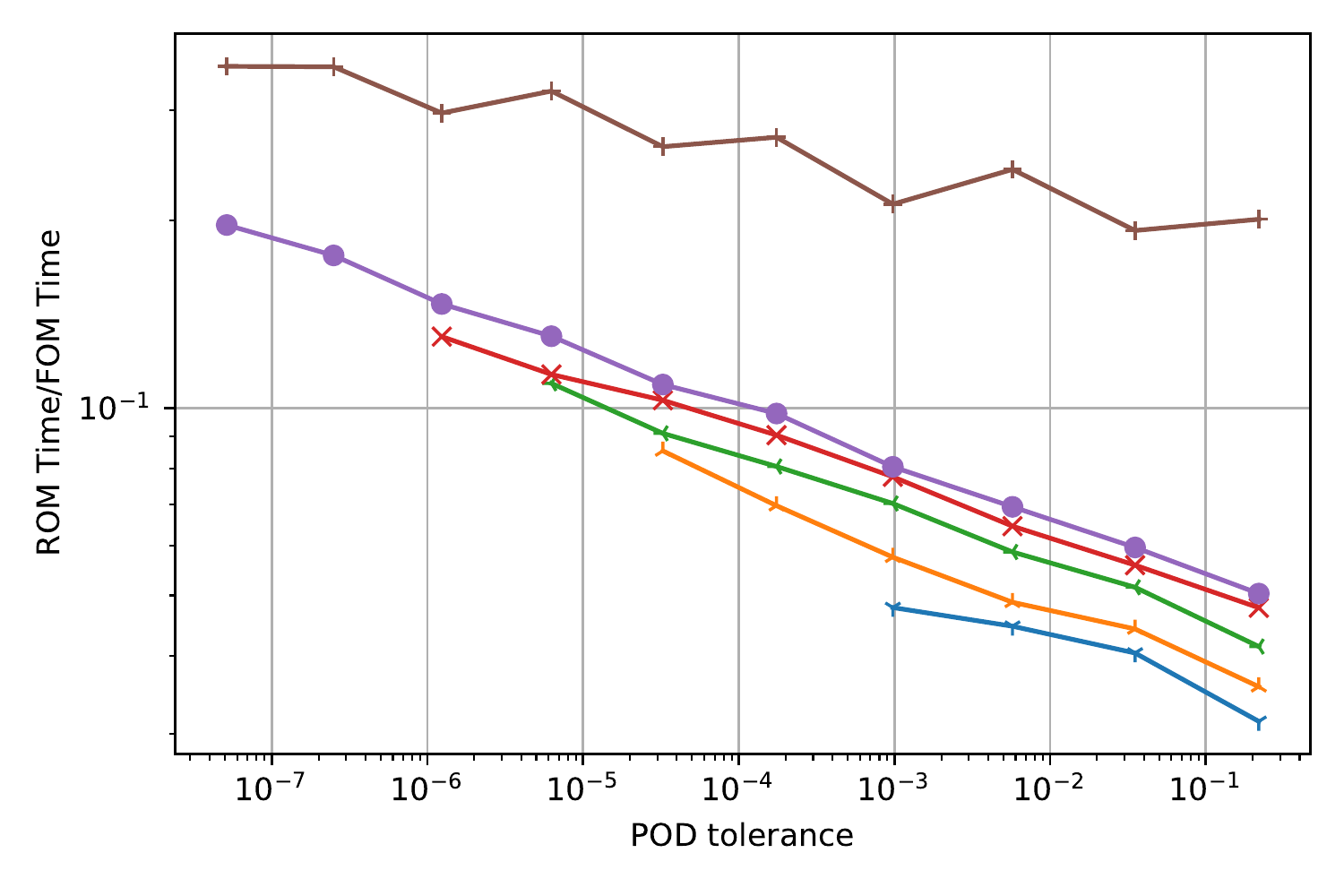}}\;
			\subfigure[ROM computational time time-parameter reduction on parameter and $\eta_0$ in the training set (legend above) \label{fig:ROMEIMparamTimePerSolin}]{\includegraphics[width=0.32\textwidth]{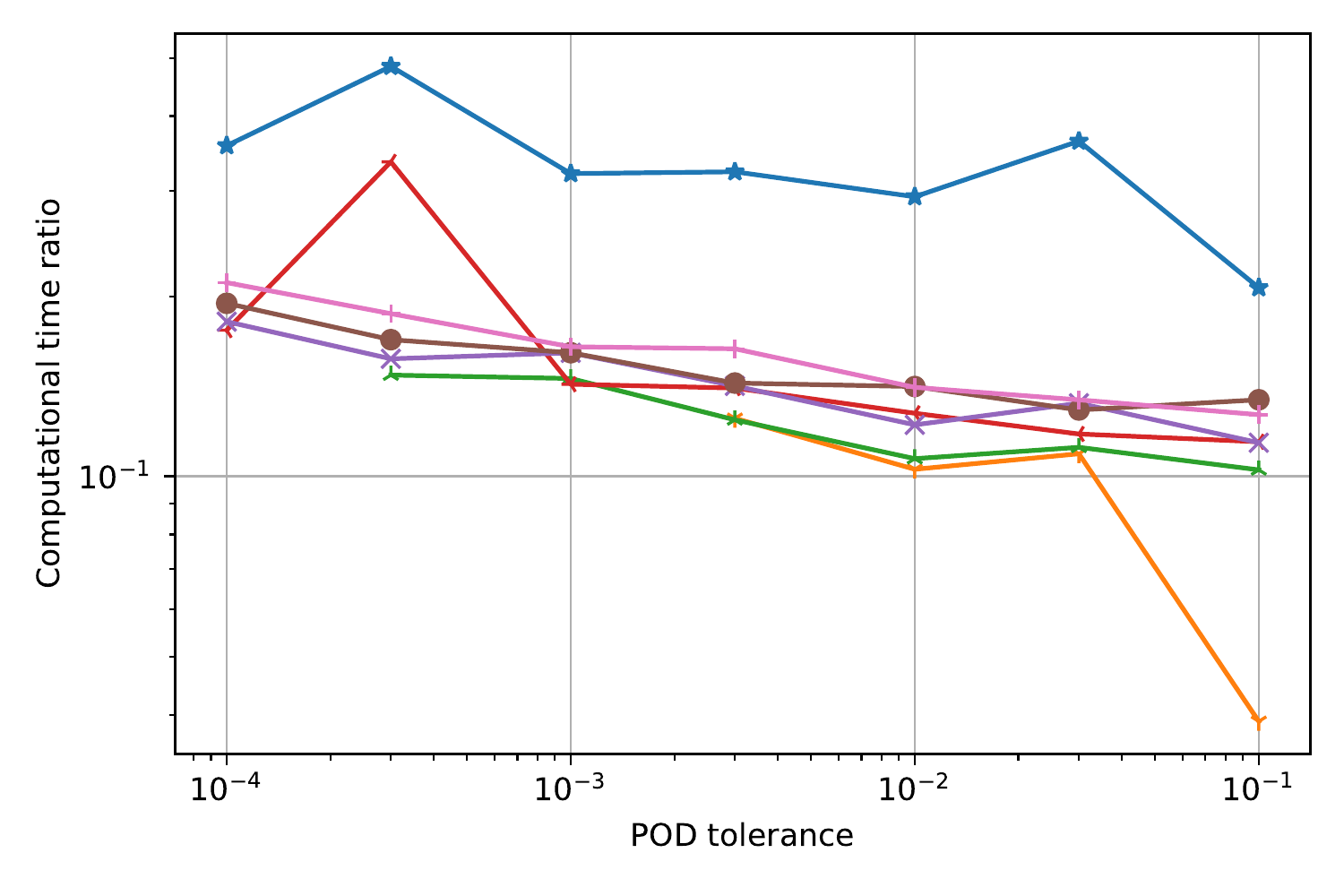}}\;
			\subfigure[ROM computational time time-parameter reduction on parameter and $\eta_0$ outside the training set (legend above) \label{fig:ROMEIMparamTimePerSolout}]{\includegraphics[width=0.32\textwidth]{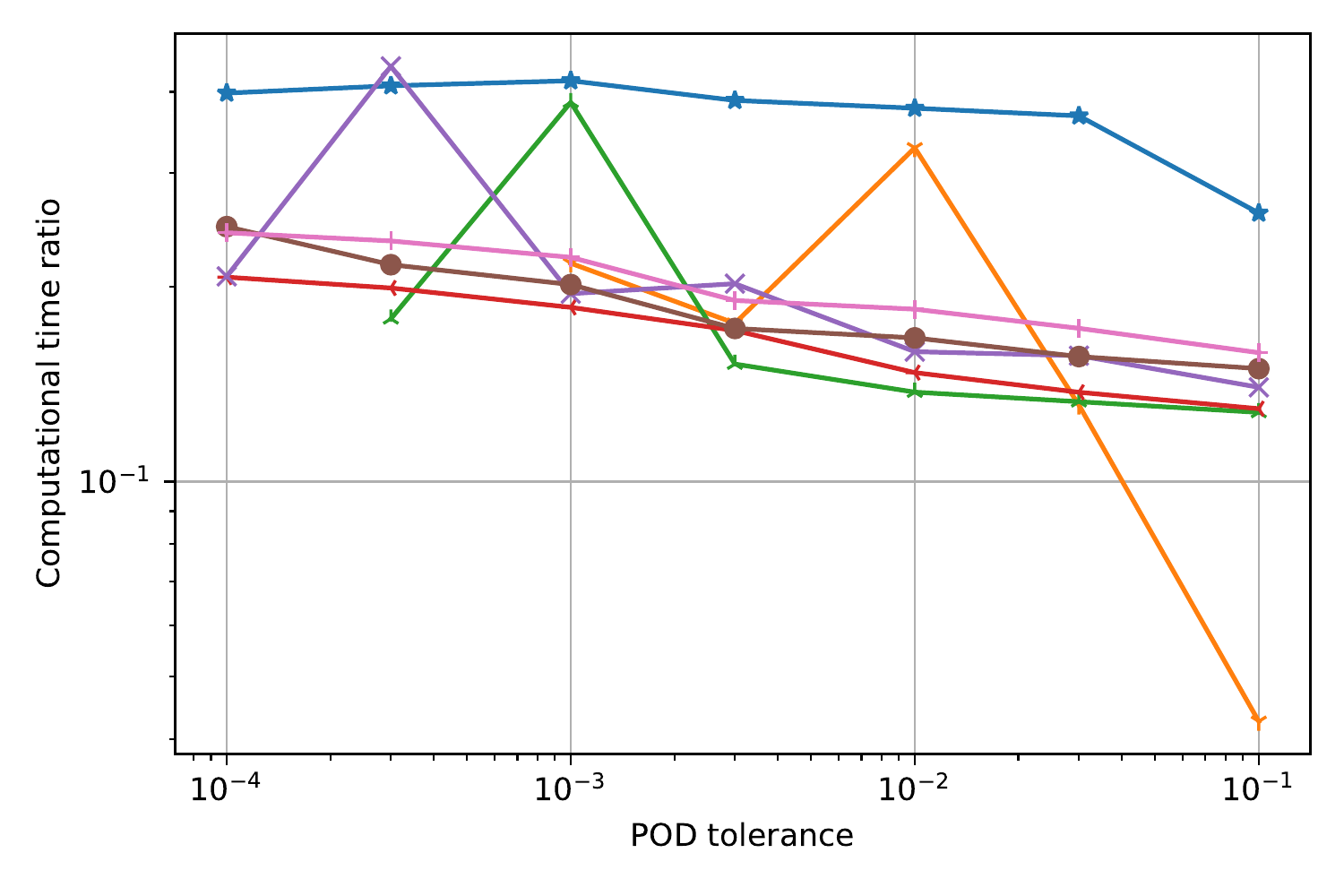}}		\\
			\subfigure[ROM solutions time-only reduction \label{fig:ROMsolPerSol}]{\includegraphics[width=0.32\textwidth]{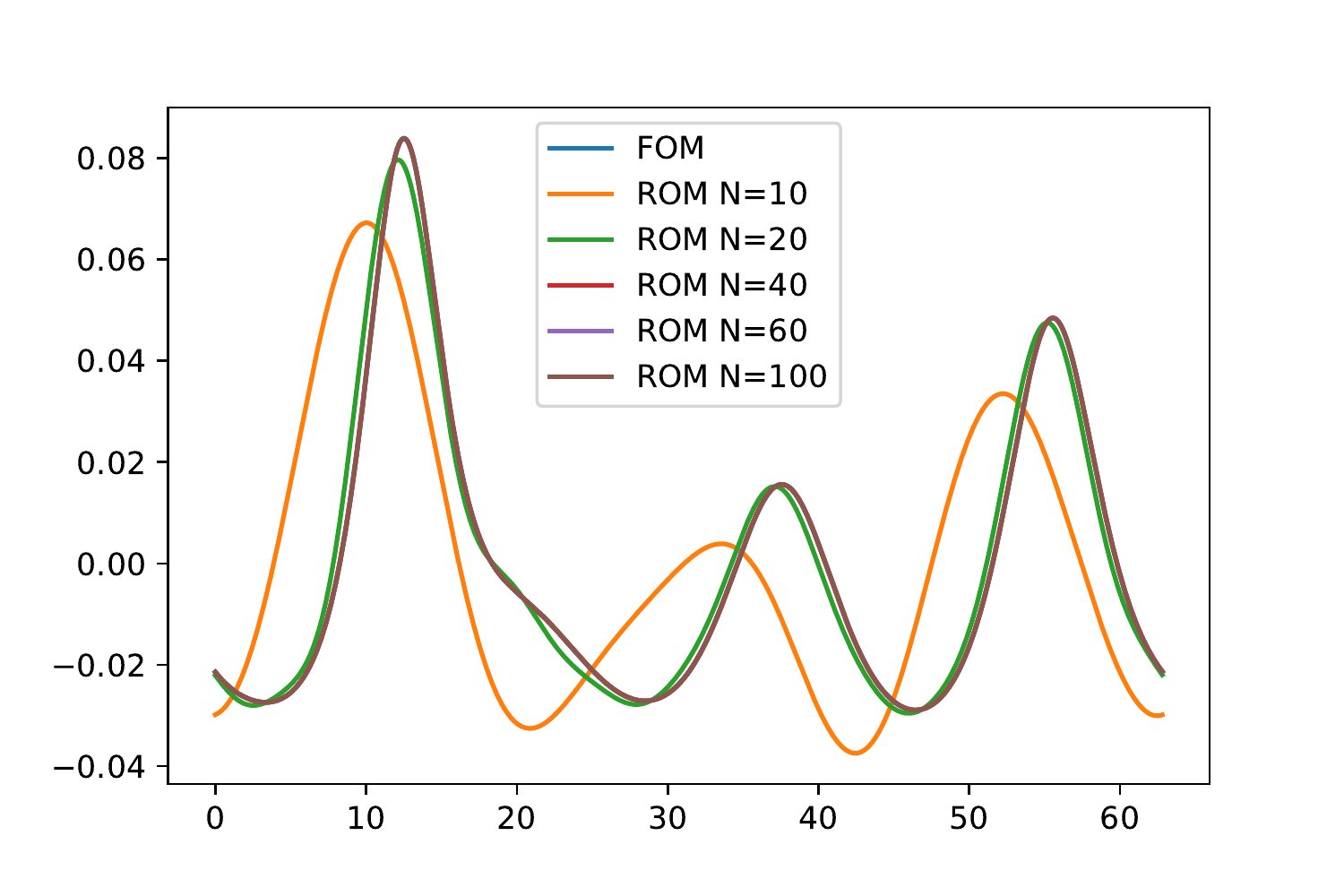}}\;
			\subfigure[ROM solutions time-parameter reduction with $tol_{POD}=0.3$, $\NRB=27$ and $tol_{EIM}=0.03$, $\NEIM=72$ on parameter and $\eta_0$ in the training set \label{fig:ROMEIMparamSimPerSolin}]{\includegraphics[width=0.32\textwidth]{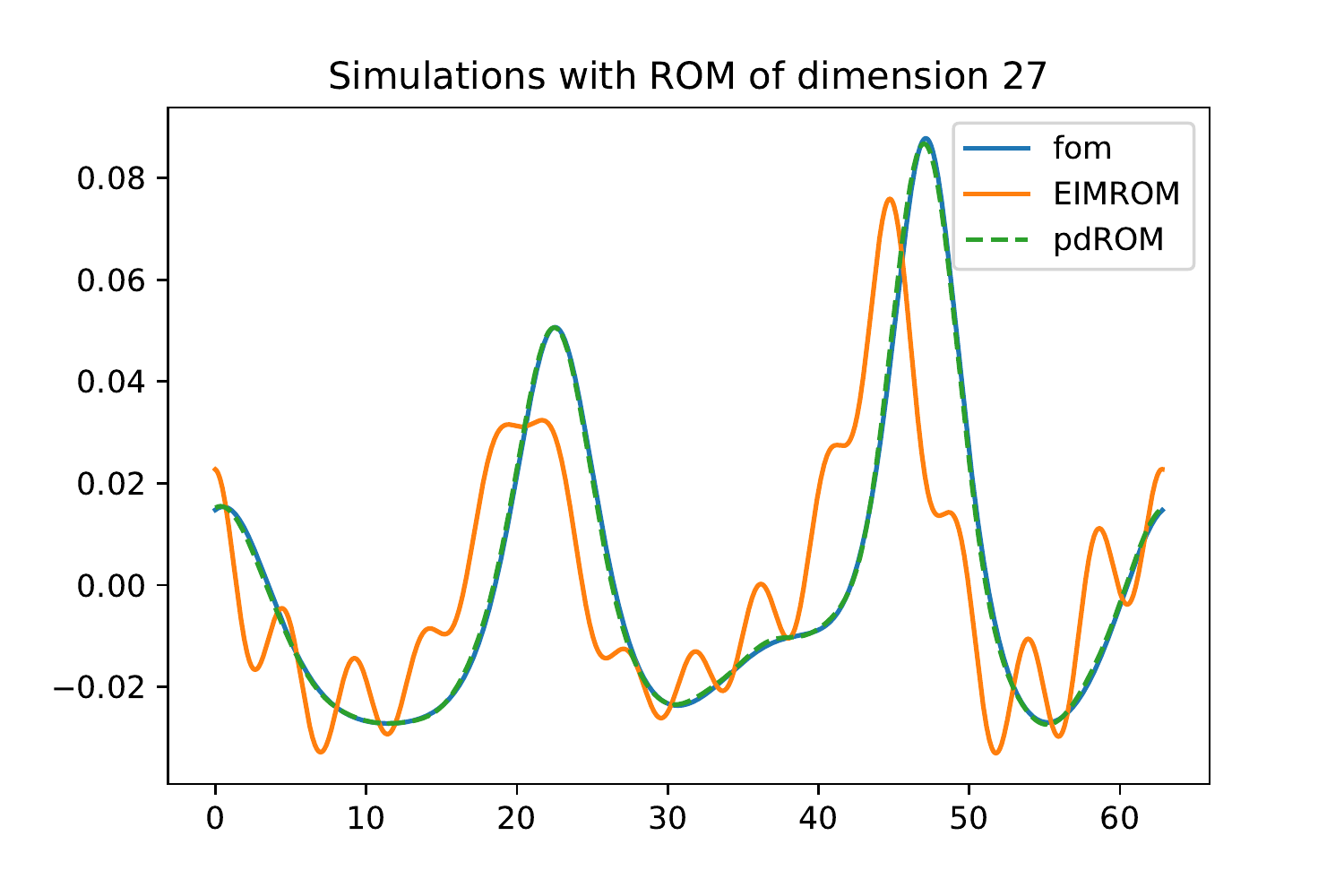}}\;
			\subfigure[ROM solutions time-parameter reduction with $\NRB=27$ and $\NEIM=72$ on parameter and $\eta_0$ outside the training set \label{fig:ROMEIMparamSimPerSolout}]{\includegraphics[width=0.32\textwidth]{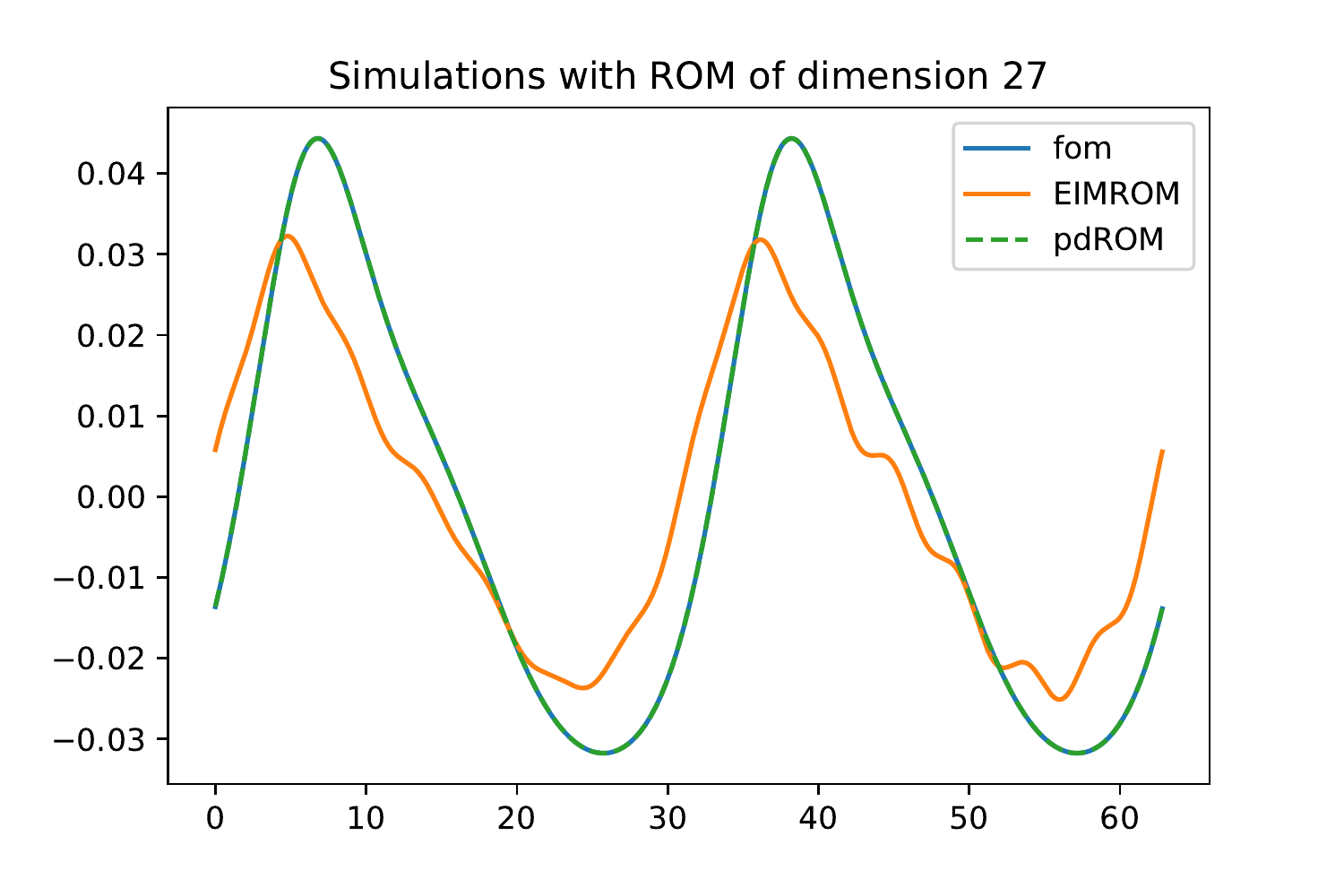}}\\
			\subfigure[ROM solutions time-only reduction \label{fig:ROMEIMsolPerSol}]{
				\includegraphics[width=0.32\textwidth]{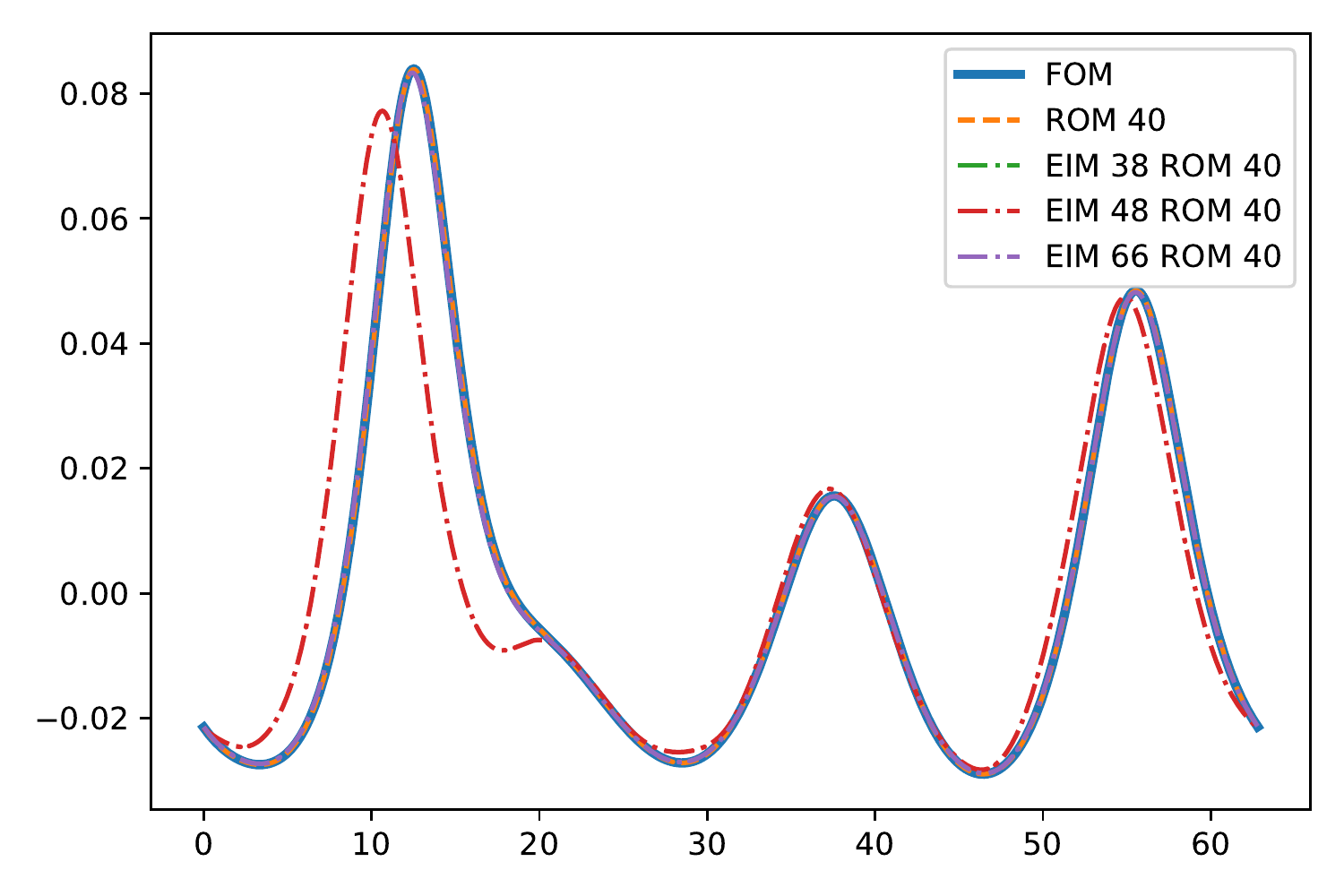}
			}\;
			\subfigure[ROM solutions time-parameter reduction with $tol_{POD}=0.3$, $\NRB=27$ and $tol_{EIM}=0.0003$, $\NEIM=119$ on parameter and $\eta_0$ in the training set \label{fig:ROMEIMparamSimPerSolin2}]{\includegraphics[width=0.32\textwidth]{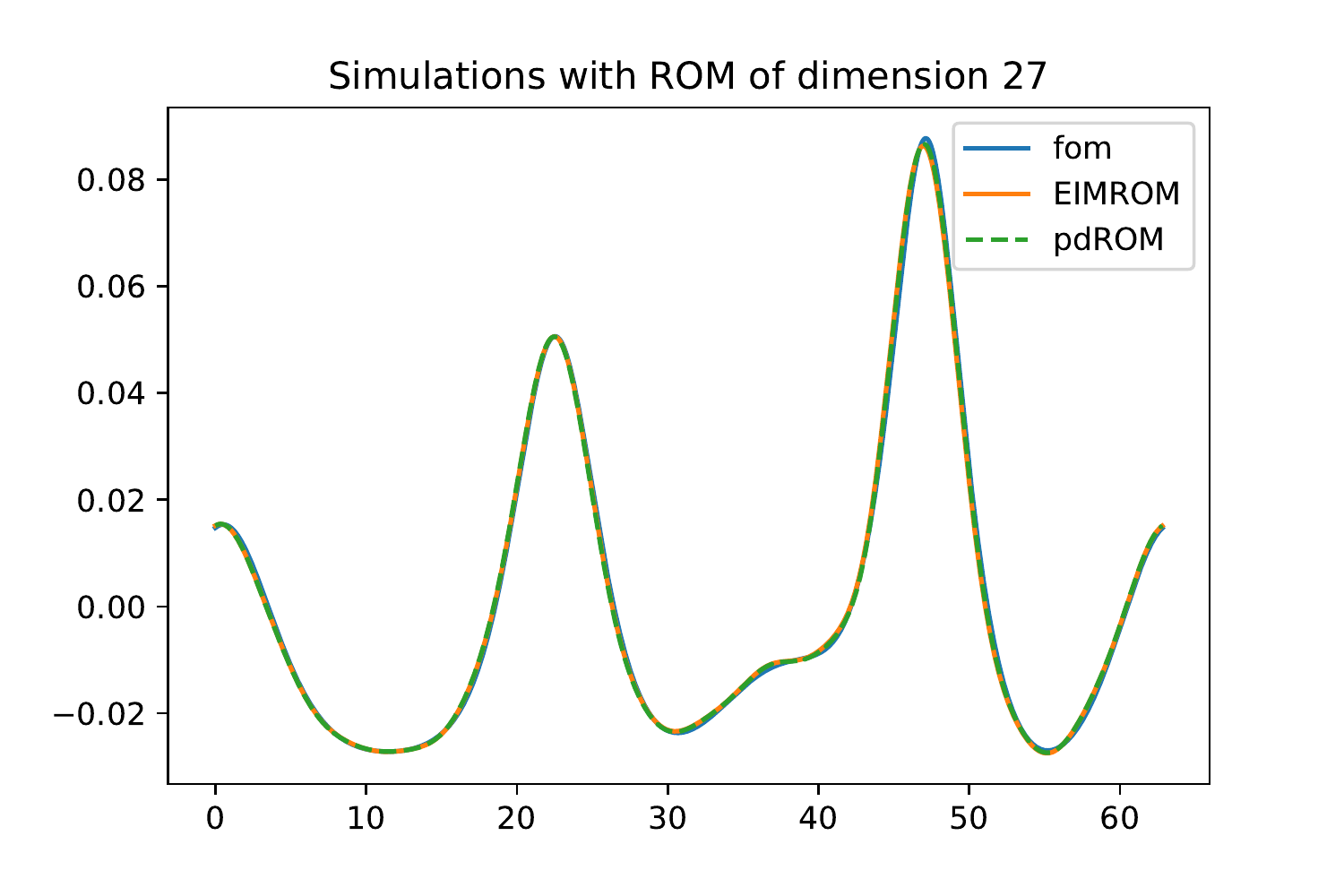}}\;
			\subfigure[ROM solutions time-parameter reduction with $\NRB=27$ and $\NEIM=119$ and $tol_{EIM}=0.0003$, $\NEIM=119$ on parameter and $\eta_0$ outside the training set \label{fig:ROMEIMparamSimPerSolout2}]{\includegraphics[width=0.32\textwidth]{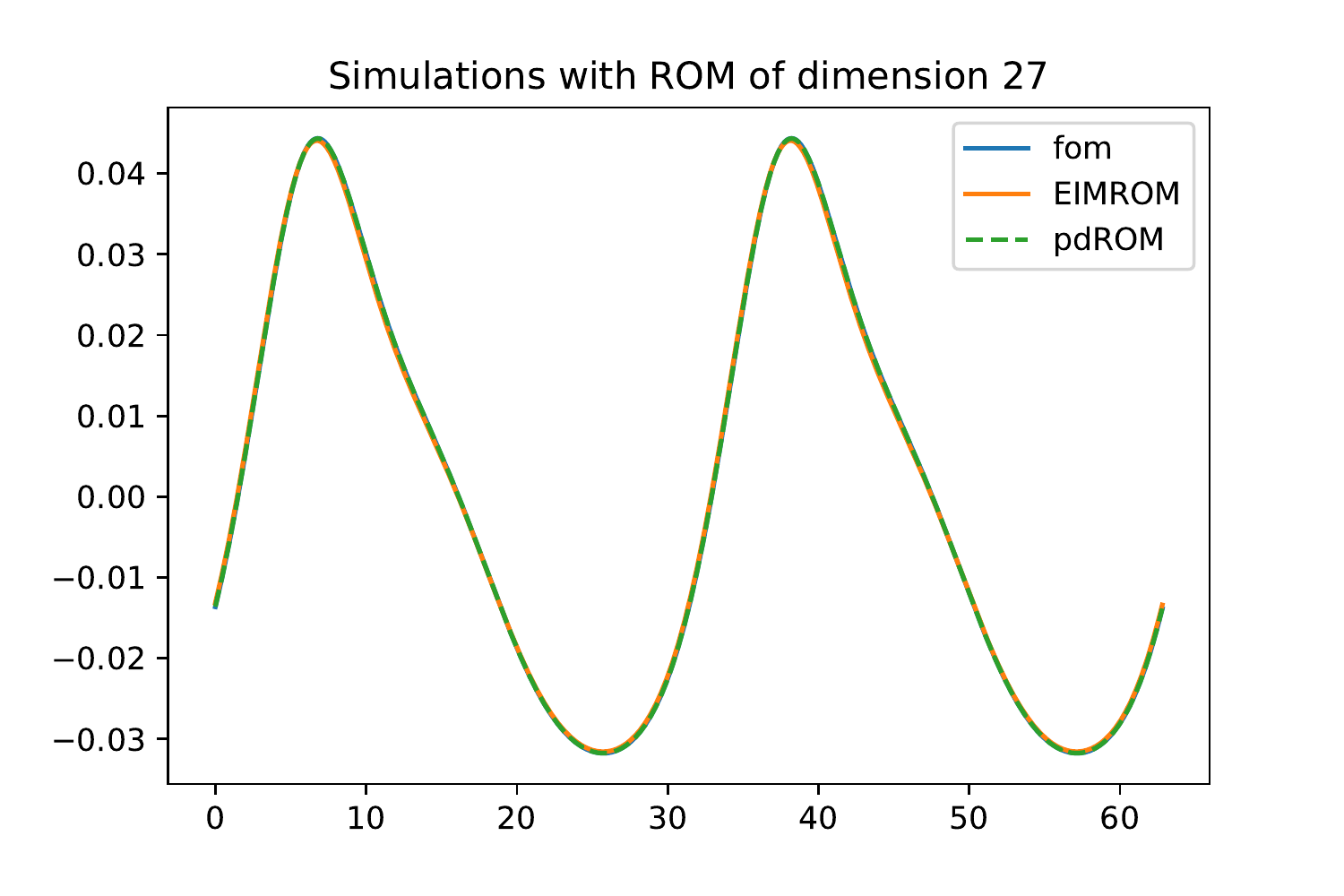}}
		\end{center}
		\caption{\textbf{Propagation of a periodic monochromatic wave.} pdROM and EIM reduction: simulation, errors and computational time, varying POD and EIM dimensions}
	\end{figure}
	\begin{figure}
		\centering
		\includegraphics[width=0.9\textwidth]{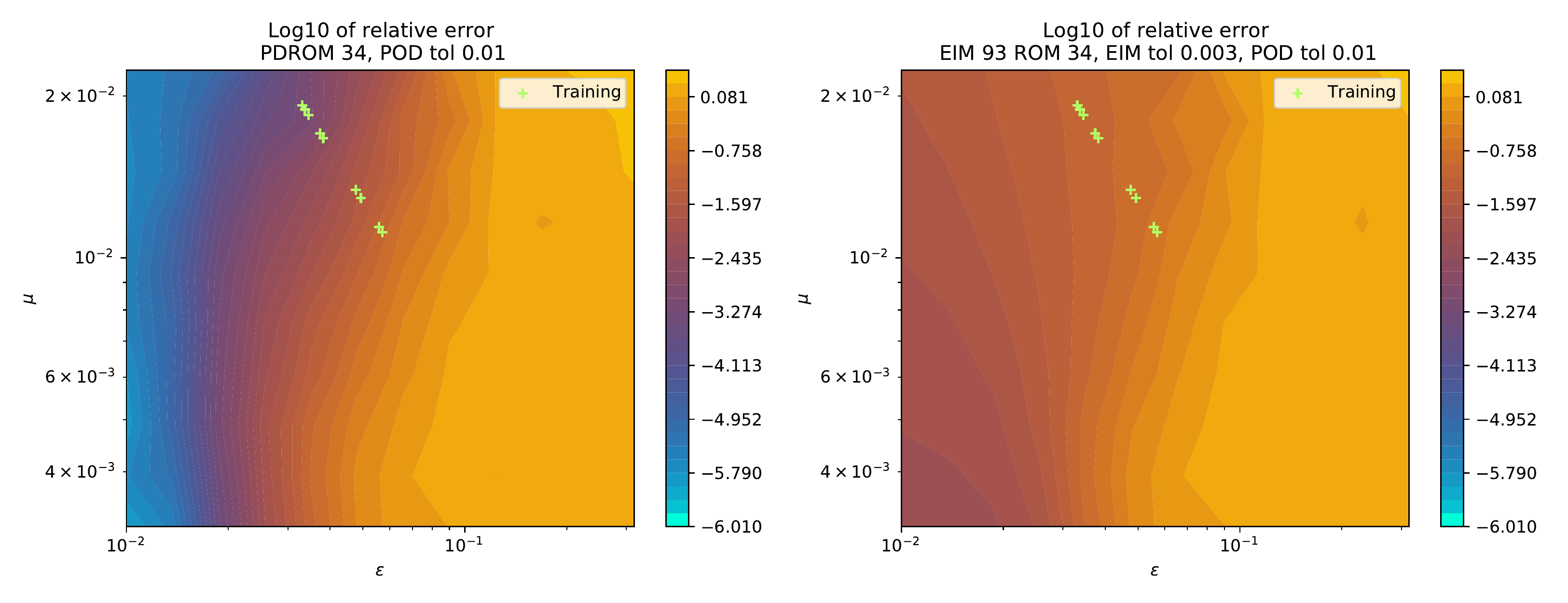}
		\caption{\label{fig:varyingPerSol}\textbf{Propagation of a periodic monochromatic wave.} Error plot varying the parameters $a_0$ and $h_0$ for pdROM $\NRB=34$ and EIMROM with $\NEIM= 93$}
	\end{figure}
	
	For this test, we use $\bar{h}_0=1$ and final time $T=200$.
	This test is fairly simple and it can be well compressed with few basis functions. Nevertheless, it already shows what can be achieved with and without EIM. We collect samples for different times of one simulation obtaining a time-only reduction. In \cref{fig:ROMEIMerrPerSol} the decay of the relative error is exponential with respect to the dimension $\NRB$ of the reduced space  (different points correspond to $\NRB =10k$ for $k=1,\dots, 10$) and with few basis functions we can reach a very small error. On the other side, this error could be worsen or even the simulation might not run until the final time if we add the EIM algorithm with not enough basis functions. The error for those methods stagnates around the tolerance of EIM and if $\NEIM\lessapprox \NRB$ the method becomes unstable.
	
	The computational time 
	for the KdV-BBM pdROM   compared to the FOM using Thomas method
	is  between 20\% and 30\% in \cref{fig:ROMEIMtimePerSol}. 
	The hyper-reduction allows to further reduce the 
	in between
	4\% and 10\% of the original computational time. In \cref{fig:ROMsolPerSol} we see the simulations run with pdROM with different dimensions of the reduced space. Already with $\NRB=20$ we obtain an accurate solution. In \cref{fig:ROMEIMsolPerSol} we see that a minimum amount of EIM basis is necessary not to have unstable solutions and to be accurate.
	
	For the second phase, we introduce the reduction in the parameter space. For  \cref{fig:ROMEIMparamErrPerSolin,fig:ROMEIMparamTimePerSolin,fig:ROMEIMparamSimPerSolin,fig:ROMEIMparamSimPerSolin2}, we test the obtained reduced spaces and the respective algorithms onto the parameter $\bar{h}_0=1$ and for final time $T=250$. 
	The error decay in \cref{fig:ROMEIMparamErrPerSolin} is similar to the one of the time-only reduced case. The largest differences stays in the differences of the number of basis functions needed to reach a certain tolerance. In the example in \cref{fig:ROMEIMparamSimPerSolin}, with $\NRB=27$ we reach a tolerance of $0.3$. 
	With time-only reduction $\NRB=10$ are enough to reach the same tolerance, though not reaching the same accuracy in the simulation see \cref{fig:ROMEIMerrPerSol}. 
	The same holds for the EIM algorithm. This is noticeable in the computational times, that increase for all algorithms given certain thresholds, see \cref{fig:ROMEIMtimePerSol,fig:ROMEIMparamTimePerSolin}, 
	this is a measure of the increased difficulty of the problem. 
	In \cref{fig:ROMEIMparamSimPerSolin} a representative simulation is shown, 
	where we see that with $\NRB=27$ basis function we get an already accurate approximation with pdROM, while $\NEIM=72$ is still not enough to obtain a good accuracy for EIMROM. 
	Increasing $\NEIM=119$ in \cref{fig:ROMEIMparamSimPerSolin2}, we can obtain an accurate solution also for EIMROM. 
	This is accordance with the observations made in \cite{peherstorfer2020stability,zimmermann2018geometric,ghavamian2017pod,argaud2017stabilization,chen2021eim}. 
	In this case, the computational time of EIMROM is 3 times faster than the pdROM, and decreasing $tol_{EIM}$ we can achieve even smaller errors. 
	A particular attention must be observed in choosing the tolerance of EIM small enough with respect to $tol_{POD}$, as one can see in \cref{fig:ROMEIMparamErrPerSolin} that the EIMROM simulations can explode if $tol_{EIM}$ is not small enough.

	In \cref{fig:ROMEIMparamErrPerSolout,fig:ROMEIMparamTimePerSolout,fig:ROMEIMparamSimPerSolout,fig:ROMEIMparamSimPerSolout2} we test the algorithms with $h_0=0.63\notin [h_{min},h_{max}]$, $a_0=0.05$ instead of 0.04 and initial conditions $u_0(x)=-1.2a_0\cos(\frac{4}{10}\pi(x-2))+\frac{a_0}{10}\cos(\frac{8}{10}\pi x)$.
	Similar conclusion can be drawn also in this case. The main difference is a larger error, due to the higher wave amplitude of the solution.
	
	In \cref{fig:varyingPerSol} we can see the relative error of pdROM $\NRB=34$ and EIMROM with $\NEIM=93$ varying $a_0$ and $h_0$ in the initial conditions and the problem settings. We plot the error as a function of $\varepsilon$ and $\mu$, in order to understand how the nonlinearity and the dispersion affect the solutions. We know that as $h_0$ gets smaller or $a_0$ gets bigger ($\varepsilon$ larger), the nonlinearity of the problem becomes more pronounced, hence leading to more oscillations and steeper gradients. In the strongly nonlinear regime ($\varepsilon = a_0/h_0 \gtrapprox 10^{-1}$) the hypothesis of having small $\varepsilon$ are not valid anymore and the error increases rapidly. On the other side, the dispersion coefficient $\mu=h_0^2/L^2$ does not affect much the error for a fixed $\varepsilon$, in \cref{fig:varyingPerSol}.
	Overall we can see that, as long as nonlinearity does not become predominant, the pdROM provides lower errors over the parameter space showing potential for a more robust approximation outside the training set.

	\subsection{Undular bore propagation}
	\renewcommand{\folderTest}{Test_damBreak}
	\begin{figure}
		\begin{center}
			\subfigure[EIM ROM errors time-only reduction \label{fig:ROMEIMerrDam}]{\includegraphics[width=0.32\textwidth]{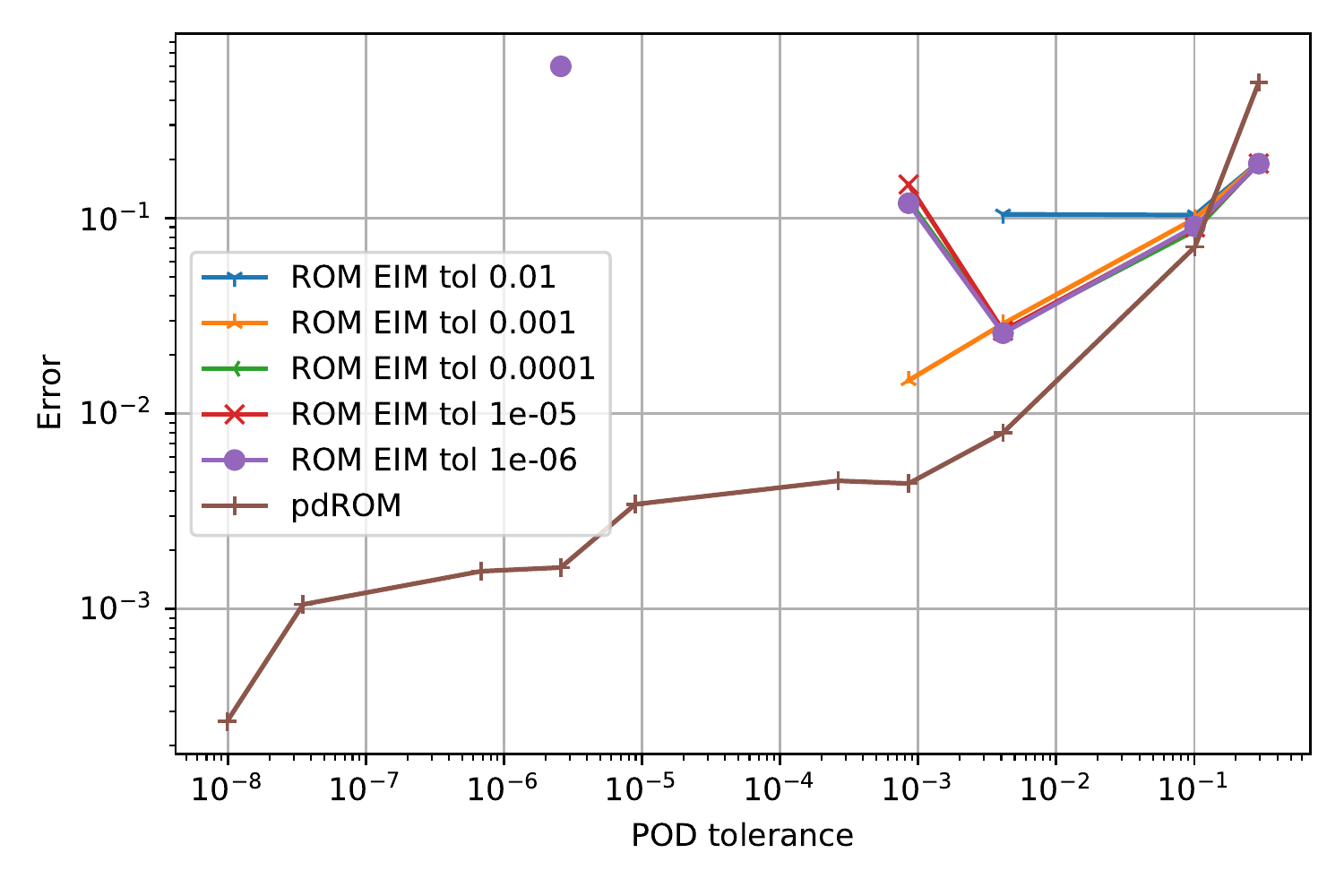}}\;
			\subfigure[ROM errors time-parameter reduction on parameter and $\eta_0$ in the training set \label{fig:ROMEIMparamErrDamin}]{\includegraphics[width=0.32\textwidth]{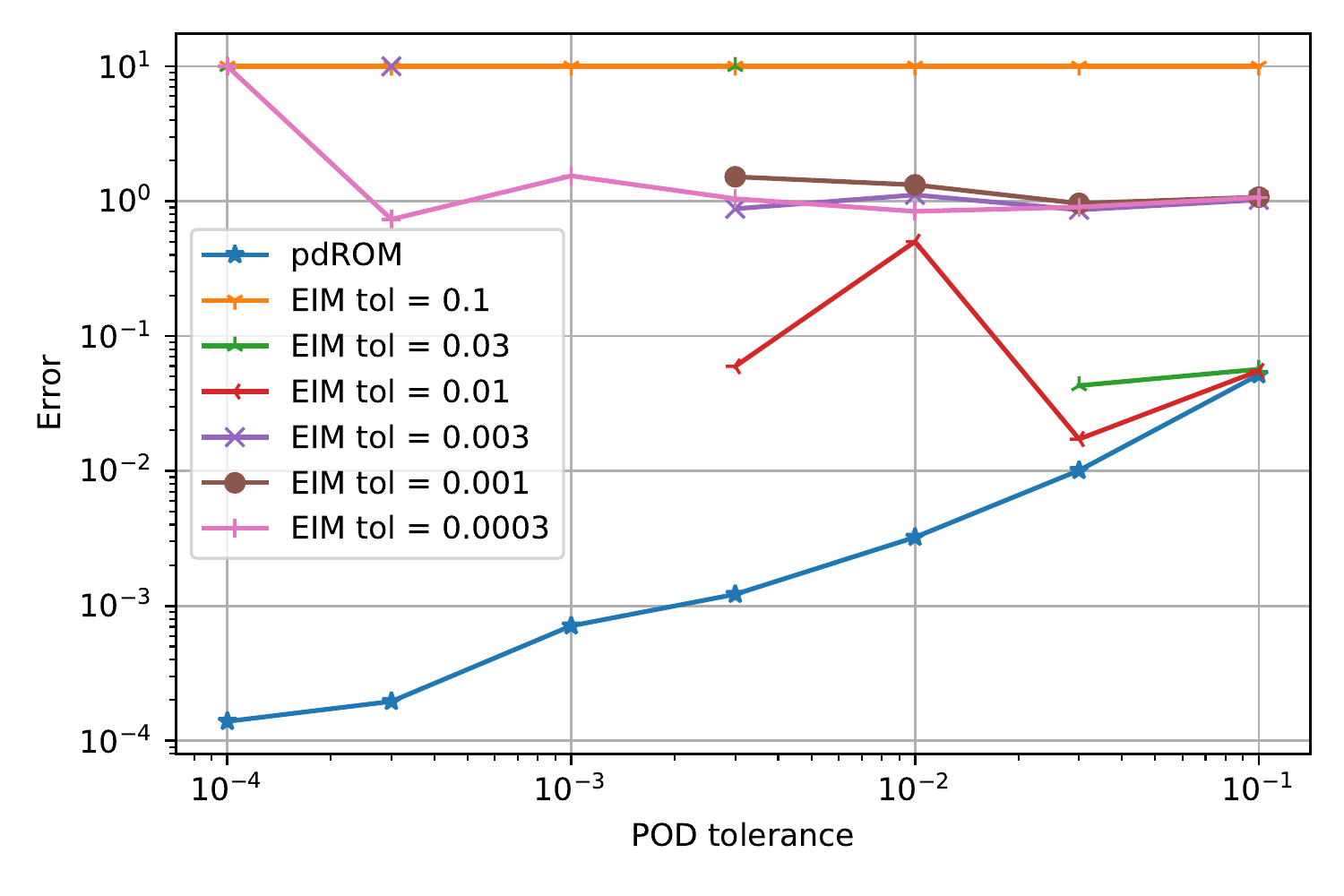}}\;
			\subfigure[ROM errors time-parameter reduction on parameter and $\eta_0$ outside the training set \label{fig:ROMEIMparamErrDamout}]{\includegraphics[width=0.32\textwidth]{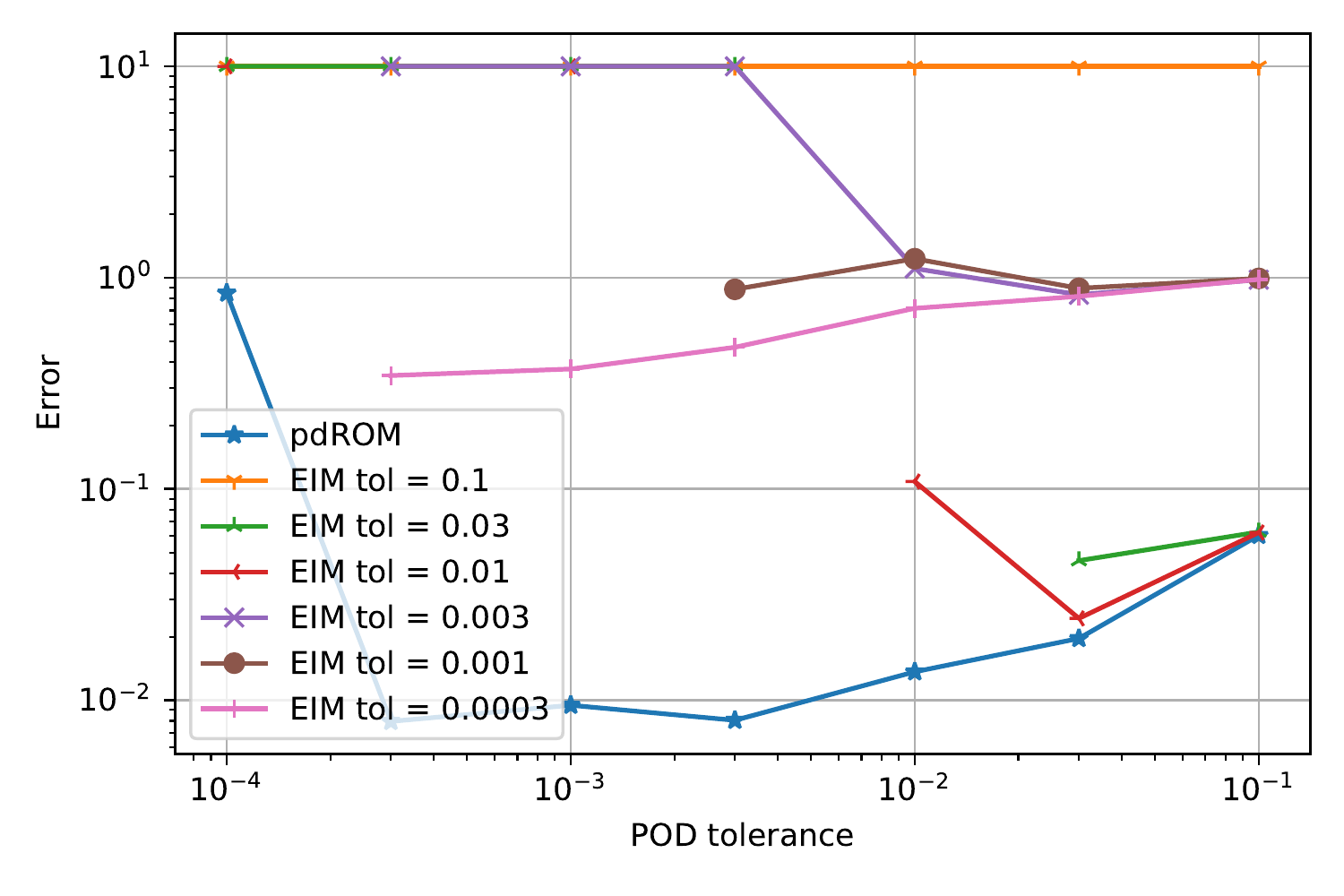}}
			
			\subfigure[EIM ROM computational times time-only reduction (legend above) \label{fig:ROMEIMtimeDam}]{\includegraphics[width=0.32\textwidth]{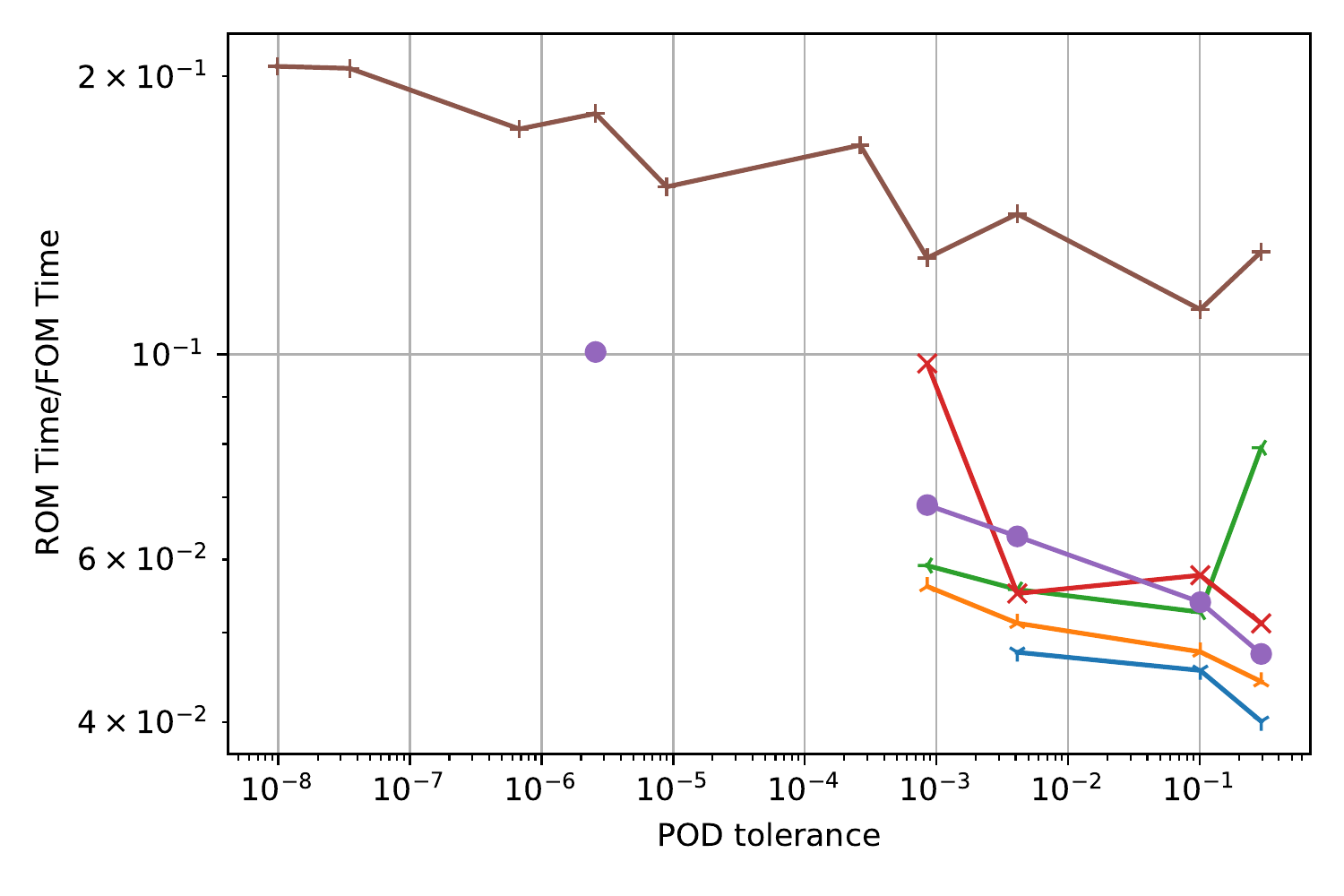}}\;
			\subfigure[ROM computational time time-parameter reduction on parameter and $\eta_0$ in the training set (legend above) \label{fig:ROMEIMparamTimeDamin}]{\includegraphics[width=0.32\textwidth]{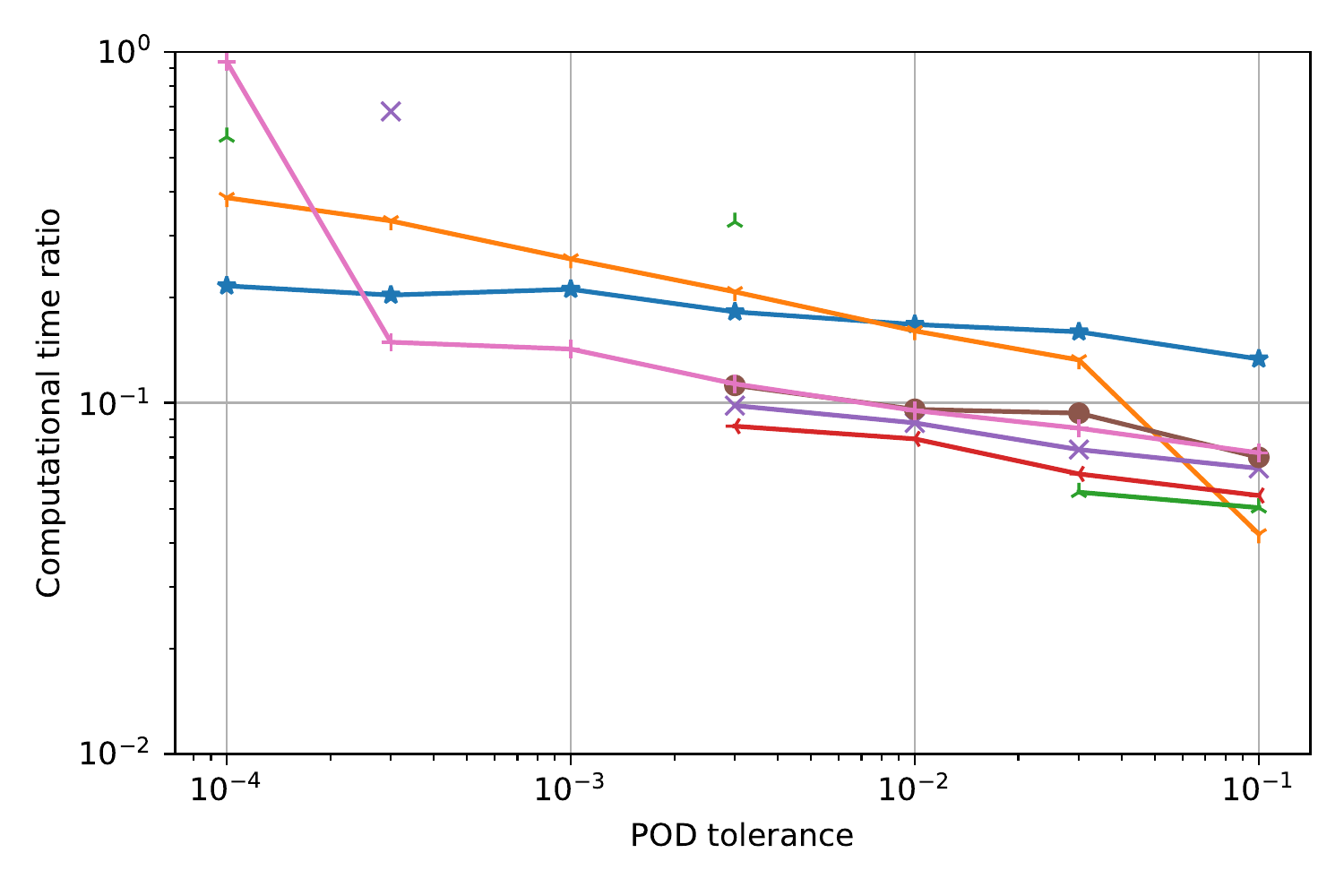}}\;
			\subfigure[ROM computational time time-parameter reduction on parameter and $\eta_0$ outside the training set (legend above) \label{fig:ROMEIMparamTimeDamout}]{\includegraphics[width=0.32\textwidth]{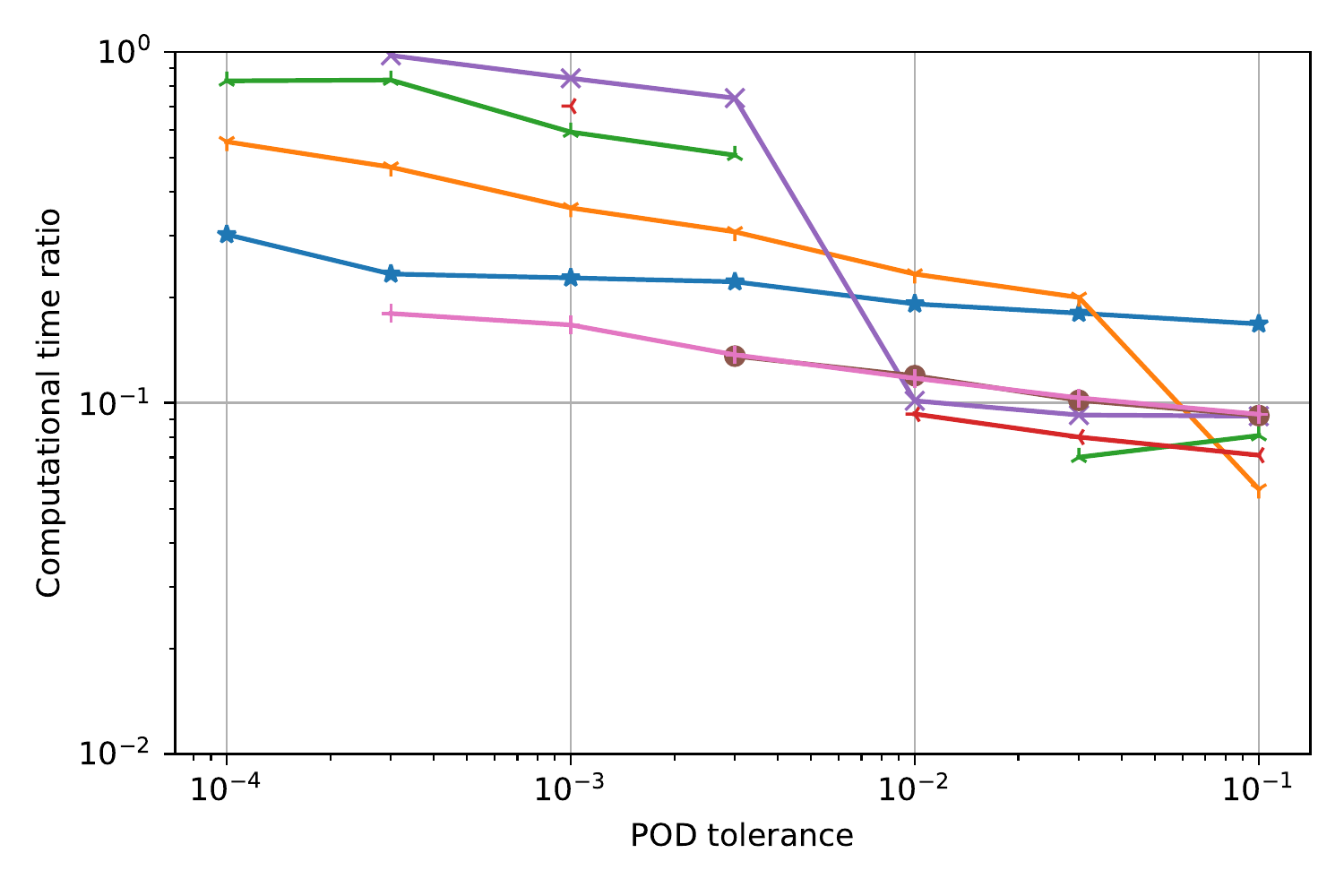}}

			\subfigure[ROM solutions time-only reduction \label{fig:ROMsolDam}]{\includegraphics[width=0.32\textwidth]{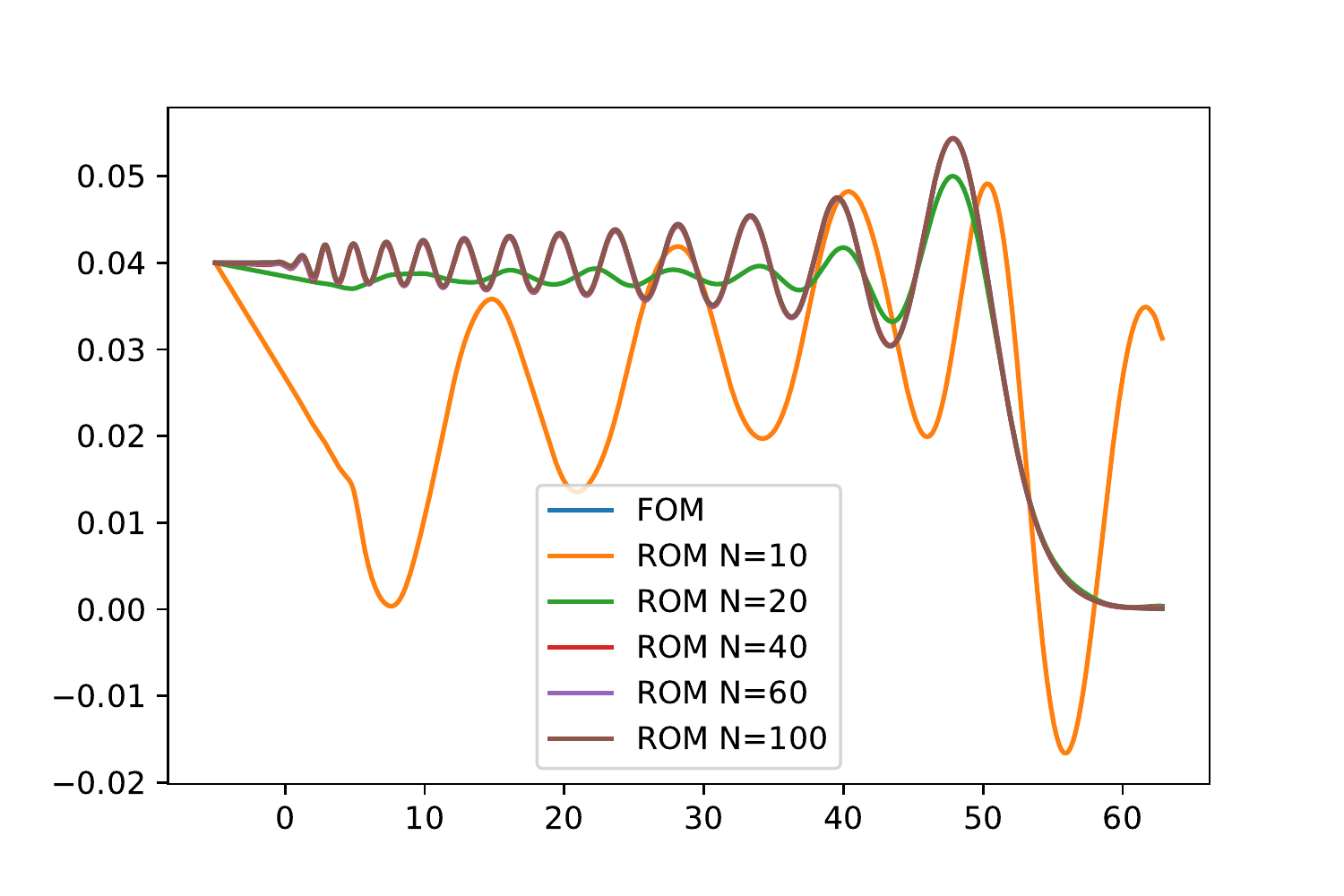}}\;
			\subfigure[ROM solutions time-parameter reduction with $tol_{POD}=0.03$, $\NRB=33$ and $tol_{EIM}=0.03$, $\NEIM=52$ on parameter and $\eta_0$ in the training set \label{fig:ROMEIMparamSimDamin}]{\includegraphics[width=0.32\textwidth]{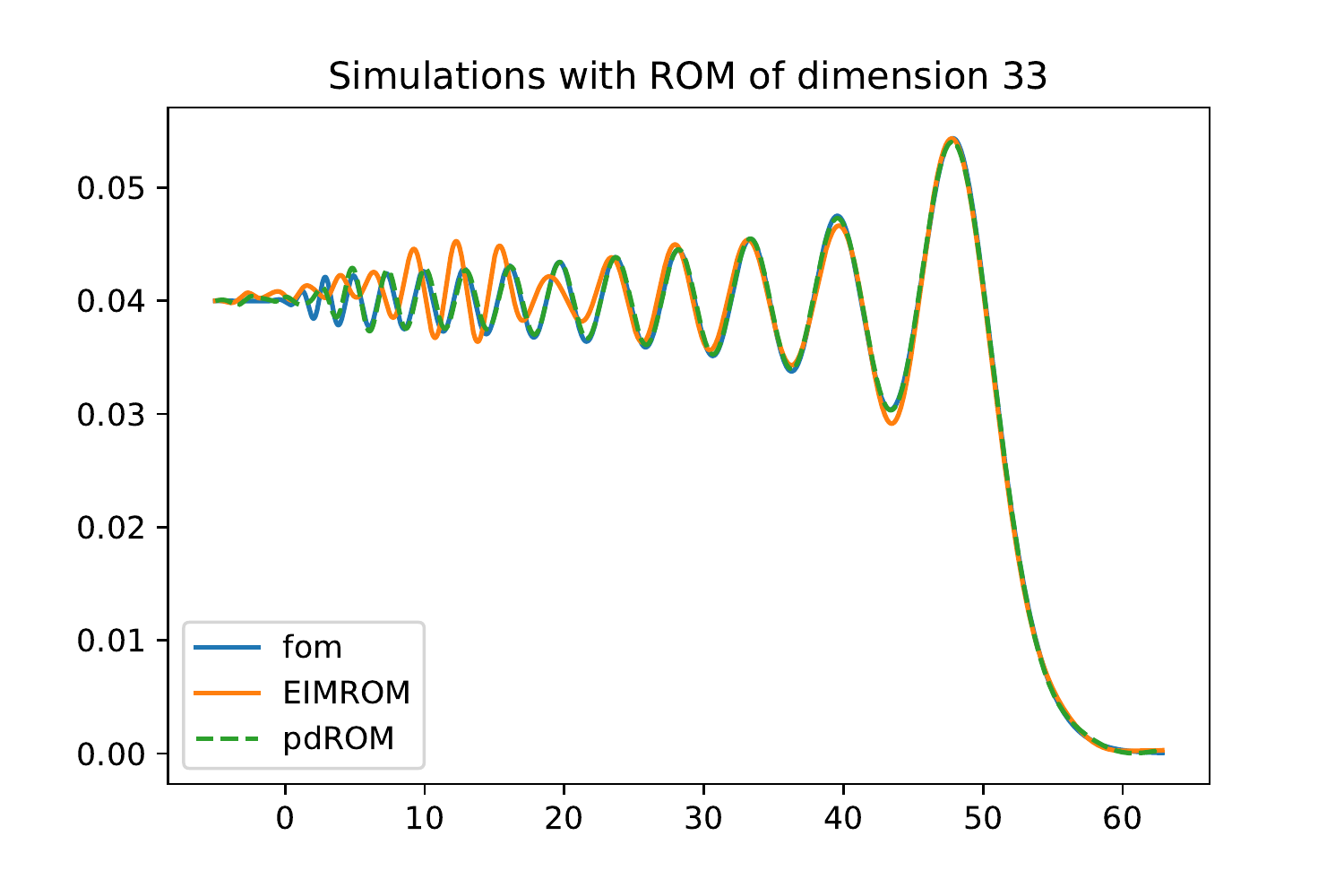}}\;
			\subfigure[ROM solutions time-parameter reduction with $tol_{POD}=0.03$, $\NRB=33$ and $tol_{EIM}=0.03$, $\NEIM=52$ on parameter and $\eta_0$ outside the training set \label{fig:ROMEIMparamSimDamout}]
			{	\includegraphics[width=0.32\textwidth]{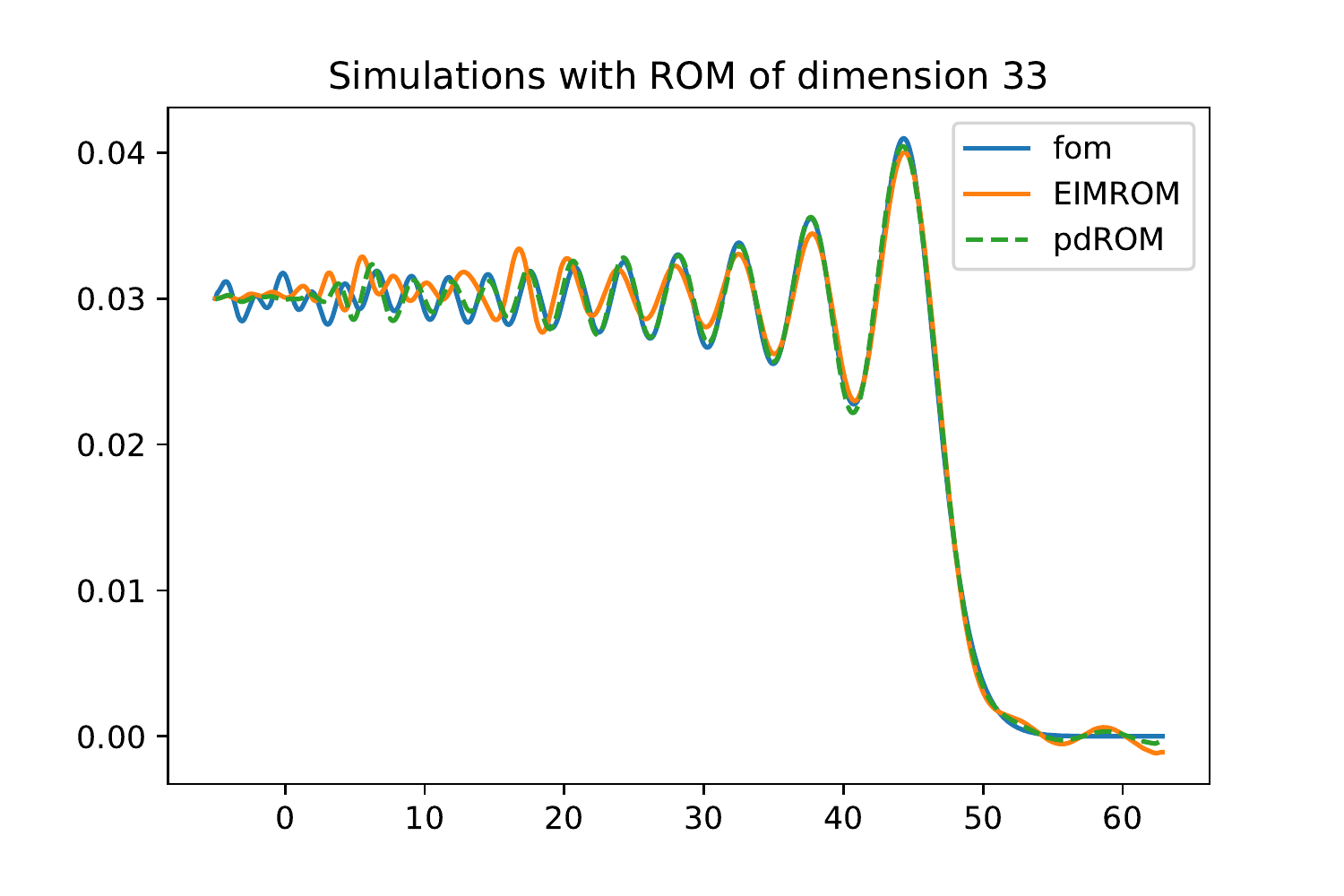}}
			
			\subfigure[EIMROM solutions time-only reduction \label{fig:ROMEIMsolDam}]{
				\includegraphics[width=0.32\textwidth]{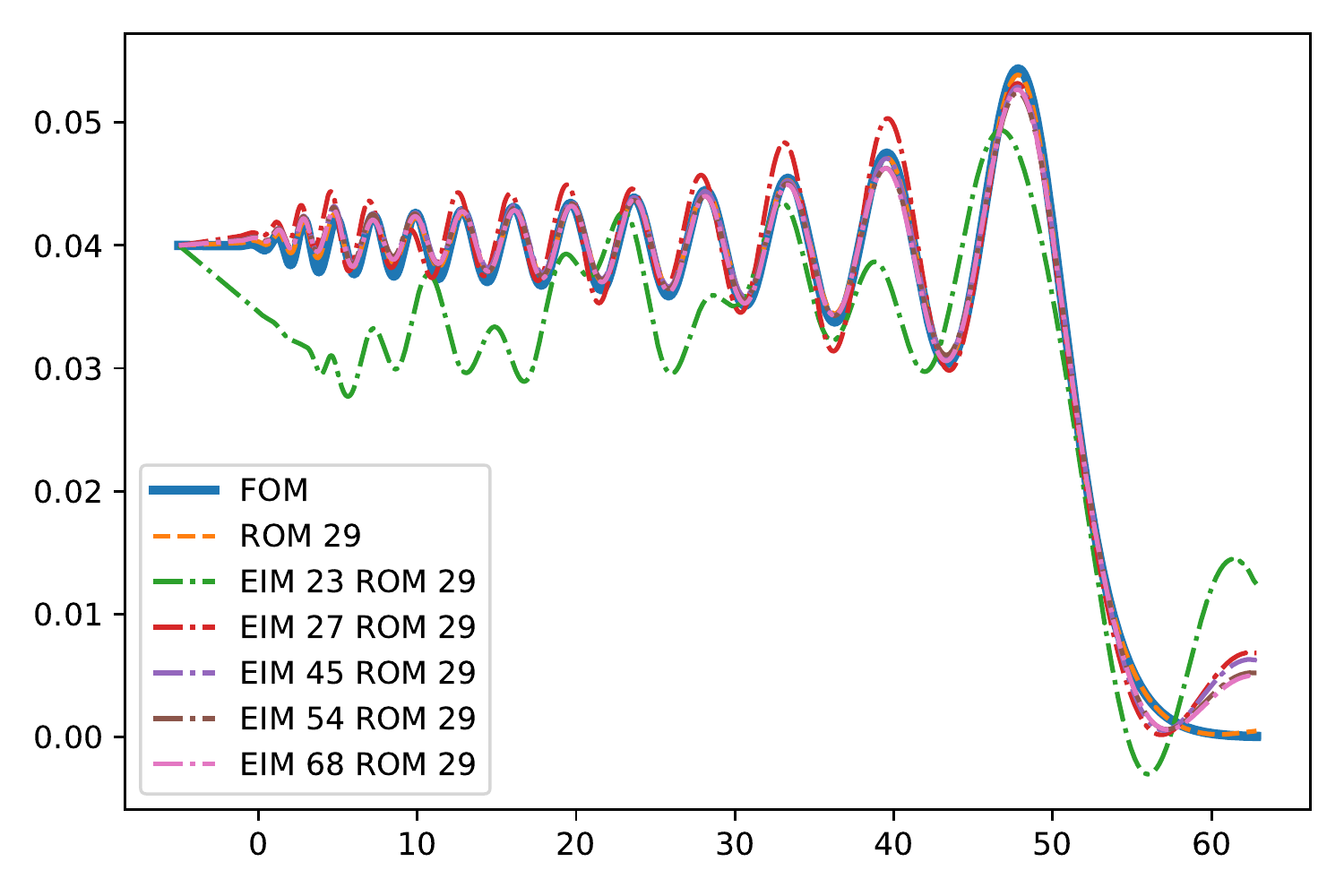}
			}\;
			\subfigure[ROM solutions time-parameter reduction with $tol_{POD}=0.003$, $\NRB=62$ and $tol_{EIM}=0.01$, $\NEIM=77$ on parameter and $\eta_0$ in the training set \label{fig:ROMEIMparamSimDamin2}]{\includegraphics[width=0.32\textwidth]{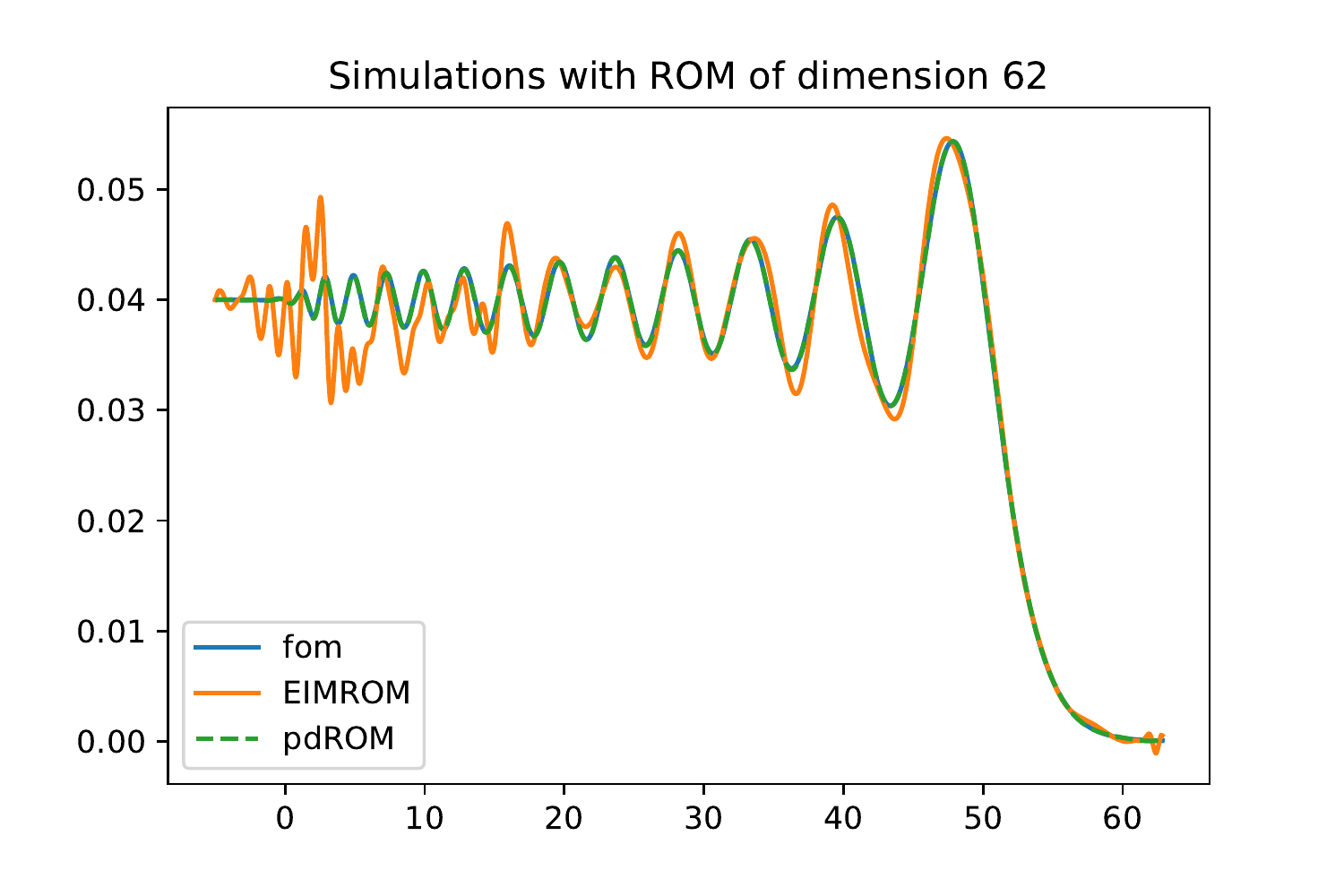}}\;
			\subfigure[ROM solutions time-parameter reduction with $tol_{POD}=0.0001$, $\NRB=103$ on parameter and $\eta_0$ outside the training set \label{fig:ROMEIMparamSimDamout2}]{\includegraphics[width=0.32\textwidth]{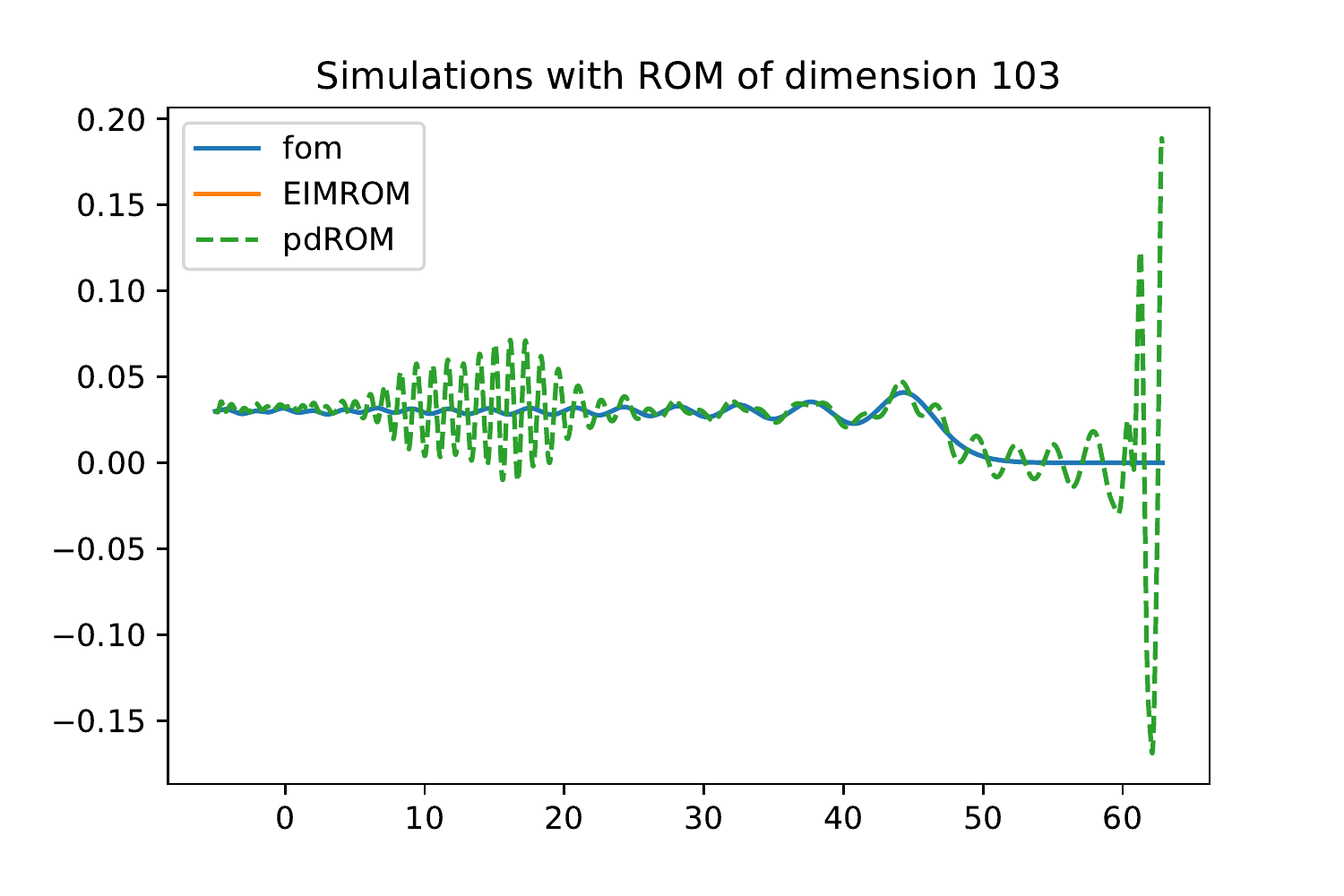}}
			
		\end{center}
		\caption{\textbf{Undular bore propagation.} pdROM and EIM reduction: simulation, errors and computational time, varying POD and EIM dimensions}
	\end{figure}
	\begin{figure}
		\centering
		\includegraphics[width=0.45\textwidth, trim={0 0 480 0},clip]{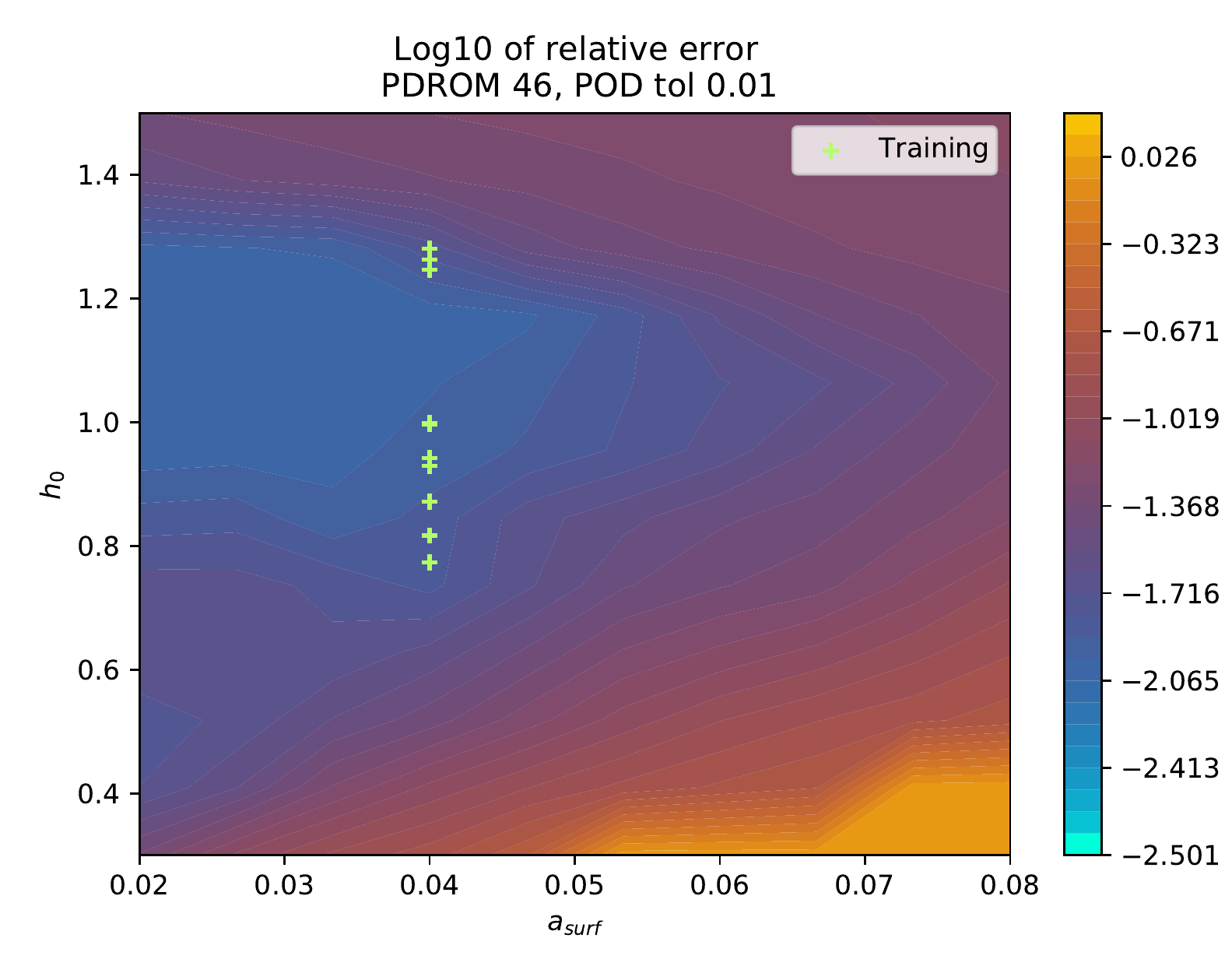}
		\caption{\label{fig:varyingDam}\textbf{Undular bore propagation.} Error plot varying the parameters $a_0$ and $h_0$ for pdROM $\NRB=46$}
\end{figure}
This is a very complex test for model reduction. 
We start by choosing $\bar{h}_0=1$, $a_0=0.04$ and final time $T=20$. In order to take into account the non homogeneous Dirichlet boundary condition at the left boundary, we subtract from the method the residual equation computed on a lift solution $u_{lift}$ that verify the Dirichlet boundary conditions. In this way, the resulting scheme can be applied with homogeneous Dirichlet boundary condition on the left. 

This tests is way more challenging than the previous one, in particular because flat zones and oscillating ones coexist and move around the domain. In particular, the EIM method easily fails in capturing this behavior without falling into Gibbs phenomena in the flat area. Indeed, in \cref{fig:ROMEIMerrDam} we see the error decay for only-time reduction both for EIMROM and pdROM and the former become very oscillatory when too many POD bases are used. The pdROM uses between 10\% and 20\% of the FOM computational time, while the EIMROM, in the few situations where it does not explode, uses around 5\% of the FOM time, see \cref{fig:ROMEIMtimeDam}. In \cref{fig:ROMsolDam} we see for different POD basis functions how we approach the exact solution, while in \cref{fig:ROMEIMsolDam} we fix the number of POD basis functions very low and we change the EIM ones. In this stable situation as the EIM bases increase we slowly converge towards the pdROM solution, even if the right part of the solution struggles with obtaining the right behavior. 

When exploring the parameter domain in the training set $h_0\in [0.7,1.3]$, we observe a similar situation. The decay of the error on a parameter inside the training domain computed at time $T=15$ for pdROM decays exponentially in \cref{fig:ROMEIMparamErrDamin}, while EIMROM is affected by oscillations and struggles with obtaining reasonable results. The computational time of the pdROM is between the 15\% and the 25\% of the FOM computational time as shown in \cref{fig:ROMEIMparamTimeDamin}, while the EIMROM, when the simulation is reasonable for few basis functions, needs only around the 5\% of the FOM computational time. 
In \cref{fig:ROMEIMparamSimDamin} we see for a tolerance of the POD of $3\cdot 10^{-2}$ we need $\NRB=33$ and the pdROM is quite close to the FOM solution, while EIMROM with same tolerance and $\NEIM=52$ start oscillating around the original position of the dam. In \cref{fig:ROMEIMparamSimDamin2} the situation worsen for the EIMROM which becomes very oscillatory for $\NRB=62$ and $\NEIM=77$ and this explains why the error is so large. 

Taking the parameter out of the training set, considering $h_0=0.63$ and $a_0=0.03$, the EIMROM becomes even more oscillatory and its error often explodes, but also the pdROM suffer of the larger flat zone in the part of the domain, leading to instabilities. The error decay is indeed affected by this and when too many POD bases are considered the error grows, see \cref{fig:ROMEIMparamErrDamout} and the simulation in \cref{fig:ROMEIMparamSimDamout2}. 
Computational times stays similar to the previous case \cref{fig:ROMEIMparamTimeDamout}, and, if properly choosing the number of basis functions, both pdROM and EIMROM  gives decent solutions, though the EIMROM oscillates around the flat area, see \cref{fig:ROMEIMparamSimDamout}.


\renewcommand{\folderTest}{Test_periodicSolitonOnBump}
\begin{figure}
	\begin{center}
		\subfigure[EIM ROM errors time-only reduction \label{fig:ROMEIMerrSolBump}]{\includegraphics[width=0.32\textwidth]{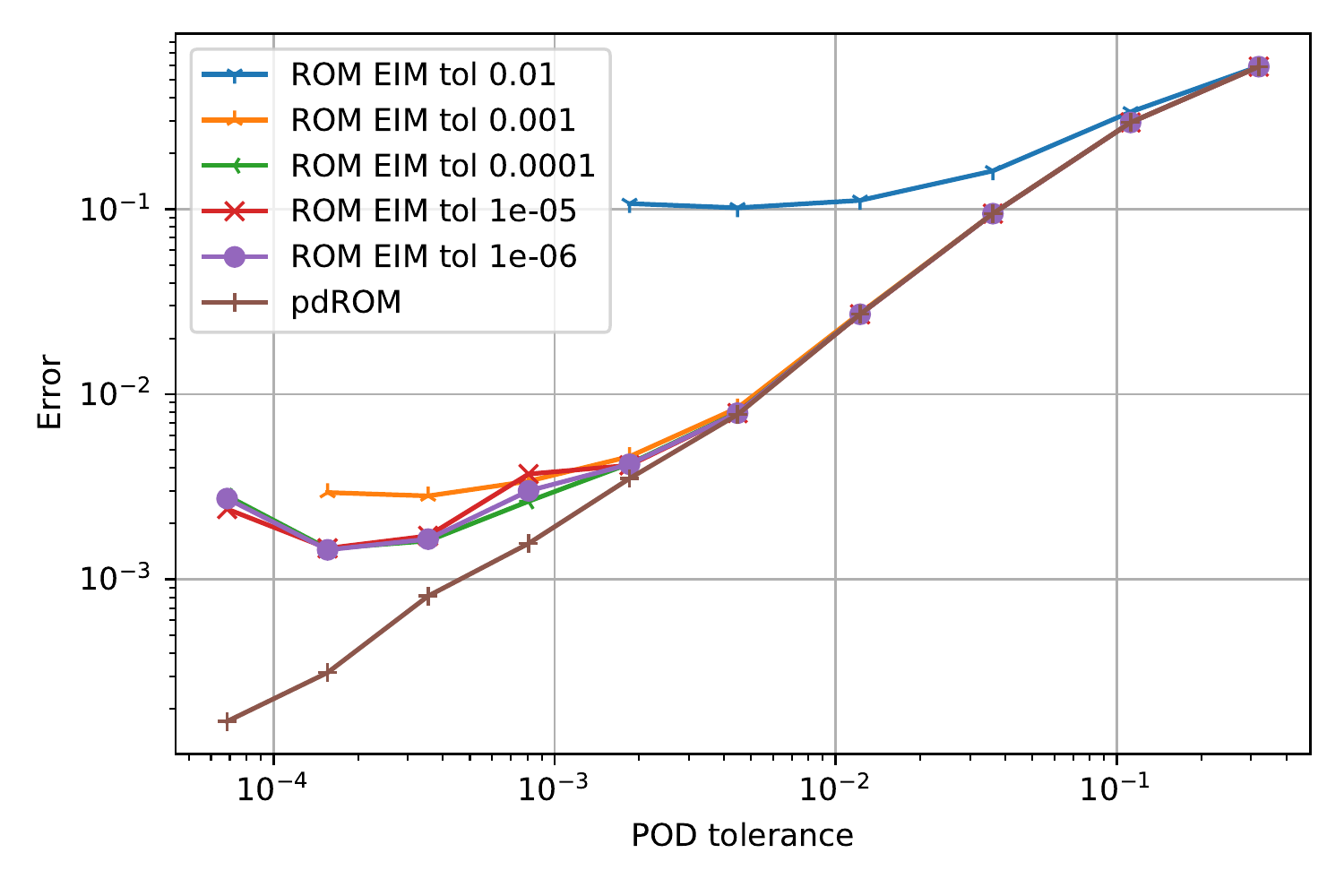}}\;
		\subfigure[ROM errors time-parameter reduction on parameter and $\eta_0$ in the training set \label{fig:ROMEIMparamErrSolBumpin}]{\includegraphics[width=0.32\textwidth]{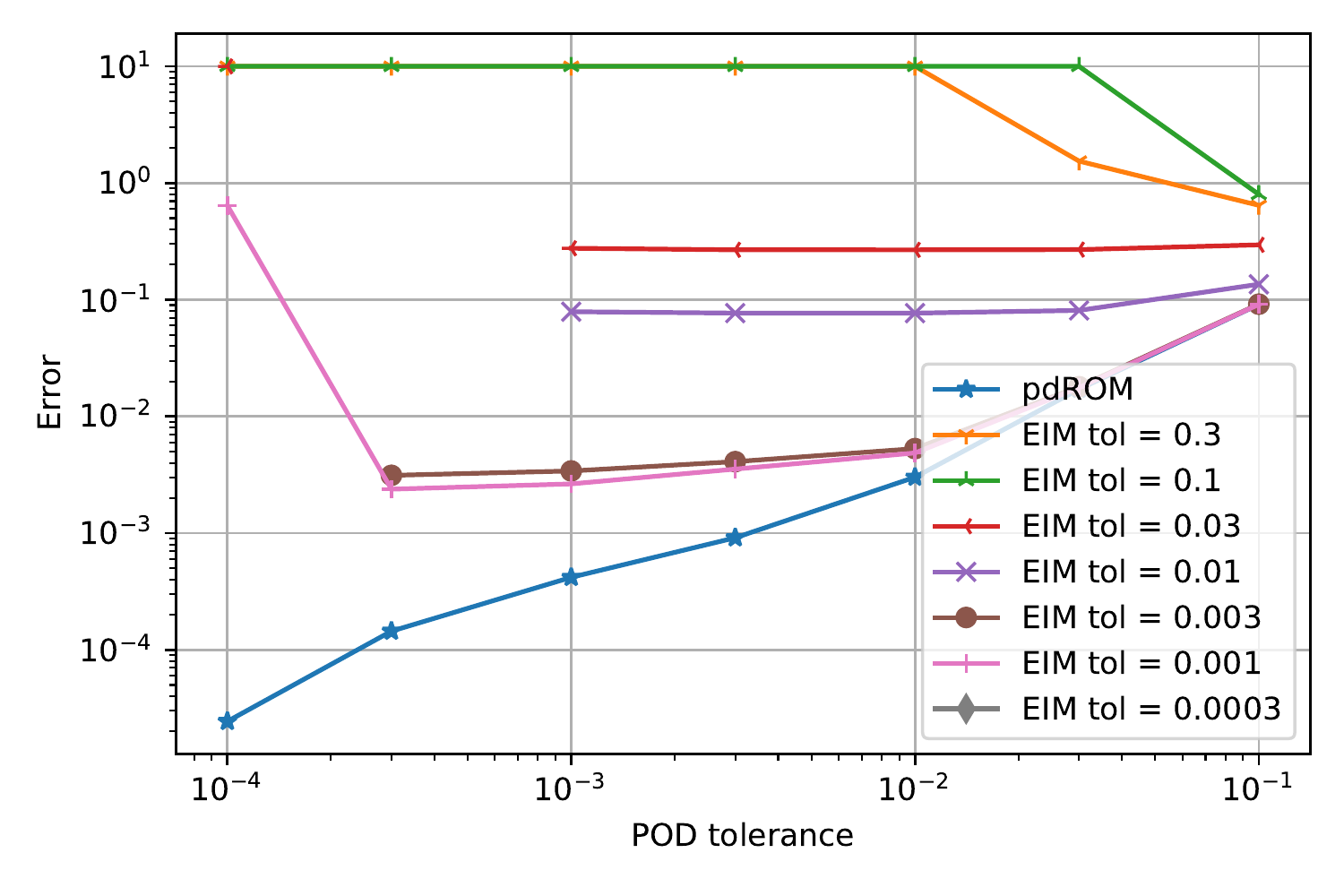}}\;
		\subfigure[ROM errors time-parameter reduction on parameter and $\eta_0$ outside the training set \label{fig:ROMEIMparamErrSolBumpout}]{\includegraphics[width=0.32\textwidth]{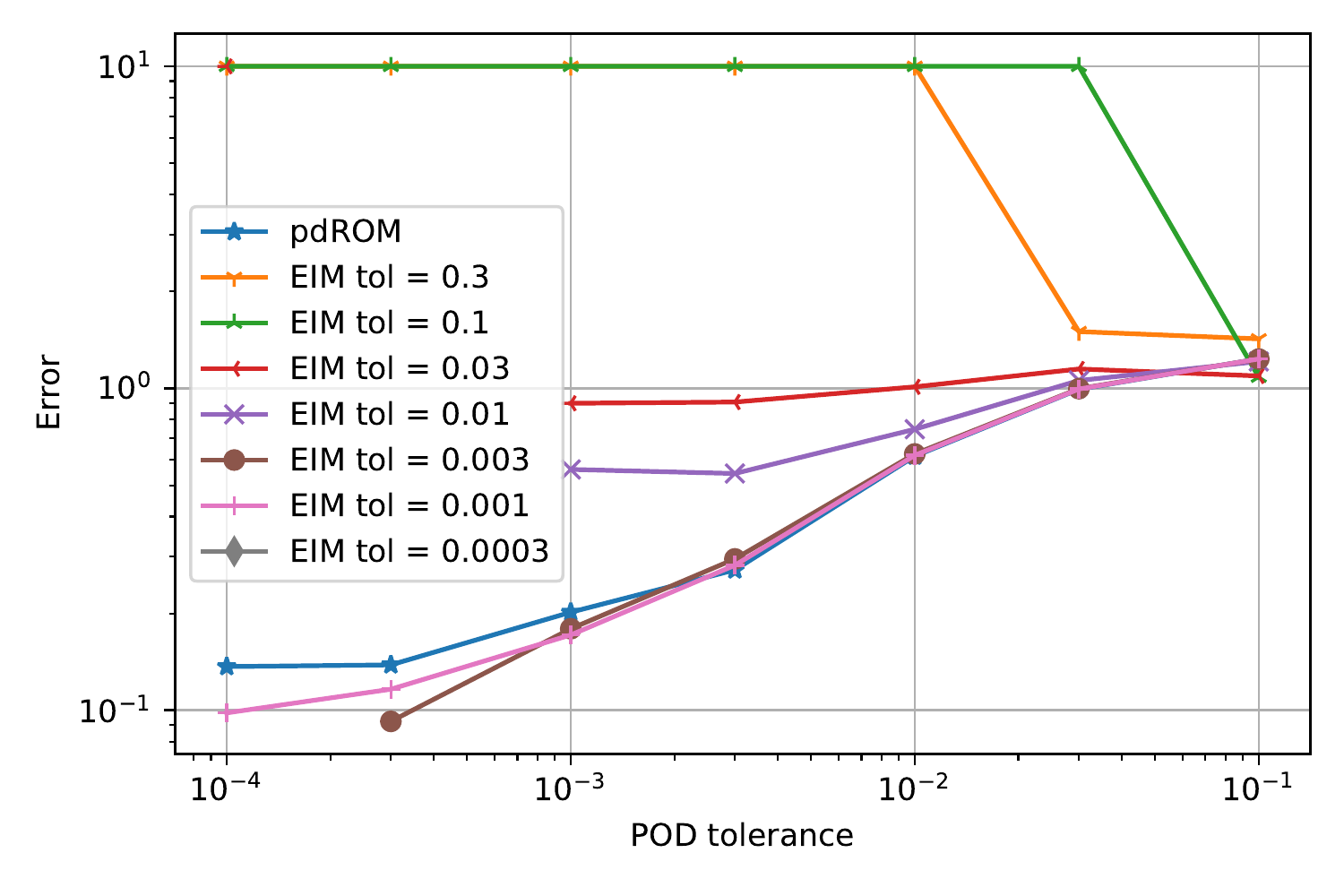}}

		\subfigure[EIM ROM computational times time-only reduction (legend above) \label{fig:ROMEIMtimeSolBump}]{\includegraphics[width=0.32\textwidth]{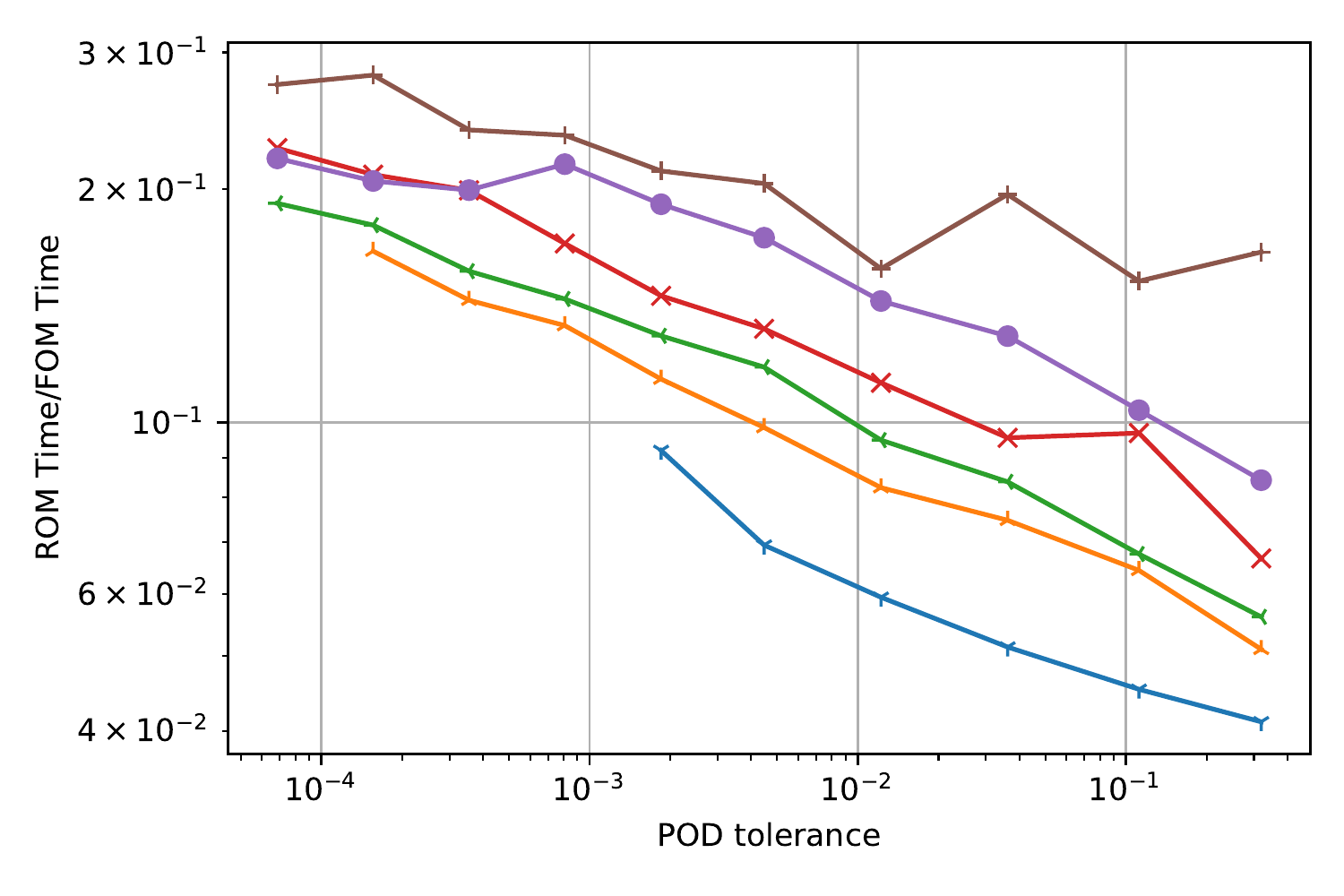}}\;
		\subfigure[ROM computational time time-parameter reduction on parameter and $\eta_0$ in the training set (legend above) \label{fig:ROMEIMparamTimeSolBumpin}]{\includegraphics[width=0.32\textwidth]{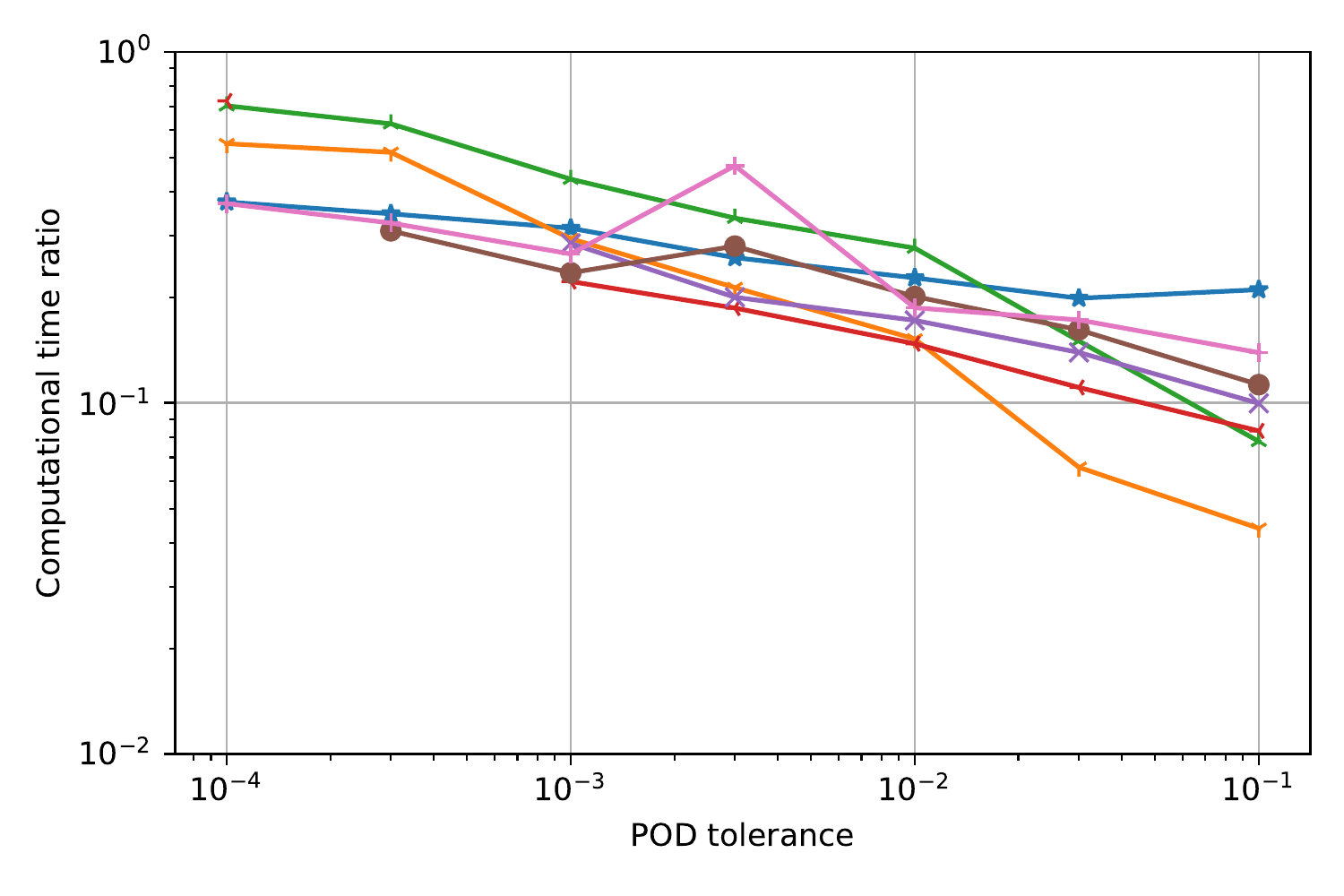}}\;
		\subfigure[ROM computational time time-parameter reduction on parameter and $\eta_0$ outside the training set (legend above) \label{fig:ROMEIMparamTimeSolBumpout}]{\includegraphics[width=0.32\textwidth]{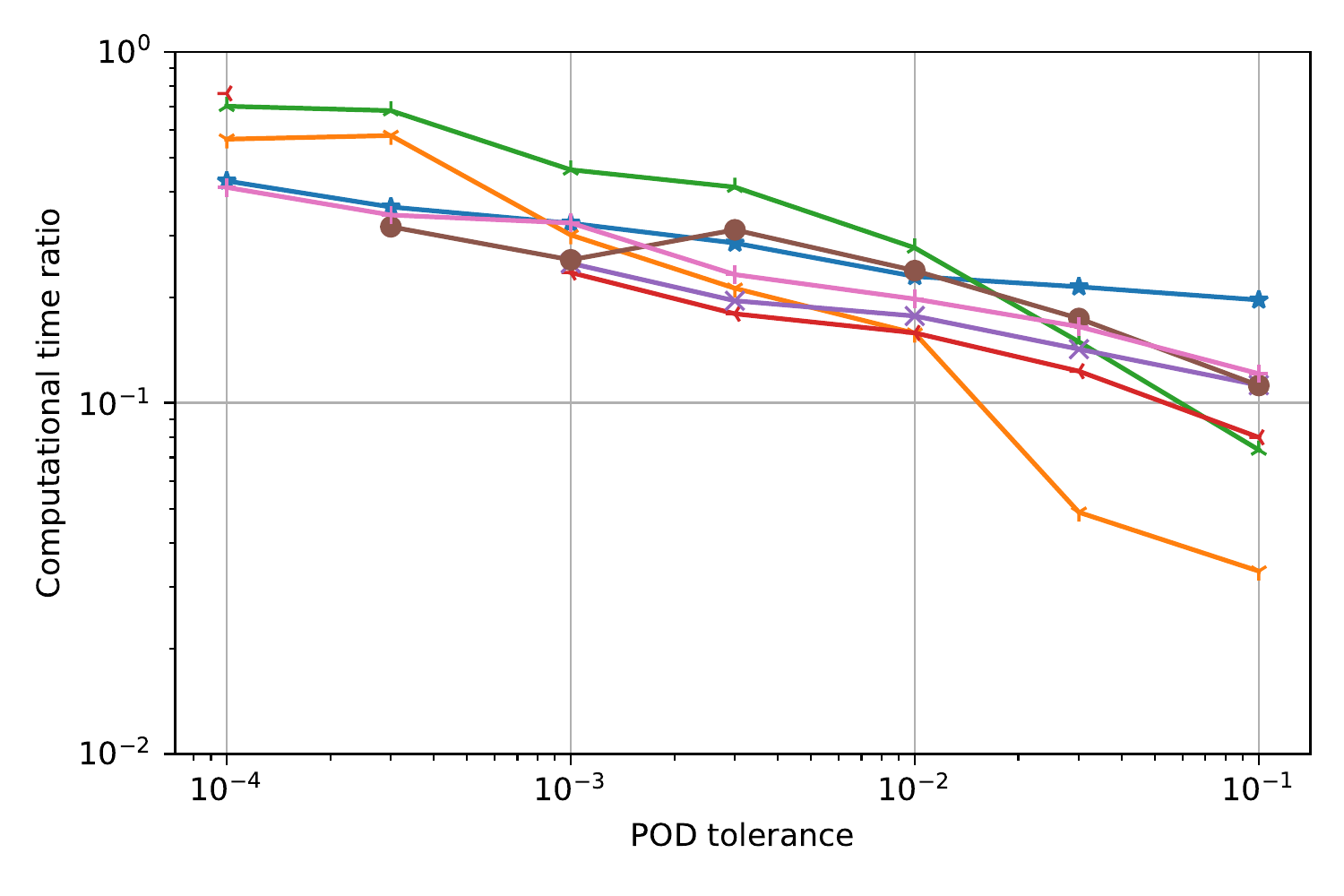}}

		\subfigure[ROM solutions time-only reduction \label{fig:ROMsolSolBump}]{\includegraphics[width=0.32\textwidth]{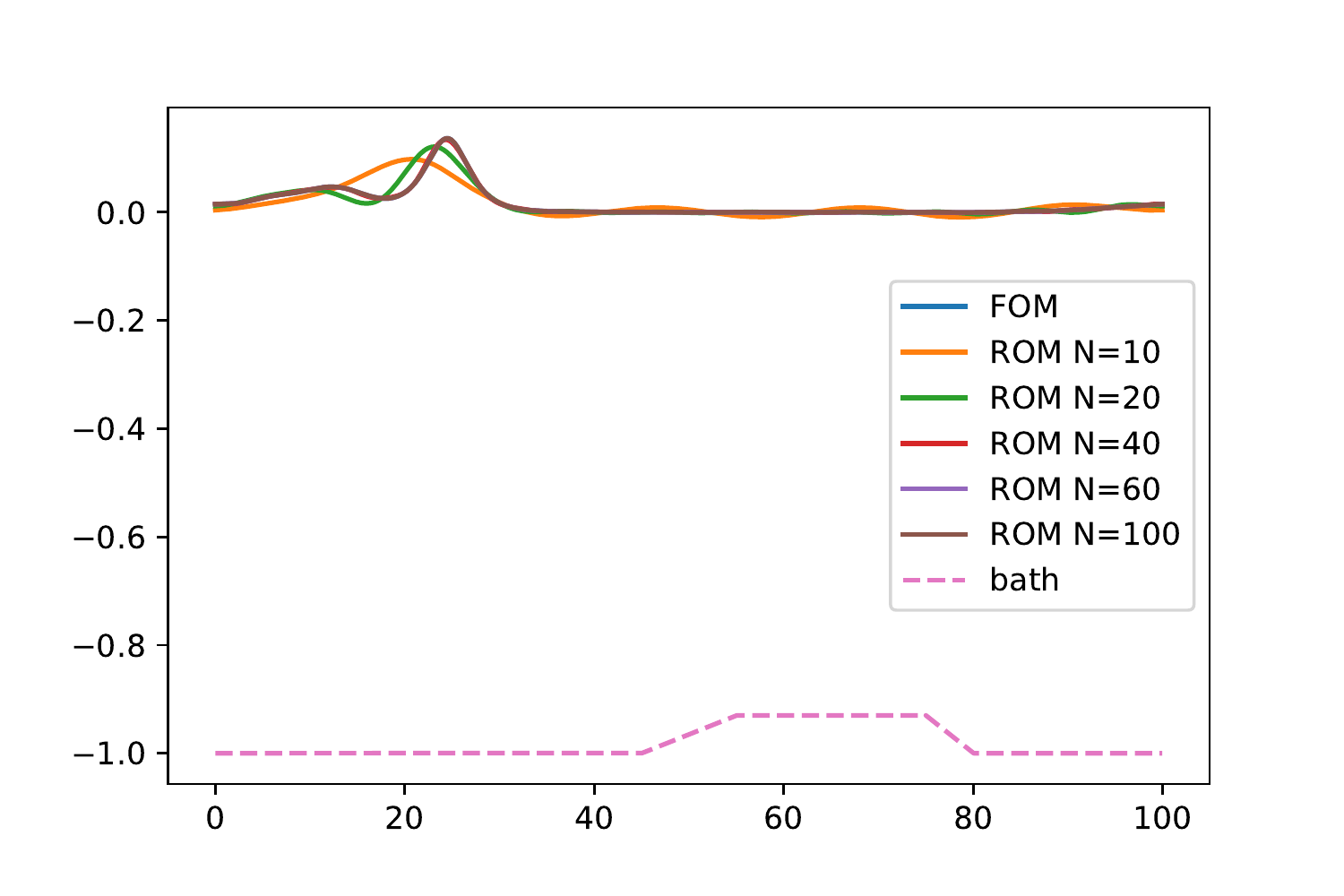}}\;
		\subfigure[ROM solutions time-parameter reduction with $tol_{POD}=0.01$, $\NRB=72$ and $tol_{EIM}=0.03$, $\NEIM=152$ on parameter and $\eta_0$ in the training set \label{fig:ROMEIMparamSimSolBumpin}]{\includegraphics[width=0.32\textwidth]{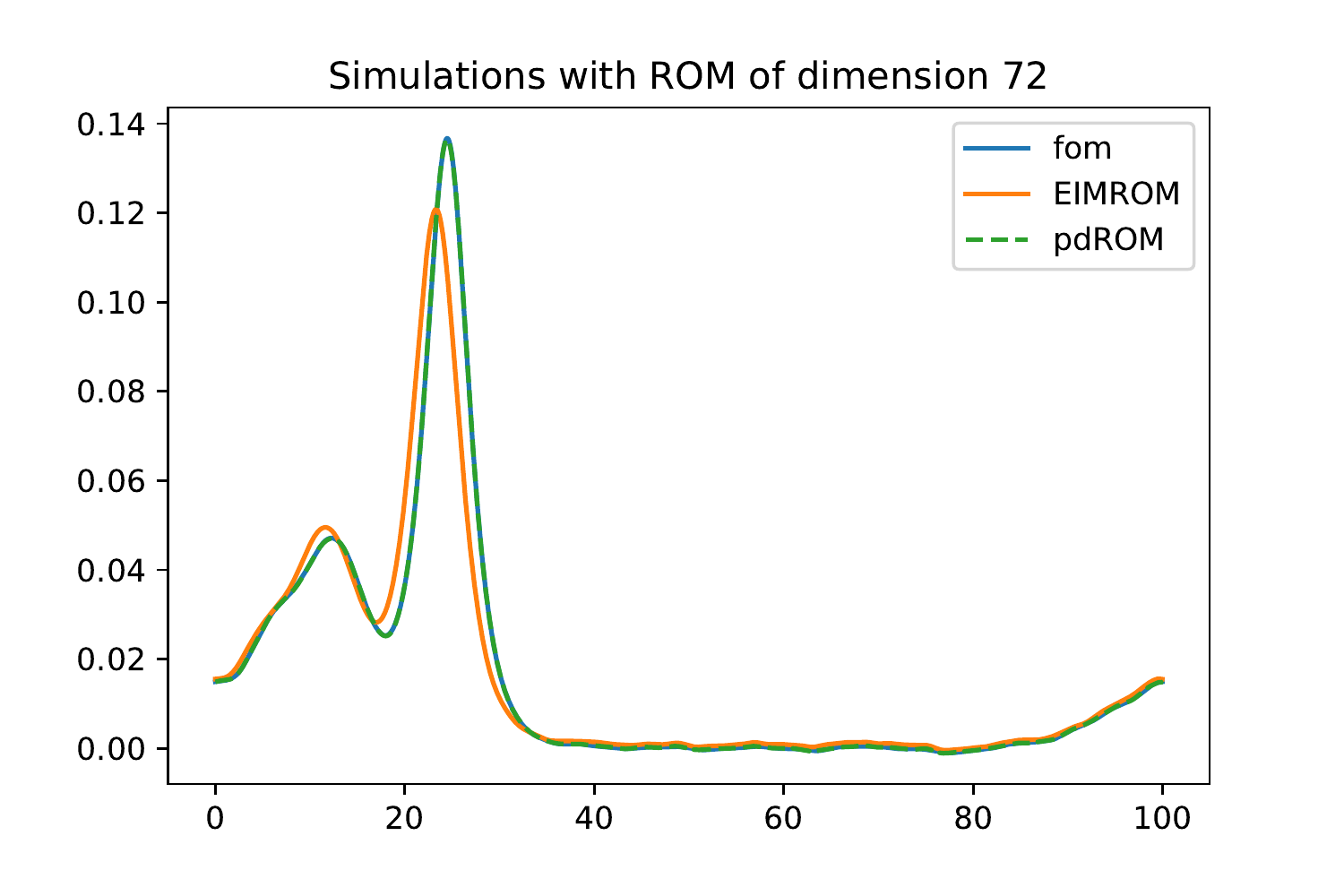}}\;
		\subfigure[ROM solutions time-parameter reduction with $tol_{POD}=0.01$, $\NRB=72$ and $tol_{EIM}=0.03$, $\NEIM=152$ on parameter and $\eta_0$ outside the training set \label{fig:ROMEIMparamSimSolBumpout}]
		{	\includegraphics[width=0.32\textwidth]{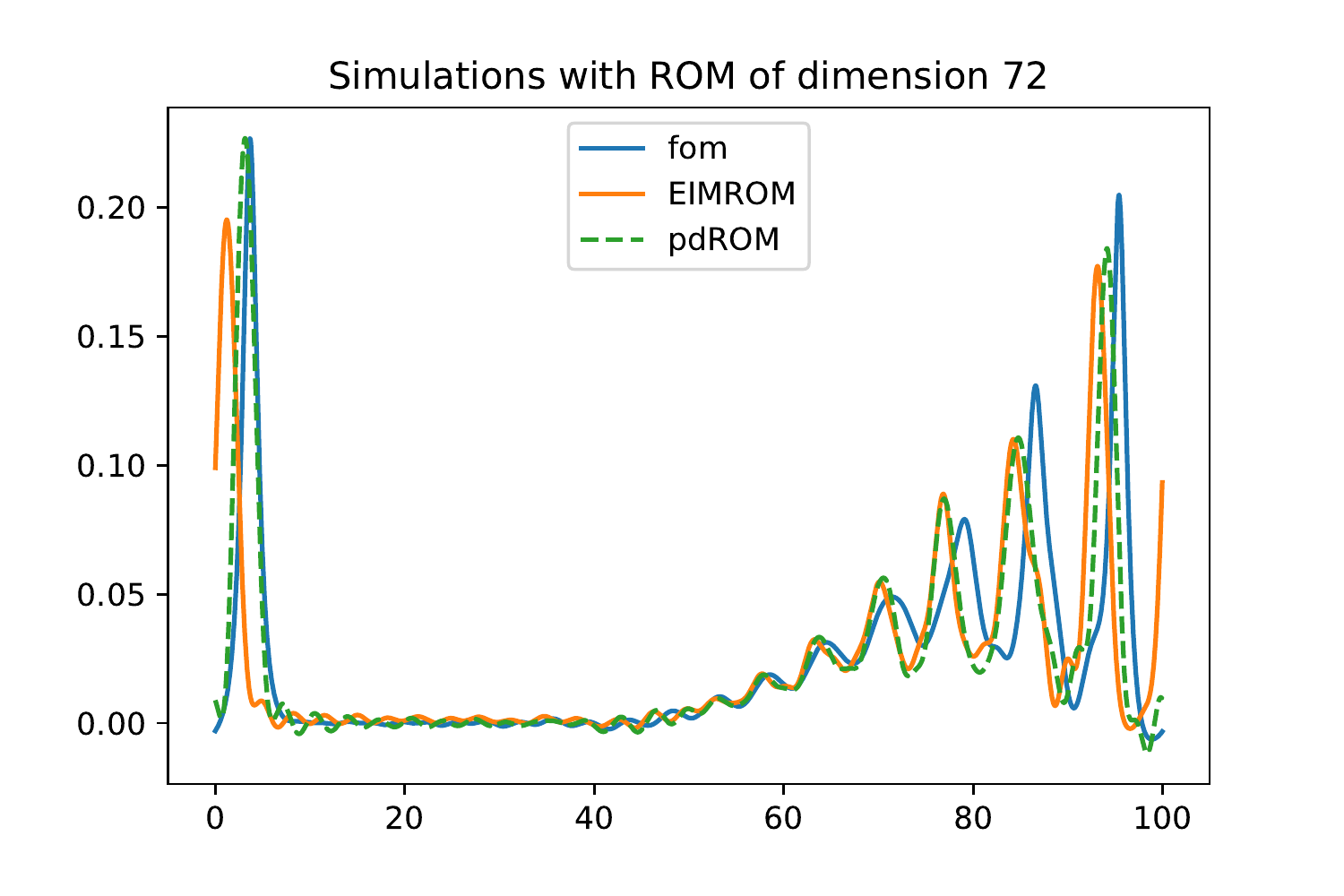}}\;

		\subfigure[EIMROM solutions time-only reduction \label{fig:ROMEIMsolSolBump}]{
			\includegraphics[width=0.32\textwidth]{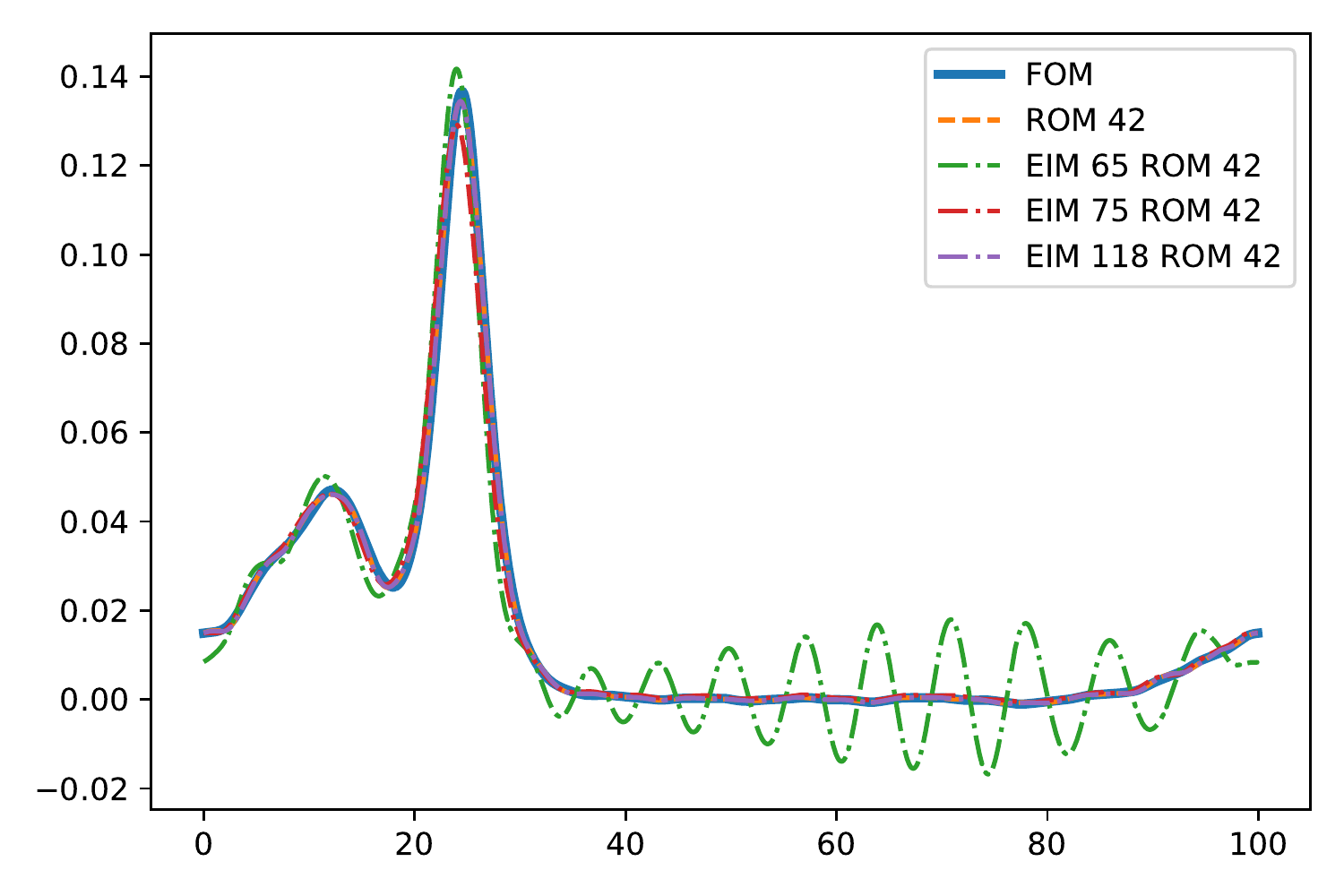}
		}\;
		\subfigure[ROM solutions time-parameter reduction with $tol_{POD}=0.001$, $\NRB=119$ and $tol_{EIM}=0.001$, $\NEIM=259$ on parameter and $\eta_0$ in the training set \label{fig:ROMEIMparamSimSolBumpin2}]{\includegraphics[width=0.32\textwidth]{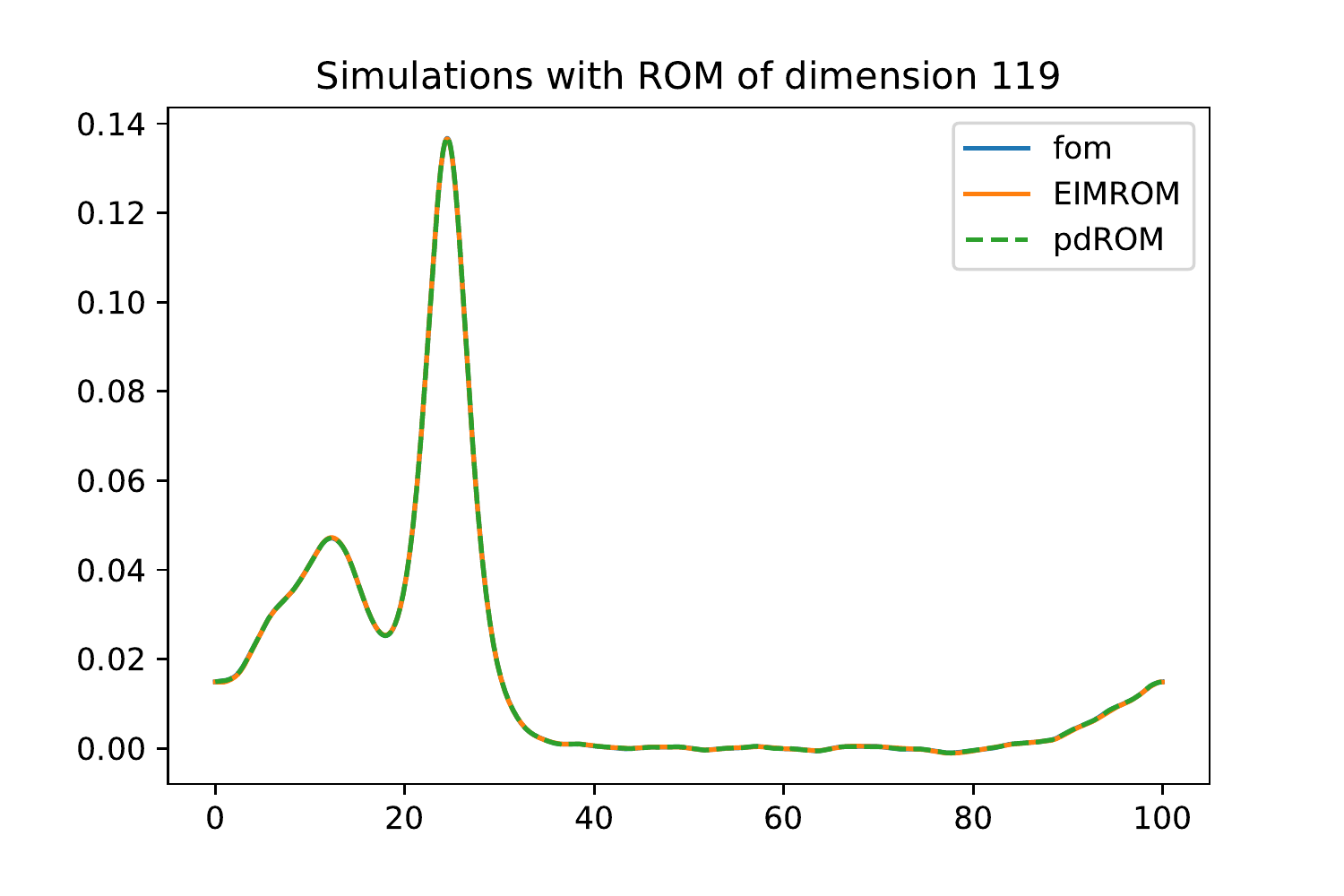}}\;
		\subfigure[ROM solutions time-parameter reduction with $tol_{POD}=0.001$, $\NRB=119$ and $tol_{EIM}=0.001$, $\NEIM=259$ on parameter and $\eta_0$ outside the training set \label{fig:ROMEIMparamSimSolBumpout2}]{\includegraphics[width=0.32\textwidth]{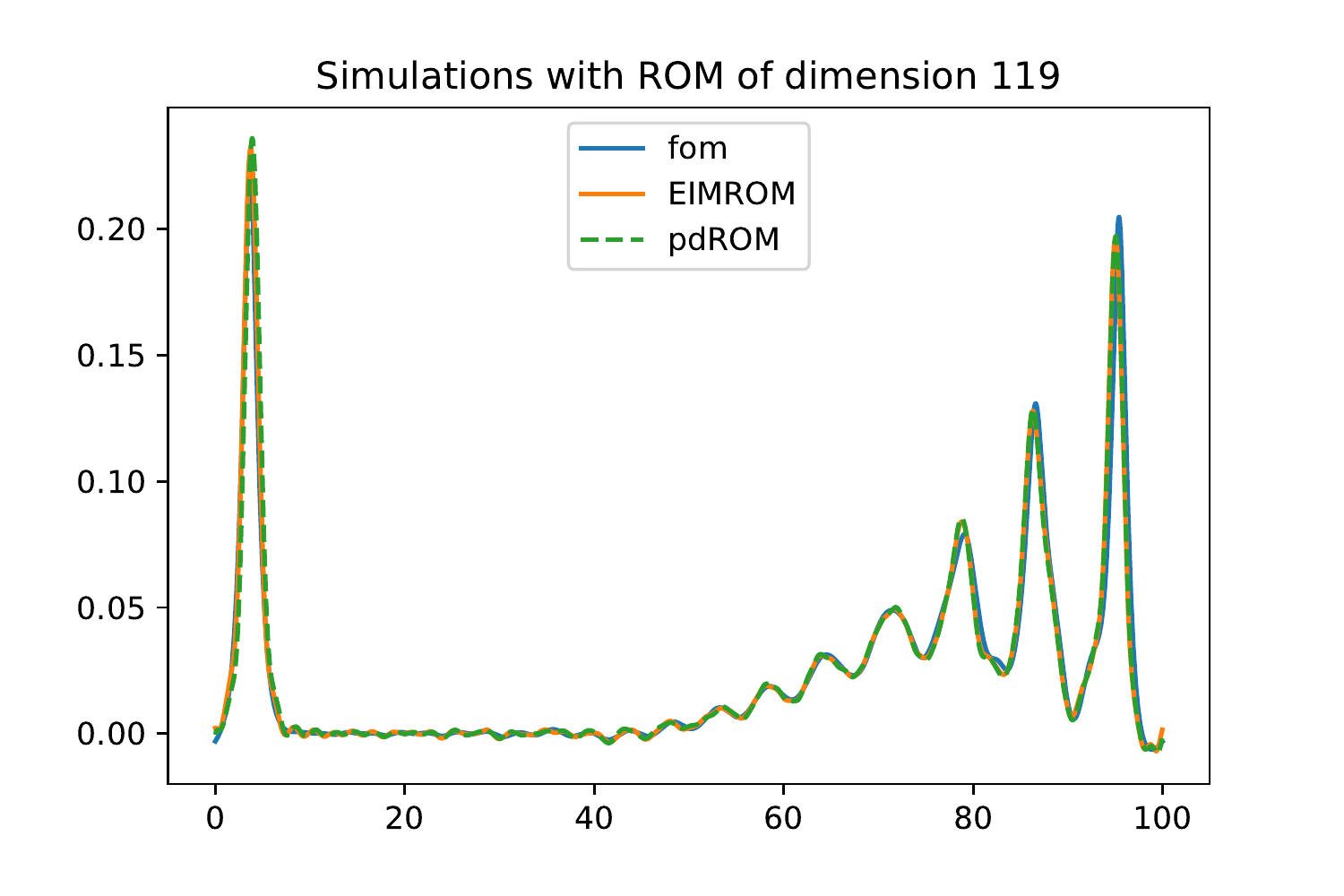}}

	\end{center}
	\caption{\textbf{Solitary waves interacting with a submerged bar.} pdROM and EIM reduction: simulation, errors and computational time, varying POD and EIM dimensions}
\end{figure}
\begin{figure}
	\centering
	\includegraphics[width=0.9\textwidth]{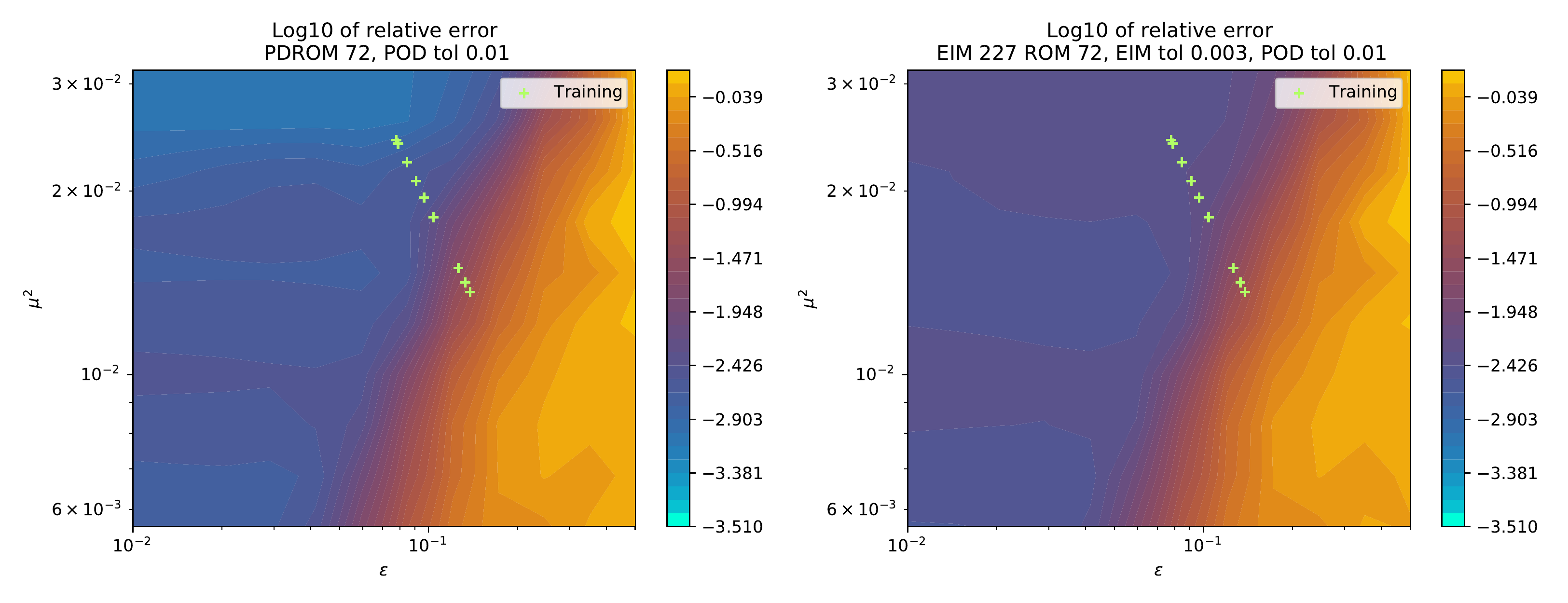}
	\caption{\label{fig:varyingSolBump}\textbf{Solitary waves interacting with a submerged bar.} Error plot varying the parameters $a_0$ and $h_0$ for pdROM $\NRB=72$ and EIMROM with $\NEIM= 227$}
\end{figure}

Clearly this test case is pushing over the limit of the method as the advection dominated part is leading the computations and the contrast between flat and oscillating region leads very easily to Gibbs phenomena. Hence, it is important to stay close to the training set area for this case and not to introduce extra interpolation error with the EIM.

In \cref{fig:varyingDam} we see how the error of the pdROM with $\NRB=46$ varies with respect to $a_0$ and $h_0$. We can see a strong gradient when $h_0$ and $a_0$ increase in the area of unstable solutions. Again, the main factor here
is the relevance of nonlinearity. As long as the waves remain dispersive and we are far from very shallow bores, we observe an error plateau which seems to indicate a certain robustness with respect to the variation
of parameters. 
It must be remarked however, that outside this area, the error grows   relatively fast also for pdROM. Some improved treatment of nonlinearity  is   necessary to handle this kind of problems when getting close the the shallow/highly nonlinear range. 
\subsection{Solitary waves interacting with a submerged bar}  

For this test, we use $\bar{h}_0=1$, final time $T=60$ and $a_0=0.1$.
This test has similar results to the \textbf{Propagation of a periodic monochromatic wave}. Its advection dominated character is more pronounced than in the \textbf{Propagation of a periodic monochromatic wave}, but the error decays quickly enough as we can see in \cref{fig:ROMEIMerrSolBump}. 
The pdROM solution with $\NRB=40$ basis function is indistinguishable from the FOM one in \cref{fig:ROMsolSolBump} and, adding the EIM approximation, we can see different results. In \cref{fig:ROMEIMsolSolBump} for the EIMROM simulation with $\NRB=42$ and $\NEIM=65$ we still do not have enough accuracy and the solution oscillates, with $\NEIM= 75$ the solution is much more stable, but we can still see differences in the peak of the wave, finally for $\NEIM=118$ observe a very accurate solution.
Similarly to before, the computational costs scale between 15\% and 30\% of FOM for pdROM and between 5\% and 25\% for EIMROM, see \cref{fig:ROMEIMtimeSolBump}.

In the second phase, we consider also here 10 simulations for randomly chosen $h_0 \in [0.7\bar{h}_0,1.3\bar{h}_0]$ and we use them to compute the reduced space with the POD. The error decay with respect to the POD tolerance in \cref{fig:ROMEIMparamErrSolBumpin} is similar to the first phase, but the dimension of the spaces are larger as we notice in the computational time plot in \cref{fig:ROMEIMparamTimeSolBumpin} up to a 40\% of the FOM computational time for pdROM. In \cref{fig:ROMEIMparamSimSolBumpin,fig:ROMEIMparamSimSolBumpin2} we see accurate pdROM simulations with $\NRB=72$ for a POD tolerance of $0.01$ and a computational time of around 25\% of the FOM time, while we need $\NEIM=259$ to get very accurate EIMROM solutions for computational costs around 20\% of the FOM ones. Again, the EIM space should be very large in order to obtain accurate and stable results \cite{peherstorfer2020stability,zimmermann2018geometric,ghavamian2017pod,argaud2017stabilization,chen2021eim}.

When testing on a parameter outside the training domain, i.e., $a_0=0.15$ and $h_0=0.63$, we have that the exact solution is more oscillatory, see \cref{fig:ROMEIMparamSimSolBumpout,fig:ROMEIMparamSimSolBumpout2}, and to catch properly the solution we need more basis functions in both spaces. Indeed, the errors are larger in \cref{fig:ROMEIMparamErrSolBumpout}, while the computational times stays similar.

In \cref{fig:varyingSolBump} we see how the error varies for pdROM with $\NRB=72$ and EIMROM with $\NEIM= 227$ for different parameters $h_0$ and $a_0$, plotting the error with respect to $\varepsilon$ and $\mu$. 
As before we see that nonlinearity really plays a major role, with the shallow water/hyperbolic limit being poorly handled by both approaches. 
Away from this limit (left part of the figures)  the EIMROM error is  larger than the pdROM one.

We remark that the computational cost of the linear system is the largest one in the FOM and, while doing a pdROM, we obtain already a great reduction. In \cref{app:cost} we compare the pdROM costs with the reduction of only the linear system for $\Phi$ (evolving the full FOM) as suggested in \cref{rem:modification_sw}. Essentially, the cost of the two methods are identical.
Moreover, changing the linear solver into something more generic results in an even stronger reduction of the computational costs. We compare the costs of the Thomas algorithm and of a generic sparse solver in \cref{app:cost}.

\section{An extension to a    Boussinesq system}\label{sec:EB}
We propose in this section a possible extension to well known dispersive shallow water model by 
Madsen and Sørensen \cite{madsen1992newForm}.
In the dimensional  form the model equations read \cite{madsen1992newForm,ricchiuto2014upwind} 
\begin{equation}\label{eq:MS}
	\begin{split}
		&\begin{cases}
			\partial_t \eta + \partial_x q=0,\\
			\partial_t q -\mathcal{T}^t[\partial_t q] + \partial_x (qu) +g h\partial_x \eta -\mathcal{T}^x[\eta]=0,
		\end{cases}\\
		\text{with }&\mathcal{T}^t [\cdot]:=B_1\bar{h}^2\partial_{xx}[\cdot] + \frac{1}{3}\bar{h}\partial_{x}\bar{h} \partial_{x}[\cdot]\text{ and }\mathcal{T}^x[\cdot]:= gB_2\bar{h}^2 \left(2\partial_{x}\bar{h} + \bar{h}\partial_{x}\right)\partial_{xx}[\cdot],
	\end{split}
\end{equation}
where, with the notation of \cref{fig:sw-notation}, $\eta(x,t)$ is the water level variation from the reference value, $q(x,t)$ is the discharge, $b_0b(x)$ is the bathymetry profile, $h_0$ is the reference value for water height, $\bar{h}(x)=h_0-b_0b(x)$ is the water height at rest, 
$h(x,t)=\eta(x,t)+\bar{h}(x)$ is the water height, $u(x,t)$ is the depth averaged speed and it is linked to $q$ by $q(x,t)=h(x,t)u(x,t)$.  
System \eqref{eq:MS}  provides a $\mathcal{O}(\varepsilon \mu^2, \mu^4)$ approximation of  the full nonlinear wave equations.
It contains the nonlinear terms belonging also to shallow water equations, and they are the only nonlinear terms of the model. 
Similarly to Boussinesq equation and to the previous scalar model there are other \textit{linear} dispersion terms which scale with constants $B_1$ and $B_2$. 
The model constants   $B_1$ and $B_2$ are set to  $B_2 = 1/15$ and $B_1=B_2+1/3$ to optimize the linear characteristics of the model  (phase and shoaling)
with respect to the full Euler equations \cite{madsen1992newForm}. \\

As before, we manipulate the equations for the purpose of their numerical solution. First we recast the model in dimensionless form:
\begin{equation}\label{eq:MSadim}
	\begin{split}
		&\begin{cases}
			\partial_t \eta + \partial_x q=0,\\
			\partial_t q - \mu^2 \mathcal{T}^t[\partial_t q] + \varepsilon \partial_x (qu) + h\partial_x \eta - \mu^2 \mathcal{T}^x [\eta]=0,
		\end{cases}	\\
		\text{with }&\mathcal{T}^t [\cdot]:= B_1\bar{h}^2\partial_{xx}[\cdot] + \frac{1}{3}\bar{h}\partial_{x}\bar{h} \partial_{x}[\cdot]\text{ and }\mathcal{T}^x[\cdot]:= B_2\bar{h}^2 \left(2\partial_{x}\bar{h}+ \bar{h}\partial_{x}\right)\partial_{xx}[\cdot].
	\end{split}
\end{equation}
We then isolate the nonlinear hyperbolic operator and write the system as 
\begin{align}
	&\begin{cases}
		\partial_t \eta + \partial_x q=0,\\
		(1-\mu^2\mathcal{T}^t)[\partial_t q  + \varepsilon \partial_x (qu) + h\partial_x \eta ] + \mu^2 \mathcal{T}^t[ \varepsilon \partial_x (qu) + h\partial_x \eta] -\mu^2 \mathcal{T}^x [\eta]=0.
	\end{cases}	
\end{align}
In  this case we do not neglect any of the small order terms to keep the original model,    often used in literature (see e.g.   \cite{b13,ricchiuto2014upwind,brb20}) and references therein.
The final model can be written as a system of 3 equations, given by
\begin{align}\label{eq:MSadimSW}
	&\begin{cases}
		\partial_t \eta + \partial_x q=0,\\
		(1-\mu^2\mathcal{T}^t)[\psi] =\mu^2\left\{  \mathcal{T}^x [\eta] -  \mathcal{T}^t[ \varepsilon \partial_x (qu) + h\partial_x \eta] \right\},\\
		\partial_t q  + \varepsilon \partial_x (qu) + h\partial_x \eta =\psi.
	\end{cases}	
\end{align}
Also for this model, we show the analogy with the model defined in \eqref{eq:generalModel}. Indeed, if we consider 
\begin{align*}
	&u:=(\eta, q)^T, \quad \partial_{x}F(u) := \left( \partial_x q, \varepsilon \partial_x (qu) +h \partial_x \eta \right)^T,\quad S(u)=0,\\ 
	&\mathcal{X}^x(u) := \left(0,-\mu^2    \left\{  \mathcal{T}^x [\eta] -  \mathcal{T}^t[ \varepsilon \partial_x (qu) + h\partial_x \eta] \right\} \right) ^T,\quad
	\mathcal{X}^t(u):= \left( 0, \mu^2 \mathcal{T}^t \right)^T,
\end{align*}
one obtains \eqref{eq:MSadimSW}.
For this model, one cannot show an exact conservation for a total energy.  For other dispersive nonlinear systems one can show the conservation of an energy of the form 
$\mathcal{E}=\frac{g h^2}{2} + \frac{hu^2}{2} + \frac{h^3(u_x)^2}{6}$  \cite{lannes2020modeling}. Differently from the BBM-KdV equation,
the relation between this energy and the main unknowns is nonquadratic, and it is not straightforward to use it to build the reduced basis and the test reduced space as the one studied before.
Based on regularity requirements, we will instead investigate a formal extension of the method proposed in the scalar case. 

\subsection{Wave generation and boundary conditions}

We will consider two families of solutions: solitary waves, monochromatic periodic waves.
Exact solitary waves for \eqref{eq:MS}   can be obtained assuming $\eta=\eta(x-Ct)$, and $q=C\eta$, which allows
to obtain  an  ODE  from the second equation in \eqref{eq:MS}. The solution is uniquely defined  given  amplitude and depth at rest  ($a_0$, $h_0$).
In particular,  setting $C_0^2=gh_0$, the celerity $C$ verifies the consistency condition 
\begin{equation}
	\frac{C^2}{C_0^2}= \frac{a_0^2}{2h_0^2} \frac{1+\frac{a_0}{3h_0}}{\frac{a_0}{h_0}-\log(1+\frac{a_0}{h_0})}
\end{equation}
Details can be found in  \cite{ricchiuto2014upwind,orszaghova2012paddle}. The ODE is integrated here  with  a built-in solver in {\tt scipy}.\\

To generate monochromatic waves, we have included an internal wave generator. 
Different definitions exist in literature  \cite{tonelli2009hybrid,nwogu1993alternative,walkley1999numerical,wei1999generation}. Following \cite{walkley1999numerical,wei1999generation}, we add to the $\eta$ equation a compact forcing term with periodic variation in time:
\begin{equation}
	\partial_t (\eta+h_{\iwg})+\partial_x q =0,
\end{equation} 
where the internal generator is defined by 
\begin{equation}
	h_\iwg (x,t):= A_\iwg f_\iwg(x)  \sin (\omega t), 
\end{equation}
with $\omega=2\pi/T$ the frequency of the signal and $T$ the period, $f_\iwg$ the spatial damping function
\begin{equation}
	f_\iwg(x) = \Gamma_\iwg e^{-(x-x_\iwg)^2/d_\iwg^2}.
\end{equation}
We refer to \cite{ricchiuto2014upwind,wei1999generation} for the tuning of the parameters and we will use the following ones
\begin{equation}
	\Gamma_\iwg = 0.185 \sqrt{g/h_0} T, \quad d_\iwg^2 = \frac{g \alpha_\iwg^2 h_0 T^2}{80},\quad T=2.525,\quad \alpha_\iwg=4,
\end{equation}
and $x_\iwg$ the center of the generator and $\alpha_\iwg\in [0,4]$ is a parameter controlling the amplitude of the wave and $\beta_\iwg$ is just a rescaling factor. 
Also $A_\iwg$ is a rescaling factor that will be changed along the simulations to test the problem for different parameters.\\

%
%

Finally,  outflow conditions are imposed through \textit{sponge layers}. Following   \cite{ricchiuto2014upwind} these terms are written as   diffusion operators which will be discretized implicitly and added to the mass matrix. At
the continuous level they have the form  
\begin{equation}
	S(v)=-\nu \partial_{xx} v,\qquad \nu(x) := \begin{cases}
		n_1\frac{1-e^{n_2\frac{x-X_{s1}}{X_{s2}-X_{s1}}}}{1-e}, &  X_{s1}\leq x \leq X_{s2},\\
		0,& \text{else.}
	\end{cases}
\end{equation}
Here $n_1$ and $n_2$ are two coefficients that should be tuned with the problem, we will use the values $n_1= 10^{-3}$ and $n_2=10$ as suggested in \cite{ricchiuto2014upwind}. The points $X_{s1}$ and $X_{s2}$ are fixed close to the boundary and $X_{s2}$ will coincide with the boundary itself. 
Also here, both these linear source terms can be added in the formulation to match \eqref{eq:generalModel} defining $\mathbb{S}^{im}$ accordingly.

\subsection{Full order discrete model}

The approach used to discretize \eqref{eq:MS} is similar to the one used in the scalar case. We combine   a linear continuous Galerkin method (see \cite{ricchiuto2014upwind}  and references therein) with an SSPRK(2,2) time integrator \cite{shu-gottlieb-1998}. 
To introduce a  little dissipation we again make use of a continuous interior penalty  operator (CIP)  \cite{burman2004edge}. 
The resulting discrete equations can be written as 
\begin{subequations}\label{eq:EBsystemFOM_SW}
	\begin{align}
		&\mass\left(\!( \boldeta+h_\iwg)^{(s)}\! -\sum_{r=0}^{s-1}\rho_r^s (\boldeta +h_\iwg)^{(r)} \right)\!+ \Delta t \sponge\boldeta^{(s)}  \!=\!- \Delta t\sum_{r=0}^{s-1}\theta_r^s\left( \mathcal{N}^\eta(  \boldeta^{(r)},\bq^{(r)}) +\diff(\boldeta^{(r)},\lambda^{(r)})\right) \\
		&(\mass - \dispersion)\boldsymbol{\psi}   =-  \sum_{r=0}^{s-1}\theta_r^s  \left( \dispersion\mass^{-1}\mathcal{N}^q(\boldeta^{(r)},\bq^{(r)}) +\mathbb{T}^x \boldeta^{(r)}\right) 	\\
		&\mass\left(\!\bq^{(s)}-\sum_{r=0}^{s-1}\rho_r^s \bq^{(r)} \right)+ \Delta t \sponge\bq^{(s)}  = \Delta t\mass \boldsymbol{\psi} -  \Delta t\sum_{r=0}^{s-1}\theta_r^s \left( \mathcal{N}^q(\boldeta^{(r)},\bq^{(r)})+ 	\diff(\bq^{(r)},\lambda^{(r)})  \right) 
	\end{align}
\end{subequations}	
%
where the index in bracket notation is  to denote the Runge--Kutta stages, i.e., $\{\eta^{(s)}, q^{(s)}\}_{s=1,2}$, and $\rho_r^s,\,\theta_r^s$ are the Runge--Kutta coefficients in the Shu--Osher formulation as in the previous sections. 
The matrix $\mass \in \R^{N\times N}$ is the usual $\mathbb P^1$ mass matrix $ \mass_{ij}= \Delta x ( 1 + 3\delta_{ij})/6 $, while the nonlinear operators $\mathcal{N}^{\eta}$ and $\mathcal{N}^q$ are given by
$$
\left[\begin{array}{c}
	\mathcal{N}^{\eta} \\ \mathcal{N}^q
\end{array}
\right]_j= \Delta x \deriv\left[\begin{array}{c}
	q \\ q^2/h
\end{array}
\right]_j+ \dfrac{g}{6}  \left[\begin{array}{c}
	0\\  (2h_j+h_{j+1})(\eta_{j+1}-\eta_j)- (2h_j+h_{j-1})(\eta_{j}-\eta_{j-1})
\end{array}
\right].
$$
The CIP stabilization $\diff$ is defined exactly as in \eqref{eq:CIP},  with $\lambda_j=u_j +\sqrt{gh_j}$. The operators  $\dispersion$ and $\mathbb{T}^x$ are the finite element approximation of the corresponding $\mathcal{T}$ operators appearing in \eqref{eq:MS},
while $\sponge$ is the linear diffusion operator obtained from the discretization of the sponge layer terms. These are reported in \ref{sec:FEM_MS} for completeness.

Also for this discrete model we can see the matching with \eqref{eq:generalModelFOM}, with few differences in the discretization. The mass matrices can be rewritten as $\mathbb{M}^u = \mathbb{M}^\Phi :=  \begin{pmatrix}
	\mathbb{M} & 0\\ 0 & \mathbb{M}
\end{pmatrix},$
the dispersion terms and the fluxes can be defined as 
\begin{equation}
	\mathbb X^t := \begin{pmatrix}
		0&0\\0&\mathbb T^t
	\end{pmatrix}, \quad \mathbb X^x u := \begin{pmatrix}
		0,& -\mu^2 \left\lbrace \mathbb{T}^x \eta - \mathbb{T}^t \mathbb{M}^{-1} \mathcal{N}^q(\boldeta, \bq) \right \rbrace
	\end{pmatrix}^T, \quad \mathbb{F}(u) := \left[\begin{array}{c}
		\mathcal{N}^{\eta} + \mathcal{J}(\boldeta, \lambda) \\ \mathcal{N}^q + \mathcal{J}(\bq, \lambda)
	\end{array}
	\right]
\end{equation}
and the source terms can be written as 
\begin{equation}
	\mathbb{S}^{im} = \begin{pmatrix}
		\mathbb{S} &0\\0&\mathbb S
	\end{pmatrix}, \qquad 
	\mathbb{S}^{ex}(u) =\left[ \begin{array}{c} \partial_t h_{\iwg}\\0 \end{array} \right] .
\end{equation}
The main difference with the KdV-BBM model is that here the explicit dispersive terms are nonlinear as they contain a flux term.
For more details on the discretization of the spatial operators we refer to \cite{ricchiuto2014upwind}.

\subsection{Benchmarks for weakly dispersive free surface waves}

\begin{figure}
	\begin{center}
		\subfigure[Solitary waves interacting with a submerged bar \label{fig:FOMfinalTest4}]{
			\includegraphics[width=0.4\textwidth, trim={25 0 20 20}, clip]{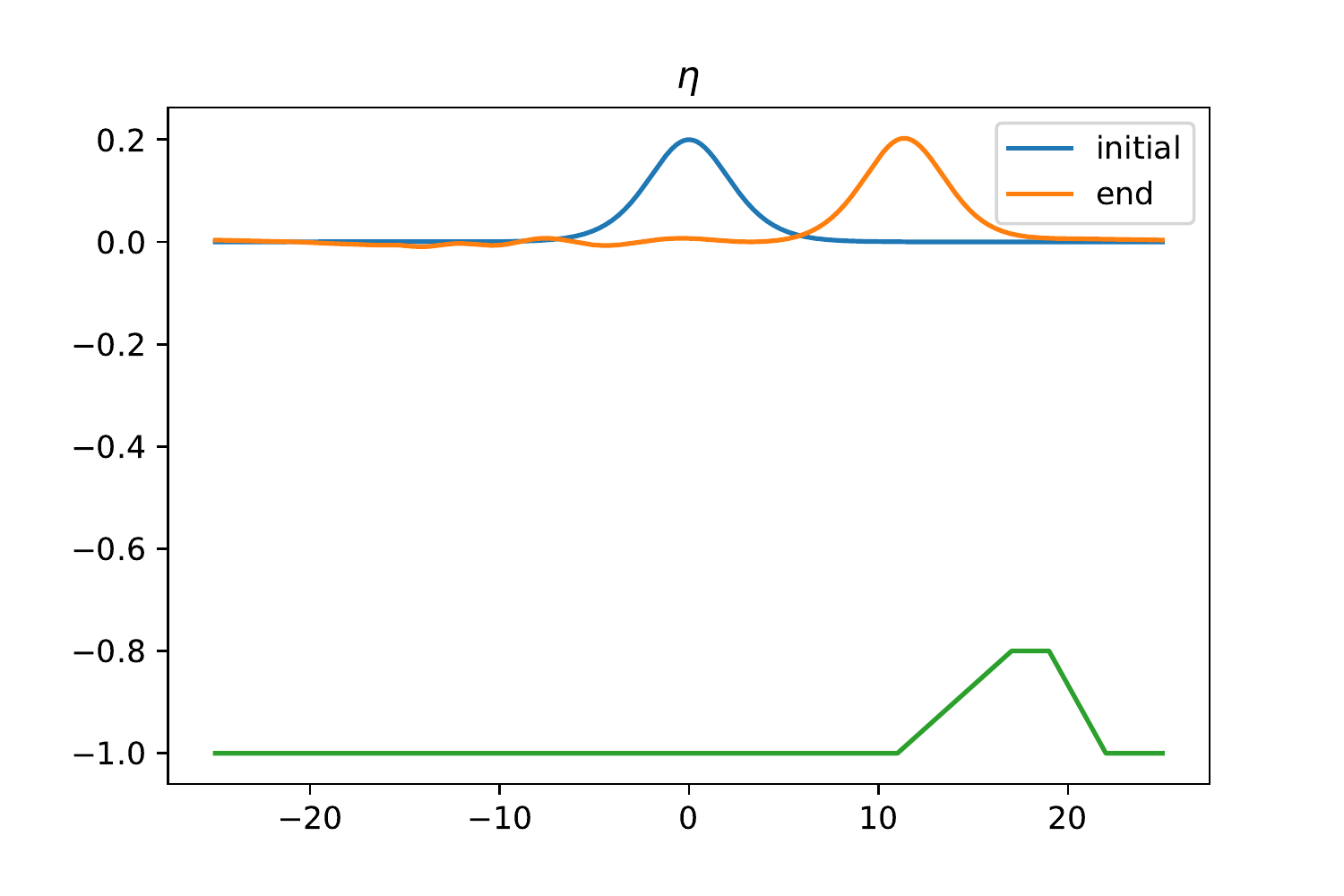}}\;
		\subfigure[Periodic waves over  a submerged bar\label{fig:FOMfinalTest5}]{
			\includegraphics[width=0.4\textwidth, trim={25 0 20 20}, clip]{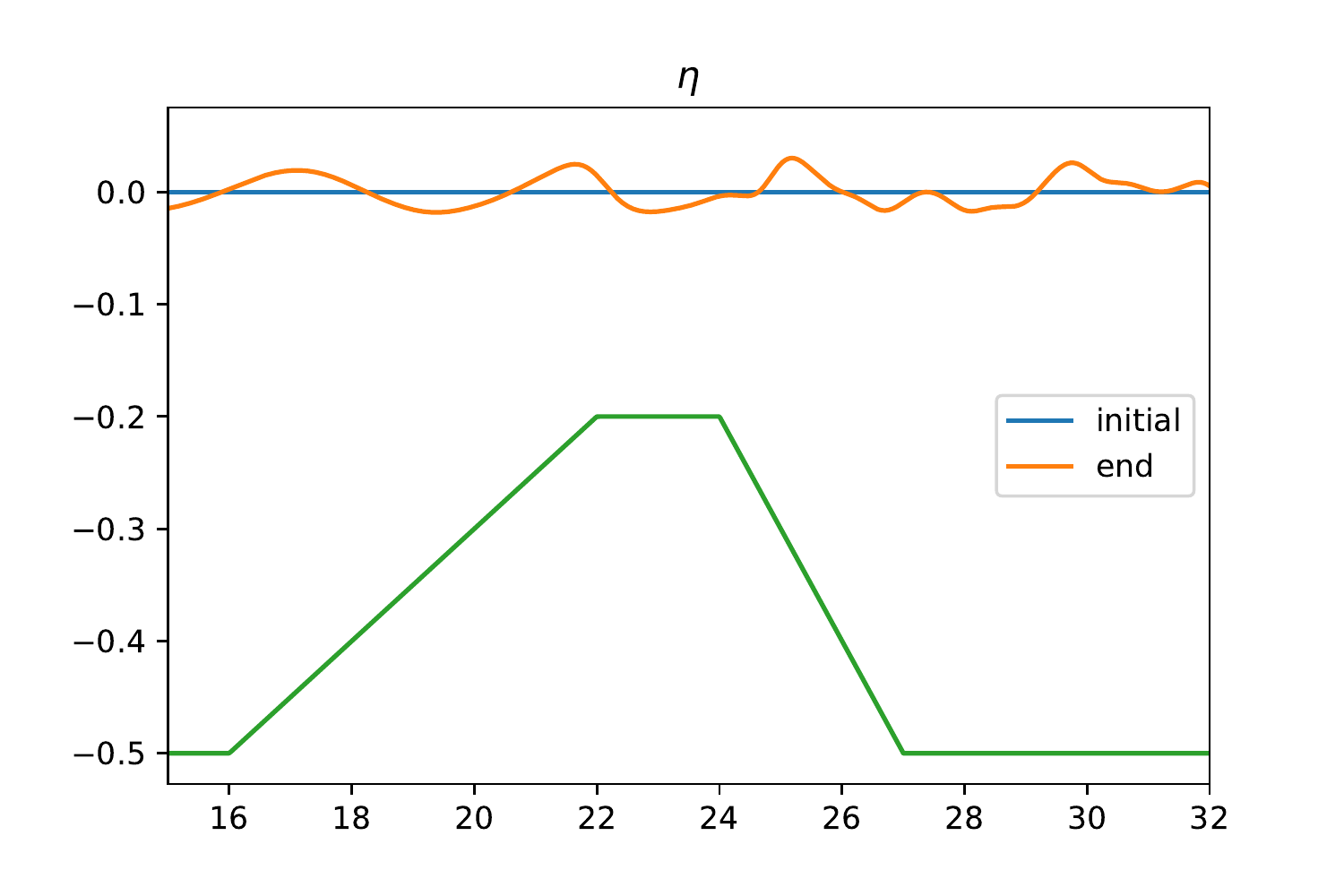}}\;
		\subfigure[Solitary waves interacting with a submerged bar: $\eta$ at different timesteps \label{fig:FOMfinalTest4b}]{
			\includegraphics[width=0.4\textwidth, trim={100 0 550 50}, clip]{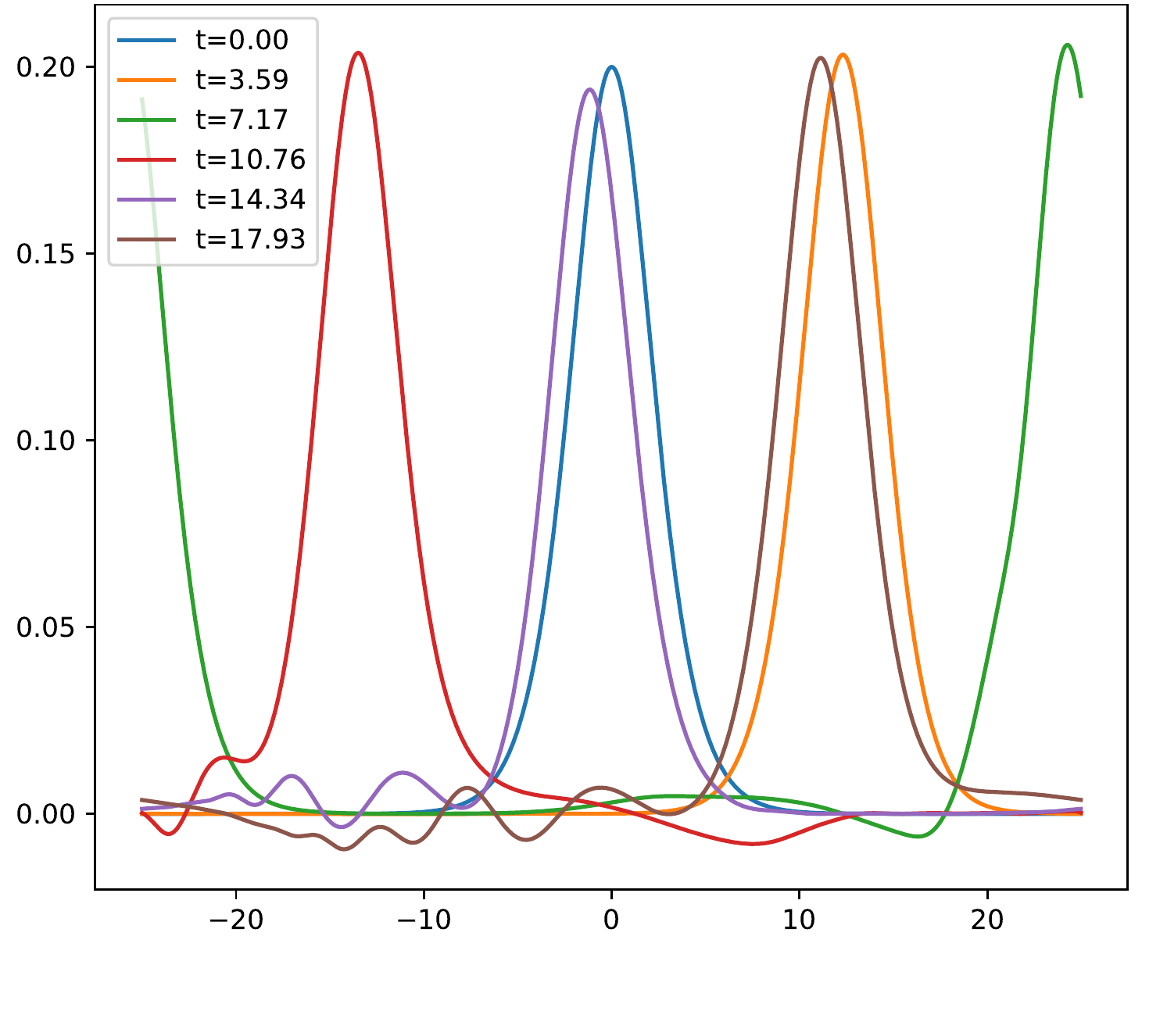}} \;
		\subfigure[Periodic waves over  a submerged bar: $\eta$ at different timesteps \label{fig:FOMfinalTest5b}]{
			\includegraphics[width=0.4\textwidth, trim={100 0 550 50}, clip]{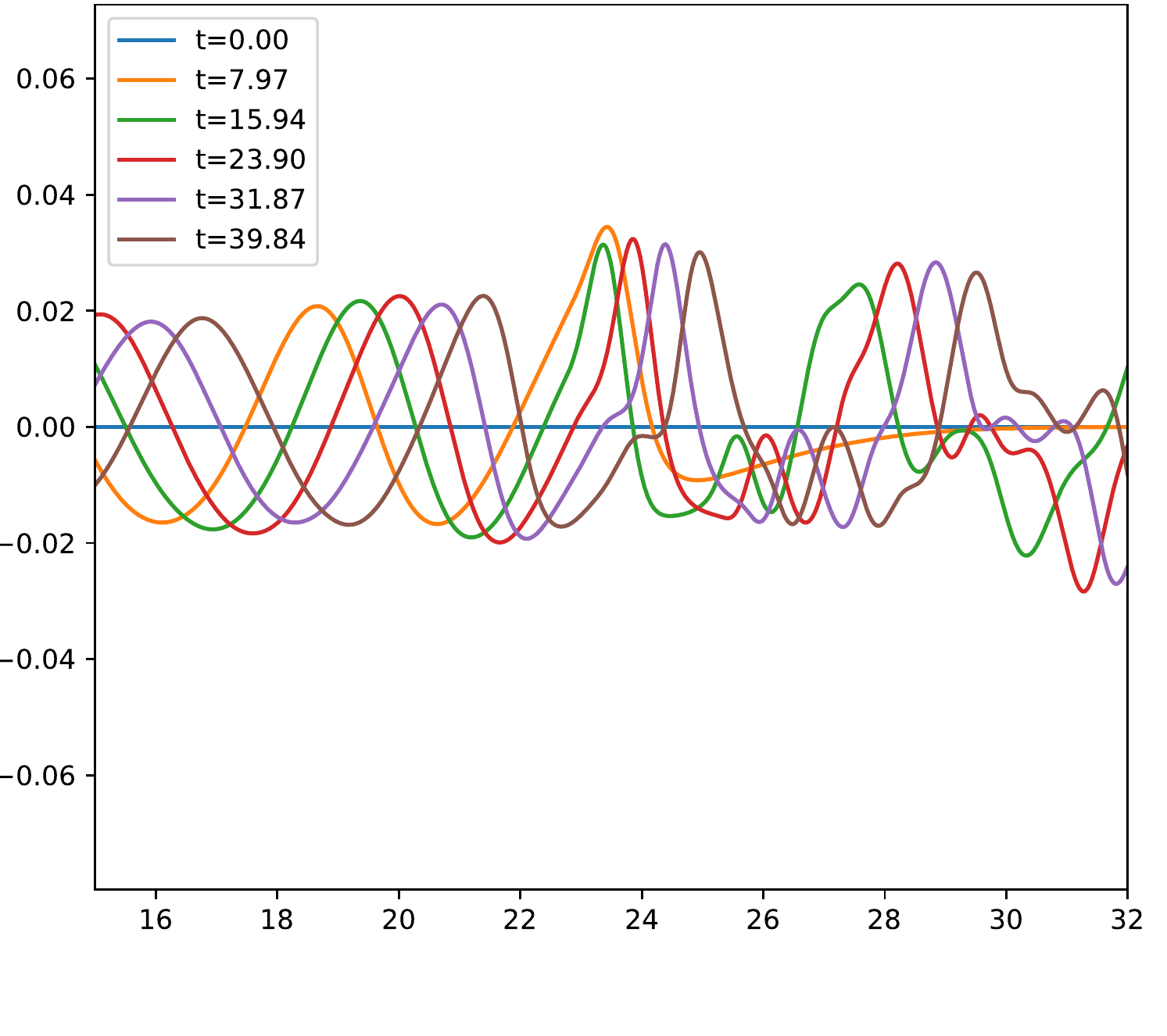}}		
		\caption{Full order model solutions for EB system tests.}\label{fig:FOMsystemEB}
	\end{center}
\end{figure}

\subsubsection{Solitary waves interacting with a submerged bar}\label{sec:EB_solitary}
We consider the simulation of solitary waves  propagating over a trapezoidal bar, defined as in 
\cref{sec:benchSolitaryTrapezoid}, but with  maximal height equal to $b_0=0.2$ and points of change of slope equal to $\lbrace 11, 17, 19, 22\rbrace$. We choose periodic boundary conditions on the domain $\Omega=[-20,30]$ and we center the initial solitary wave at the initial condition in $x=5$. 
Even though this is a traveling solution, the shape of the soliton is smooth  and large enough not to encounter problems in the reduction due to the advection character of the solution.
We use the following parameters $a_0=0.2$, $h_0=1$, $g=9.81$, $\CFL = 0.5$, $T=18$.  Simulations at various times are displayed in \cref{fig:FOMfinalTest4,fig:FOMfinalTest4b}.

\subsubsection{Periodic waves over  a submerged bar}\label{sec:EB_periodic}
We consider monochromatic waves  propagating in the domain $\Omega = [0,35]$ on a trapezoidal submerged bar  of maximum height $0.3$.  An internal wave generator is positioned at $x_{\iwg}=10$. Sponge layers are
set on both ends of the domain. The wave parameters $a_0$ and $h_0$ will be changed along the simulations. For the first test we use $a_0=0.027$, as in \cite{ricchiuto2014upwind} and $h_0=0.5$ (original was $h_0=0.4$ that we will leave for the final simulations). As before, $g=9.81$ and $\CFL = 0.5$, while we set $T=40$. 
FOM simulations for different times are depicted in \cref{fig:FOMfinalTest5,fig:FOMfinalTest5b}. This test is particularly complicated to be reduced, as one can heuristically see from the simulations. Indeed, the gradients get steep after the trapezoid and travels all along the right domain.

\begin{table}
	\centering
	\caption{Estimate percentage of computational time for solving the linear systems for the three benchmark problems: comparison of the implementation of Thomas algorithm all in \texttt{numba} framework and sparse solver by \texttt{scipy} with fluxes computed in \texttt{numba}  \label{tab:eb_system} }
	\begin{tabular}{|c||c|c||c|c|}
		\hline 
		& \multicolumn{2}{|c||}{\texttt{numba }and Thomas} &\multicolumn{2}{|c|}{\texttt{numba }and \texttt{scipy}} \\ \hline
		$N_h$ & \ref{sec:EB_solitary} & \ref{sec:EB_periodic}& \ref{sec:EB_solitary} & \ref{sec:EB_periodic} \\ \hline 
		2000 & 75.8\%& 39.1\%& 74.3\%& 77.7\% \\
		4000 & 72.1\%&  32.8\% &  87.4\%&  80.8\% \\
		8000 & 69.1\%& 27.7\% & 94.5\%& 81.1\% \\
		\hline
	\end{tabular} \vspace{2mm}
	\caption{Estimate of ratio of computational time of the linear system over computational time of the nonlinear fluxes for the three benchmark problems: comparison of the implementation of Thomas algorithm all in \texttt{numba} framework and sparse solver by \texttt{scipy} with fluxes computed in \texttt{numba} \label{tab:eb_flux_system}}
	\begin{tabular}{|c||c|c||c|c|}
		\hline 
		& \multicolumn{2}{|c||}{\texttt{numba }and Thomas} &\multicolumn{2}{|c|}{\texttt{numba }and \texttt{scipy}} \\ \hline
		$N_h$ & \ref{sec:EB_solitary} & \ref{sec:EB_periodic}& \ref{sec:EB_solitary} & \ref{sec:EB_periodic} \\ \hline 
		2000 &6.8& 2.25&39.8& 31.5 \\
		4000 &6.7& 2.18 &41.1 & 31.5\\
		8000 &6.9& 2.15 &39.0 & 32.1\\
		\hline
	\end{tabular}
\end{table}
Here, we compare as for the BBM-KdV benchmarks the computational cost of the system solution with respect to the whole simulation. As before, we observe that the cost of the system is a large part of the simulation. In \cref{tab:kdv_system} we list the time percentage devoted to the solution of the linear system, while in \cref{tab:kdv_flux_system} we show the ratio between the solver time and the flux evaluation.
For the very optimized Thomas algorithm in \texttt{numba} we observe that the linear solver is not so expensive as in the previous cases, in particular for the non periodic boundary condition case.
On the other side, for a generic sparse solver the computational costs of the linear solver are huge with respect to the other operations and it is noticeable comparing it with the cost of the explicit flux evaluation.

\subsection{Reduction of the enhanced Boussinesq model}
For the reduction of the FOM   we proceed analogously to the scalar case.  The main difference is that, in absence of a clear energy statement allowing a linear reduction, each  variable  has its own reduced basis    associated to the minimization of the $L^2$ norm 
as in classical Galerkin projection methods.  We discuss   the reduction main steps  hereafter.

\subsubsection{Projection on reduced spaces and linear reduction}
We consider two reduced spaces for the two variables: $V^\eta \in \R^{\NFE \times \NRB^\eta}$ and $V^q \in \R^{\NFE \times \NRB^q}$. The spaces are obtained again with the POD algorithm, using the same tolerance on both variables. Then, we project the equations \eqref{eq:EBsystemFOM_SW} onto these spaces, obtaining
\begin{equation}\label{eq:EBsystemProj}
	\begin{split}
		W^{\eta ,T} (\mass+ \Delta t \sponge) V^\eta \hat{\boldeta}^{(s)} 	=&W^{\eta ,T} \mass \left(  \sum_{r=0}^{s-1} \rho_r^s (V^\eta \hat{\boldeta}^{(r)} +h_{\iwg}^{(r)}) - h_{\iwg}^{(s)} \right)\\
		-&\Delta t \sum_{r=0}^{s-1} \theta_r^s W^{\eta ,T}\left( \NLRHS^\eta\left( V^\eta \hat{\boldeta}^{(r)}, V^q \hat{\bq}^{(r)}\right)   +\diff(V^\eta \hat{\boldeta}^{(r)},\lambda^{(r)})  \right)\\
		W^{q ,T}(\mass - \dispersion)V^q\hat{\boldsymbol{\psi}}  =  &	-\sum_{r=0}^{s-1}\theta_r^s W^{q ,T}\left( \dispersion \mass^{-1}\mathcal{N}^q(V^\eta\hat{\boldeta}^{(r)},V^q\hat{\bq}^{(r)})+\mathbb{T}^x V^\eta\hat{\boldeta}^{(r)}  \right)  \\
		W^{q ,T}(\mass + \Delta t  \sponge)V^q \hat{\bq}^{(s)}   = &  W^{q ,T}\mass V^q  \sum_{r=0}^{s-1}\rho_r^s \hat{\bq}^{(r)}    +    \Delta tW^{q ,T}\mass^q V^q\hat{\boldsymbol{\psi}} \\ -&
		\Delta t\sum_{r=0}^{s-1}\theta_r^s  W^{q ,T}\left( \mathcal{N}^q(V^\eta\hat{\boldeta}^{(r)},V^q\hat{\bq}^{(r)})  +  \diff(V^q \hat{\bq}^{(r)},\lambda^{(r)})  \right) 
	\end{split}
\end{equation}
for every stage $s$ of SSPRK2. The reduced variables are denoted with the $\hat{\cdot}$ symbol.
As already remarked several times, we have no energy allowing to define a minimization problem leading to   an efficient linear   reduction strategy. Hence we will make use here of classical $L^2$ projection arguments, and  
consider $\Theta$ to be the identity, and $W=V$.\\

Similarly to the scalar case, we can precompute the projected operators and use them in the reduced simulation to save computational time. In this first stage, we obtain the pdROM algorithm. Let us define $\hat{\mass}^\eta :=W^{\eta ,T} \mass V^\eta$, $\hat{\mass}^q :=W^{q ,T} \mass  V^q$, $\hat{\dispersion}:=W^{q,T}\dispersion V^q$, $\hat{\mathbb{T}^x}:=W^{q,T}\mathbb{T}^x V^\eta$, $\hat{\sponge}^\eta := W^{\eta ,T} \sponge  V^\eta$,  $\hat{\sponge}^q := W^{q ,T} \sponge  V^q$, recalling that $h_{\iwg} = \mu_{\iwg}(t, a_0) f_{\iwg}$, we also define $\hat{f}_{\iwg} := W^{\eta ,T}f_{\iwg}$. We obtain substituting in \eqref{eq:EBsystemProj}
\begin{subequations}\label{eq:EBsystemROM}
	\begin{align}
		(\hat{\mass}^\eta + \Delta t \hat{\sponge}^\eta) \hat{\boldeta}^{(s)}= &\hat{\mass}^\eta \sum_{r=0}^{s-1} \rho_r^s \hat{\boldeta}^{(r)}  -  \left(\mu_{\iwg}^{(s)}-\sum_{r=0}^{s-1} \rho_r^s\mu_{\iwg}^{(r)}\right)\hat{\mass}^\eta \hat{f}_{\iwg}  -\Delta t \sum_{r=0}^{s-1} \theta_r^s \hat{\NLRHS}^{\eta,(r)},\\
		(\hat{\mass}^q - \hat{\dispersion} )\hat{\boldsymbol{\psi}}  = &-\sum_{r=0}^{s-1}\theta_r^s  \left( \hat{\dispersion} \hat{\mass}^{q,-1}   \hat{\NLRHS}^{q,(r)} +\hat{\mathbb{T}^x} \hat{\boldeta}^{(r)} \right),   \label{eq:EBROM_psi}\\
		(\hat{\mass}^q + \Delta t \hat{\sponge}^q) \hat{\bq}^{(s)} 
		&= \hat{\mass}^q  \sum_{r=0}^{s-1}\rho_r^s \hat{\bq}^{(r)}   + \Delta t\hat{\mass}^q \hat{\boldsymbol{\psi}} -\Delta t   \sum_{r=0}^{s-1}\theta_r^s    \hat{\NLRHS}^{q,(r)},\label{eq:EBROM_q}
	\end{align}
\end{subequations}
where the reduced nonlinear fluxes are defined as
\begin{equation}
	\begin{split}
		\hat{\NLRHS}^{\eta,(r)}:= &W^{\eta ,T}\left( \NLRHS^\eta\left( V^\eta \hat{\boldeta}^{(r)}, V^q \hat{\bq}^{(r)}\right)   +\diff(V^\eta \hat{\boldeta}^{(r)},\lambda^{(r)})  \right),\\\
		\hat{\NLRHS}^{q,(r)}:=&W^{q ,T}\left( \mathcal{N}^q(V^\eta\hat{\boldeta}^{(r)},V^q\hat{\bq}^{(r)})  +  \diff(V^q \hat{\bq}^{(r)},\lambda^{(r)})  \right) .
	\end{split}
\end{equation} 
Notice that there has been a further projection inside the matrix multiplication 
\begin{equation}
	W^{q,T}\dispersion \mass^{-1} \NLRHS^q \approx W^{q,T}\dispersion V^q V^{q,T} \mass^{-1} W^{q} W^{q,T} \NLRHS^q = \hat{\dispersion} \hat{\mass}^{q,-1} \hat{\NLRHS}^q, 
\end{equation}
in order to have an efficient implementation of the operations without the need of too many operators.
Hence, all the matricial operations are reduced and only the nonlinear fluxes have computational costs that scale as an $\mathcal{O}(\NFE)$.
The solution of the reduced system can be done after the simple assembly of the reduced mass and sponge matrices, with dimension $\NRB$.
As before, this method will be denoted with pdROM in the benchmarks section and all the terms can be matched with \eqref{eq:ROM_general} following the same definitions of the FOM discretization.

\begin{remark}[A more efficient implementation]
	All of the above manipulations retain the structure of a perturbation of the nonlinear shallow  water solver via the use of the auxiliary variable $\boldsymbol{\psi}$.  This is very interesting from the practical point of view for several reasons as in perspective it could be used
	to enhance an existing shallow water solver. It also allows to more easily embed certain wave breaking closures (see e.g. \cite{filippini2016flexible,kr18}).
	However, this reduced model can be more efficiently implemented using a global form that does not isolate the shallow water equations. 
	In particular, instead of \eqref{eq:EBROM_q} and \eqref{eq:EBROM_psi} we could use the following
	\begin{equation}\label{eq:EBsystemROM-fast}
		(\hat{\mass}^q - \hat{\dispersion}  +\Delta t \hat{\sponge}^q  )  \hat{\bq}^{(s)}- (\hat{\mass}^q - \hat{\dispersion} )\sum_{r=0}^{s-1}\rho_r^s \hat{\bq}^{(r)}   + \sum_{r=0}^{s-1}\theta_r^s  \Delta t\left(  \hat{\NLRHS}^{q,(r)} +\hat{\mathbb{T}^x} \hat{\boldeta}^{(r)} \right) =0.
	\end{equation}
	Avoiding the explicit evaluation of  $\boldsymbol{\psi}$ allows to avoid several projections and reduce the online cost.
\end{remark}

\subsubsection{Hyper-reduction}
The final reduction that can be applied is the EIM on the nonlinear fluxes $\NLRHS^\eta,\,\NLRHS^q$. The same procedures presented in the scalar case can be applied in this case, obtaining an extra level of approximation which allows to get rid of all the computations that scale as an $\mathcal{O}(\NFE)$. Same caveats hold for this problem: the tuning of the tolerance of the EIM algorithm must be done accordingly to the POD one otherwise instabilities and Runge phenomena might appear, see \cite{peherstorfer2020stability,zimmermann2018geometric,ghavamian2017pod,argaud2017stabilization,chen2021eim}.

\section{Simulations for enhanced Boussinesq}\label{sec:simulationsEB}
In this section we study the simulations for the EB model. We follow the \textit{modus operandi} of \cref{sec:modusOperandi} We will show results mainly for the variable $\eta$, omitting the ones for $q$ for brevity. 
In the second stage we will let not only $h_0$ vary to form the training set, but also $a_0$, we will specify their range in the different tests.
To make a comparison with a reasonable solver, being the problem more complicated, we will use a sparse \texttt{scipy} solver for the FOM, where the tridiagonal structure of the matrices is not exploited \textit{a priori}.

\subsection{Solitary waves over a submerged bar}   
\renewcommand{\folderTest}{Test_systemSolitonPeriodic}
\begin{figure}
	\begin{center}
		\subfigure[ROM solutions time-only reduction \label{fig:ROMsolSysSolPer}]{\includegraphics[width=0.32\textwidth,trim={20 0 460 20},clip]{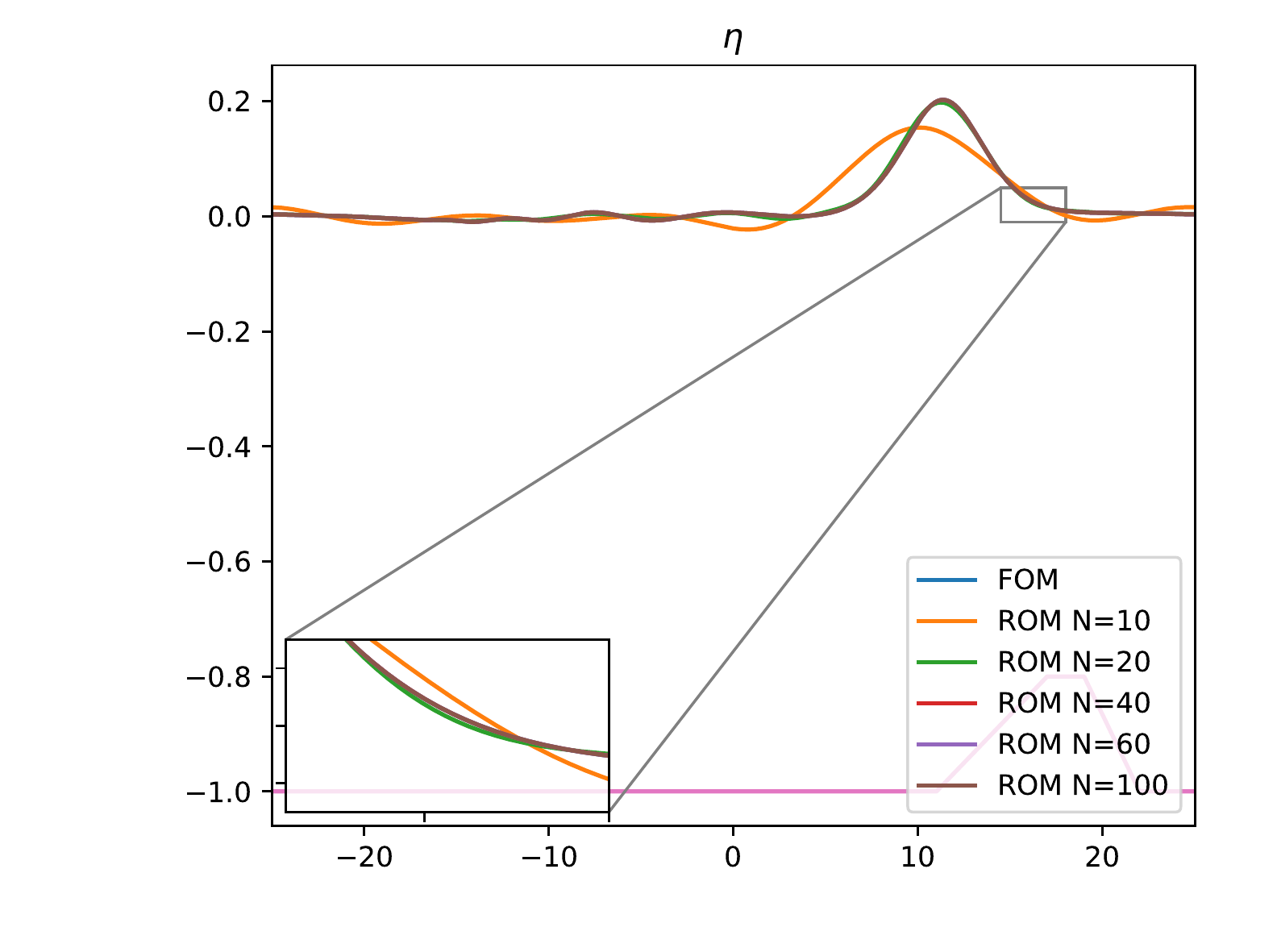}}\;
		\subfigure[EIM ROM errors time-only reduction \label{fig:ROMEIMerrSysSolPer}]{\includegraphics[width=0.32\textwidth]{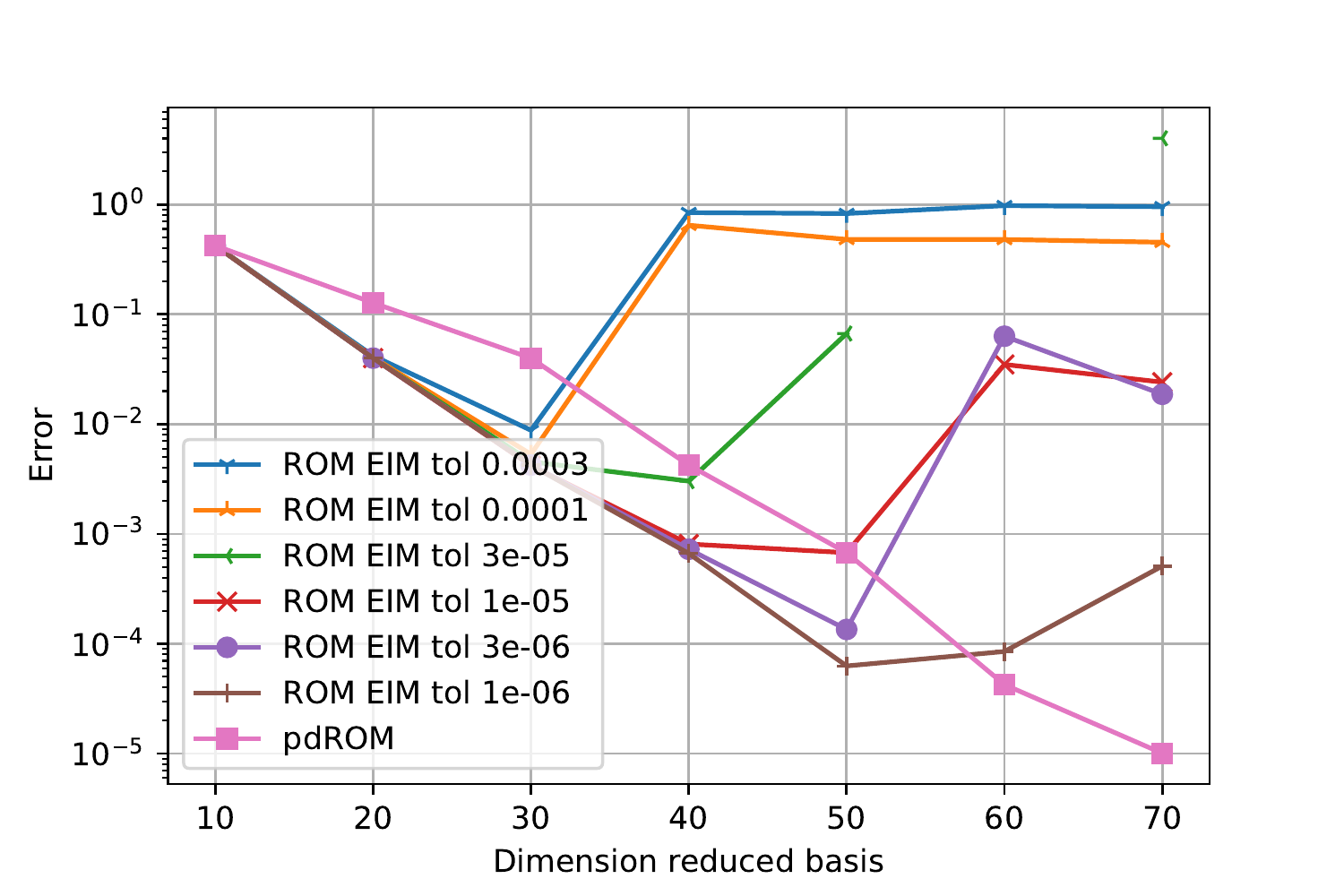}}\;
		\subfigure[EIM ROM computational times time-only reduction \label{fig:ROMEIMtimeSysSolPer}]{\includegraphics[width=0.32\textwidth]{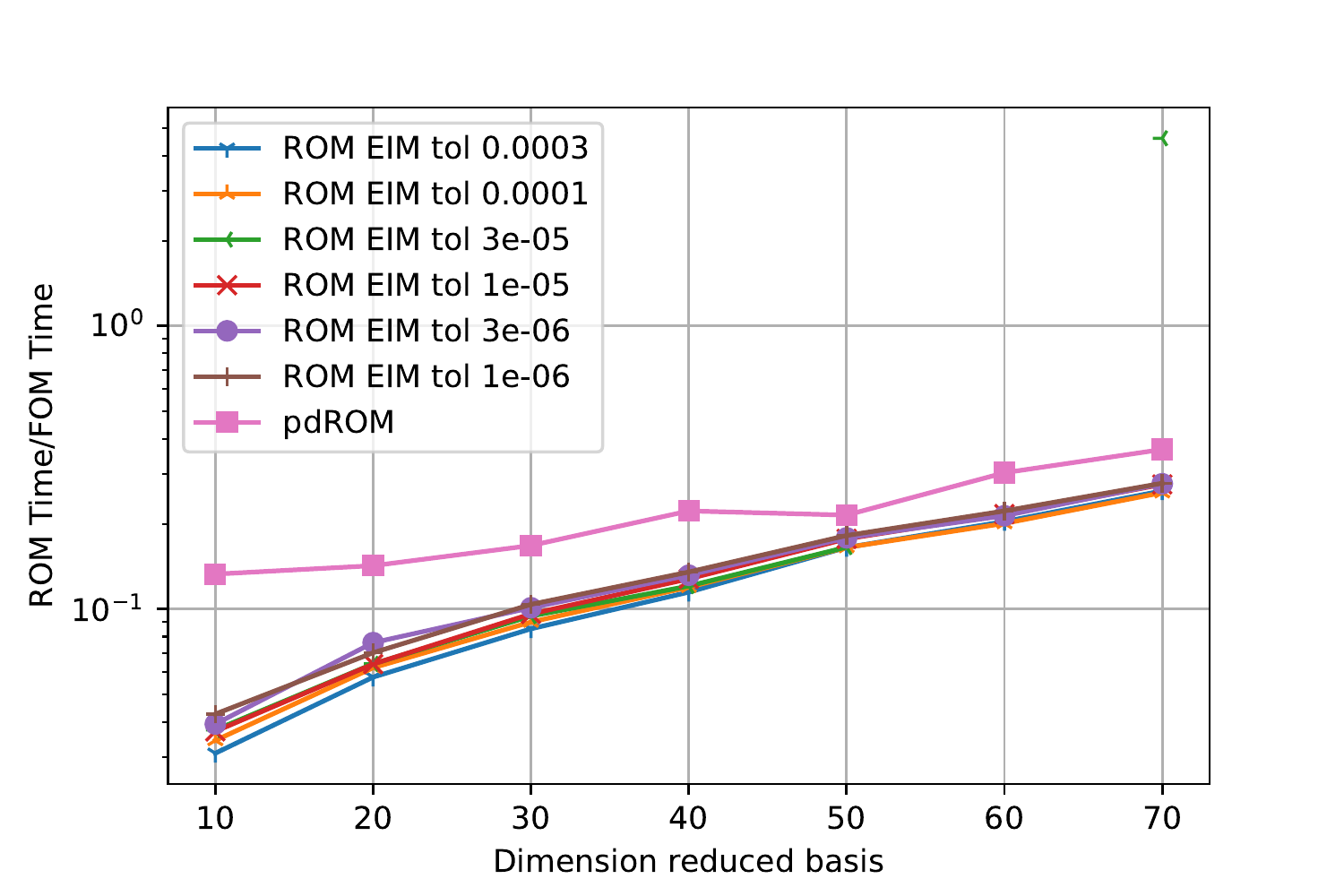}}
	\end{center}	\caption{\textbf{Solitary waves over a submerged bar.} pdROM and EIM (time-only) reduction: simulation, errors and computational time, varying POD and EIM dimensions}
\end{figure}
\begin{figure}
	\begin{center}
		\subfigure[ROM errors time-parameter reduction on parameter in the training set \label{fig:ROMEIMparamErrSysSolPerin}]{\includegraphics[width=0.37\textwidth,trim={20 0 360 35},clip]{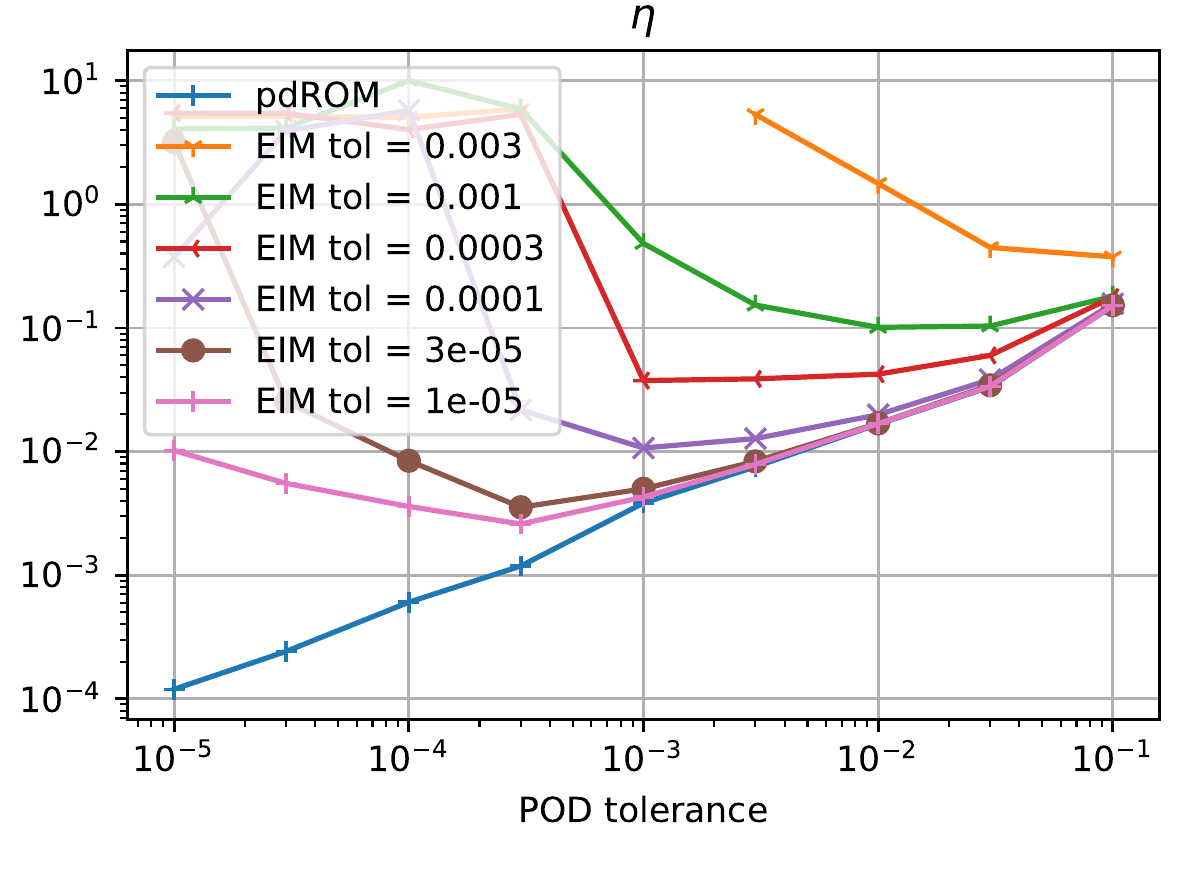}}\;
		\subfigure[ROM computational time time-parameter reduction on parameter and $u_0$ in the training set \label{fig:ROMEIMparamTimeSysSolPerin}]{\includegraphics[width=0.34\textwidth]{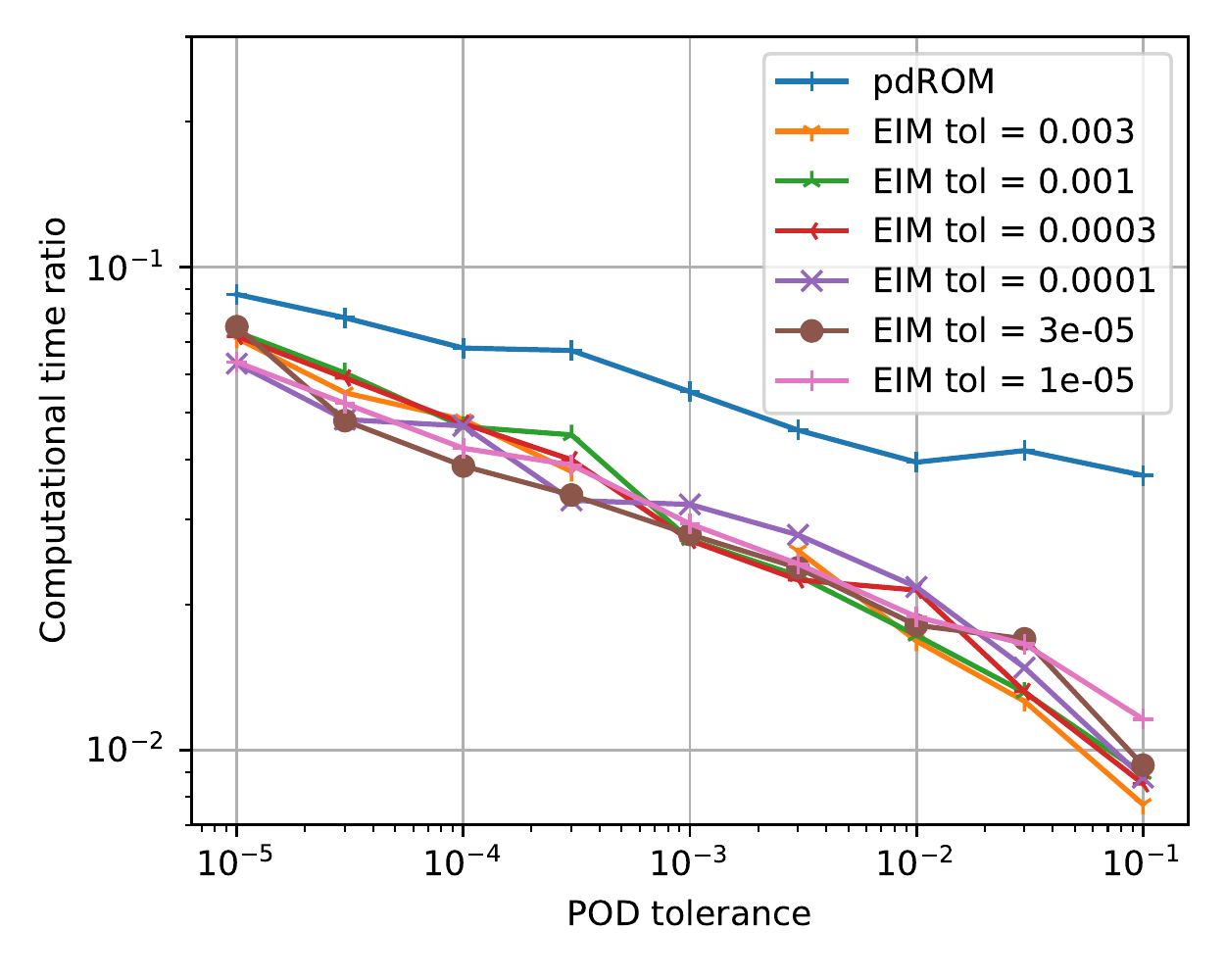}}\;
		\subfigure[\label{fig:ROMEIMparamSimSysSolPerpin}ROM solutions time-parameter reduction in the training set]
		{\includegraphics[width=0.26\textwidth,trim={20 0 285 20},clip]{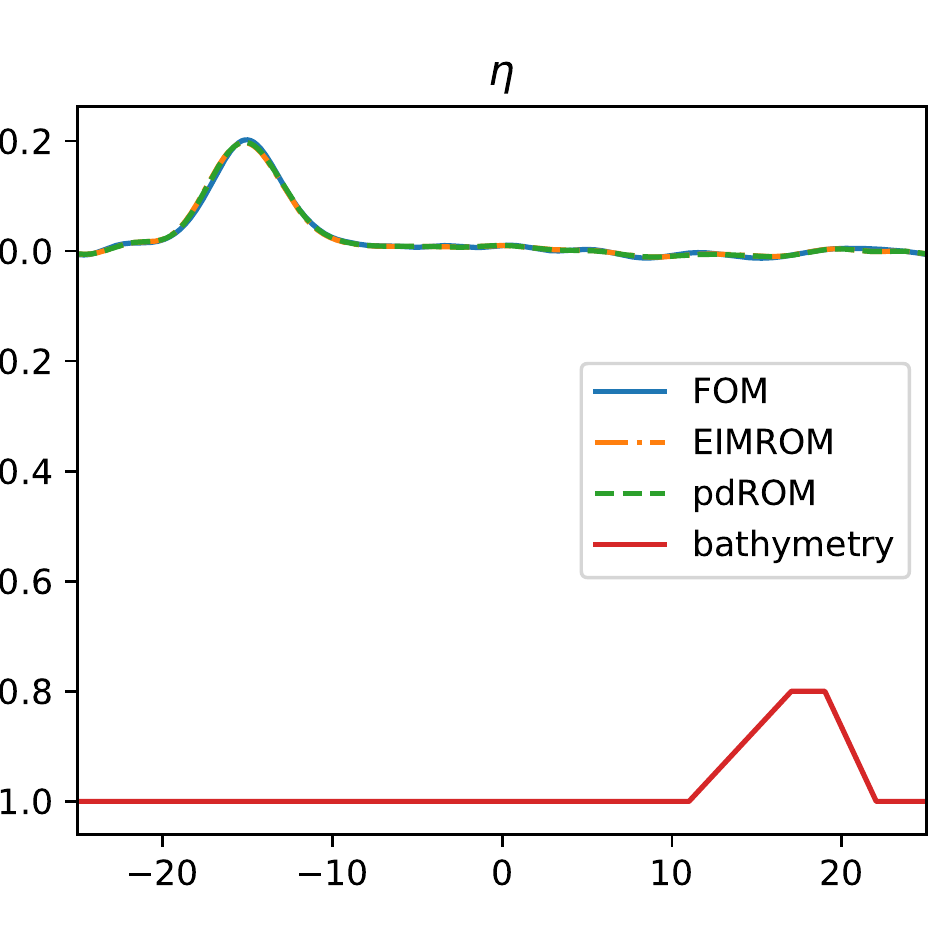}}
	\end{center}	\caption{\textbf{Solitary waves over a submerged bar.}  pdROM and EIM (time-parameter) reduction: simulation, errors and computational time, varying POD and EIM dimensions. Parameter inside the training domain. }
\end{figure}
\begin{figure}
	\begin{center}
		\subfigure[ROM errors time-parameter reduction on parameter and $u_0$ outside the training set \label{fig:ROMEIMparamErrSysSolPerout}]{\includegraphics[width=0.37\textwidth,trim={20 0 360 35},clip]{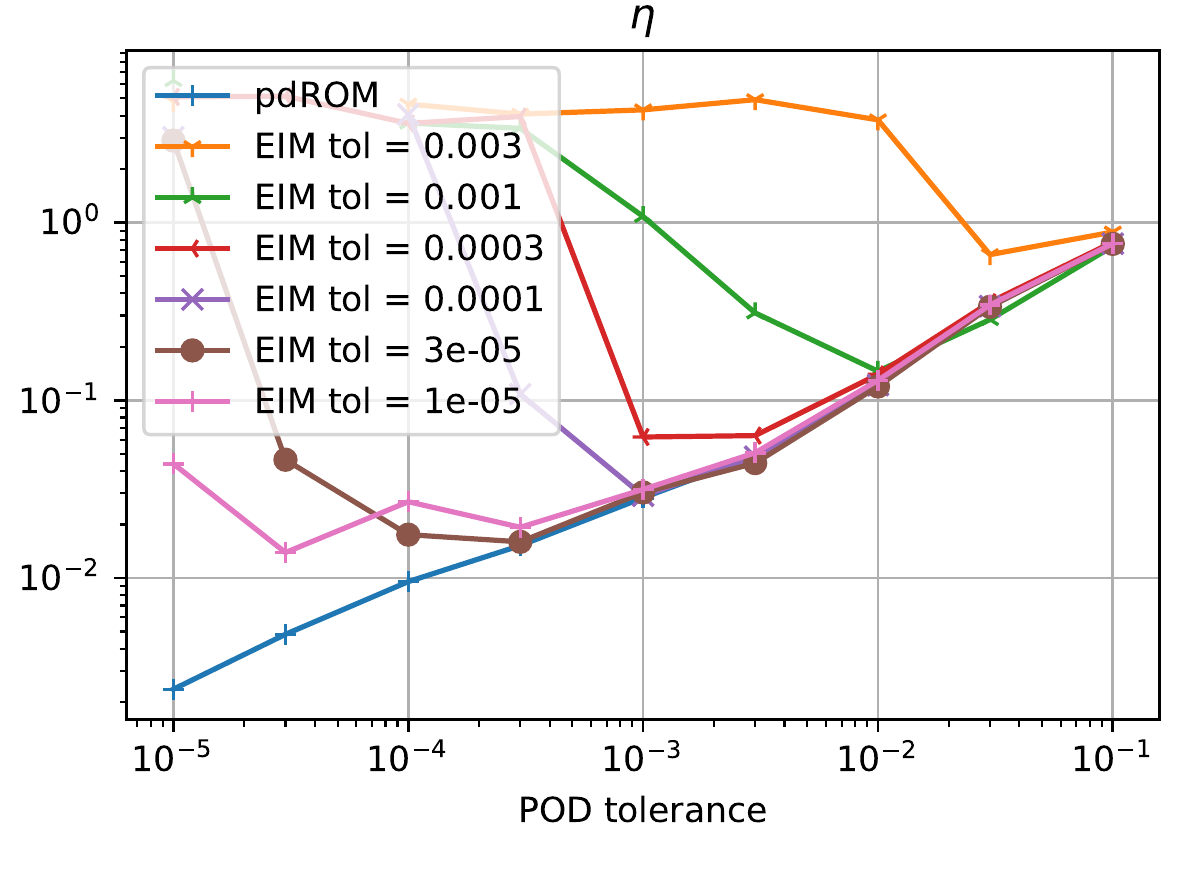}}\;
		\subfigure[ROM computational time time-parameter reduction on parameter and $u_0$ outside the training set \label{fig:ROMEIMparamTimeSysSolPerout}]{\includegraphics[width=0.34\textwidth]{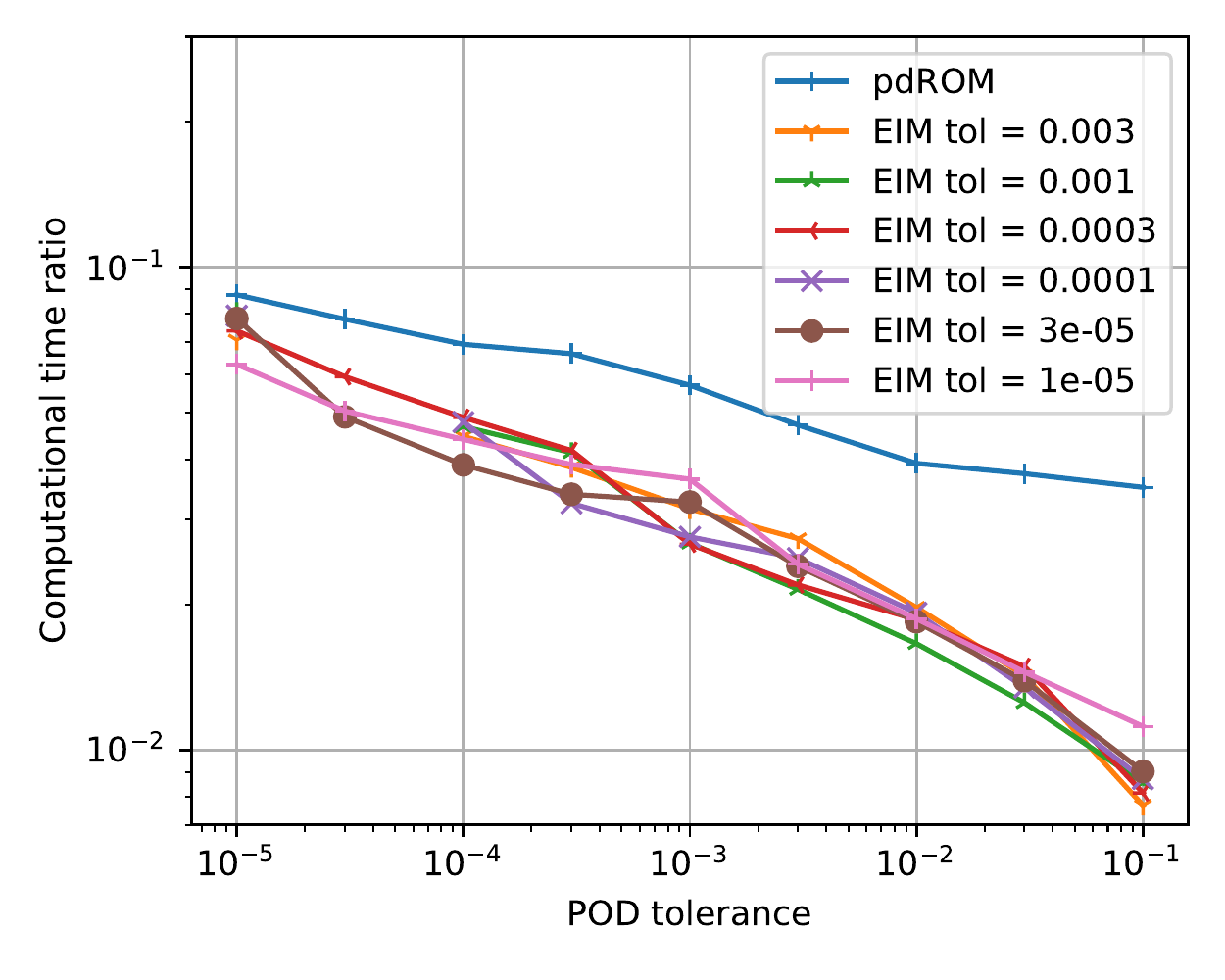}}\;
		\subfigure[\label{fig:ROMEIMparamSimSysSolPerout}pdROM and EIMROM solutions time-parameter reduction outside the training set]
		{\includegraphics[width=0.26\textwidth,trim={20 0 285 20},clip]{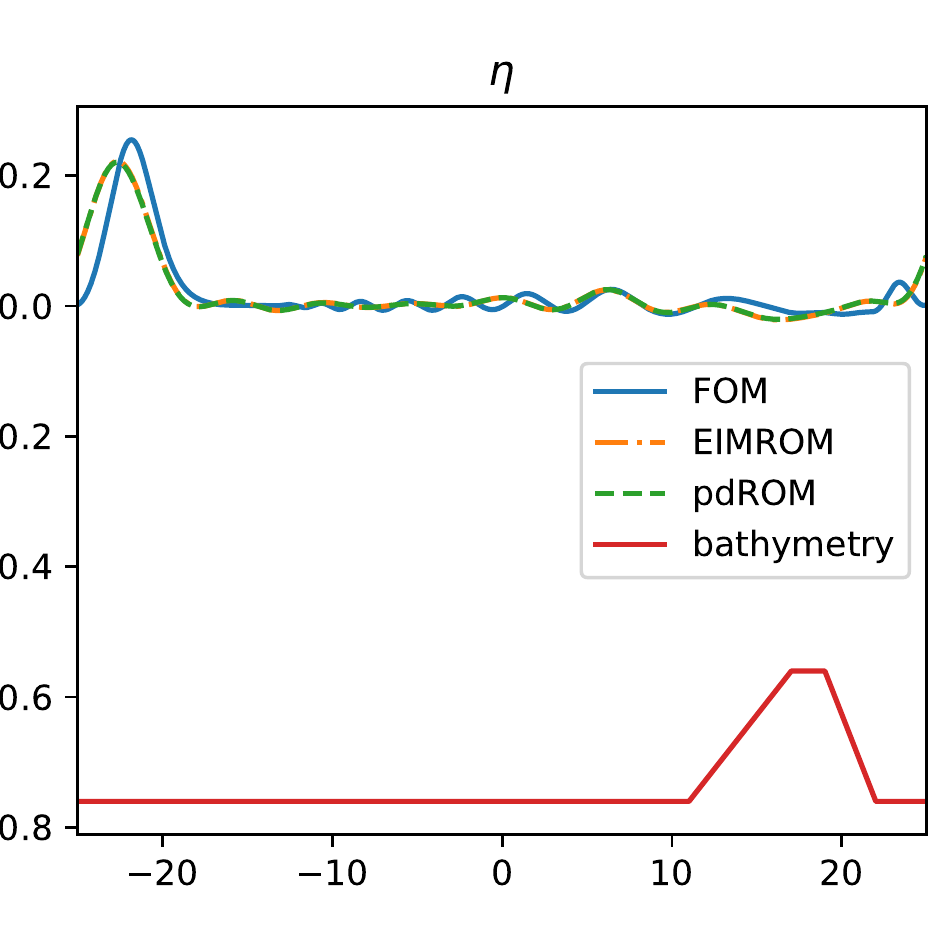}}
		\caption{ \textbf{Solitary waves over a submerged bar.}  pdROM and EIM (time-parameter) reduction: simulation, errors and computational time, varying POD and EIM dimensions. Parameter outside the training domain.}		
	\end{center}
	
\end{figure}

\begin{figure}
	\centering
	\includegraphics[width=0.9\textwidth]{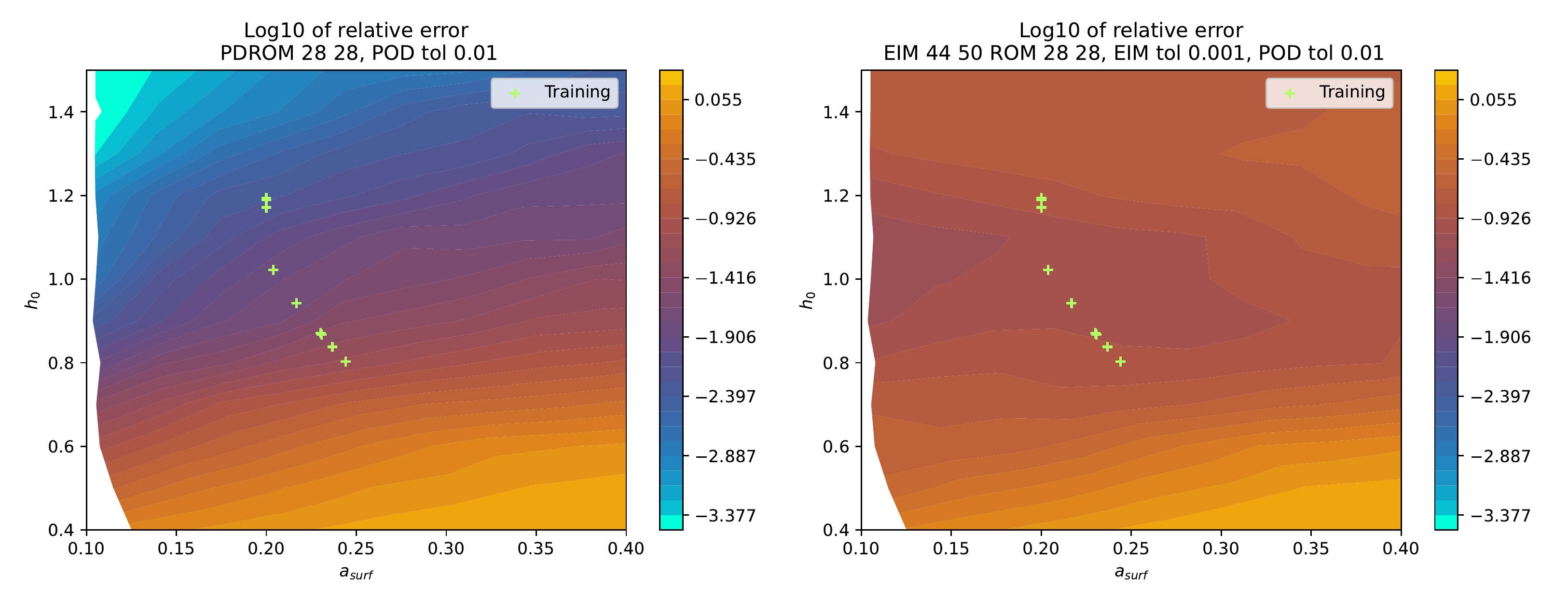}	
	\caption{\label{fig:varyingSysSolPer}\textbf{Solitary waves over a submerged bar.} Error plot varying the parameters $a_0$ and $h_0$ for pdROM $\NRB^\eta=31$ and $\NRB^q=30$ and EIMROM with $\NEIM^\eta= 40$ and $\NEIM^q=51$}
\end{figure}
For this test we set $a_0=0.2$ and $\bar{h}_0=1$. For the first phase we just compress the solutions for a simulation until $T=25$ with the POD. The problem is relatively simple and the wave does not show much dispersion. We can see with different choices of $\NRB^\eta=\NRB^q$ the behavior of the pdROM approximated solutions in \cref{fig:ROMsolSysSolPer}. Adding also the EIM to the algorithm we obtain expected results as we can see both in the error behavior in \cref{fig:ROMEIMerrSysSolPer}, and in computational times of \cref{fig:ROMEIMtimeSysSolPer}, where we can achieve around 4\% to 10\% of FOM computational time for pdROM and between 1\% and 8\% for EIMROM. In this case, the system of the FOM requires more computational time than the scalar case one, as it is composed of more terms.

Now, we consider 10 snapshots with randomly chosen parameters $h_0\in [ 0.8,1.2]$ and $a_0 \in [0.16,0.24]$.   This training set is used to compute the POD and EIM basis functions. Simulating pdROM and EIMROM for $a_0=0.2$ and $h_0=1$, with relatively few basis functions the error is already low, see \cref{fig:ROMEIMparamErrSysSolPerin} and the computational costs stay in the same range between 1\% and 10\% of FOM computational times, see \cref{fig:ROMEIMparamTimeSysSolPerin}. In \cref{fig:ROMEIMparamSimSysSolPerpin} we see that  with only $\NRB^\eta=20$, $\NRB^q=20$ for 0.03 POD tolerance the pdROM already approximates very well the solution and adding the EIM with $\NEIM^\eta=47$ and $\NEIM^q=60$ for 0.0003 EIM tolerance results in an accurate approximation of the FOM in just 7\% of the FOM computational time.

Considering $h_0=0.76$ and $a_0=0.252$ slightly outside the training set we can already see larger errors (one order of magnitude) in \cref{fig:ROMEIMparamErrSysSolPerout} and comparable computational times in \cref{fig:ROMEIMparamTimeSysSolPerout}. In \cref{fig:ROMEIMparamSimSysSolPerout} we see that the FOM is a bit more oscillating than in the previous case and that the pdROM and EIMROM, for the same parameters used above, struggle more with approximating accurately the solution, still being not too far from the FOM solution.

For $\NRB^\eta=28,\,\NRB^q=28$ relative to a POD tolerance of 0.01 and $\NEIM^\eta = 44,\,\NEIM^q=50$ for EIM tolerance of 0.001, we see in \cref{fig:varyingSysSolPer} how the parameters $a_0$ and $h_0$ influence the relative error of the pdROM and EIMROM approximations. The structure is similar to the scalar case, where increasing the nonlinearity of the problem the oscillations rise and the error follows. Again, we see that the EIMROM do not obtain such an accurate solution as pdROM for problems close to the training set, though using much less computational time (around 3 times faster than pdROM).

\subsection{Monochromatic waves on a submerged bar}  
\renewcommand{\folderTest}{Test_systemSourceTrapezoid}
\begin{figure}
	\begin{center}
		\subfigure[ROM solutions time-only reduction \label{fig:ROMsolSysBump}]{\includegraphics[width=0.32\textwidth,trim={20 0 460 20},clip]{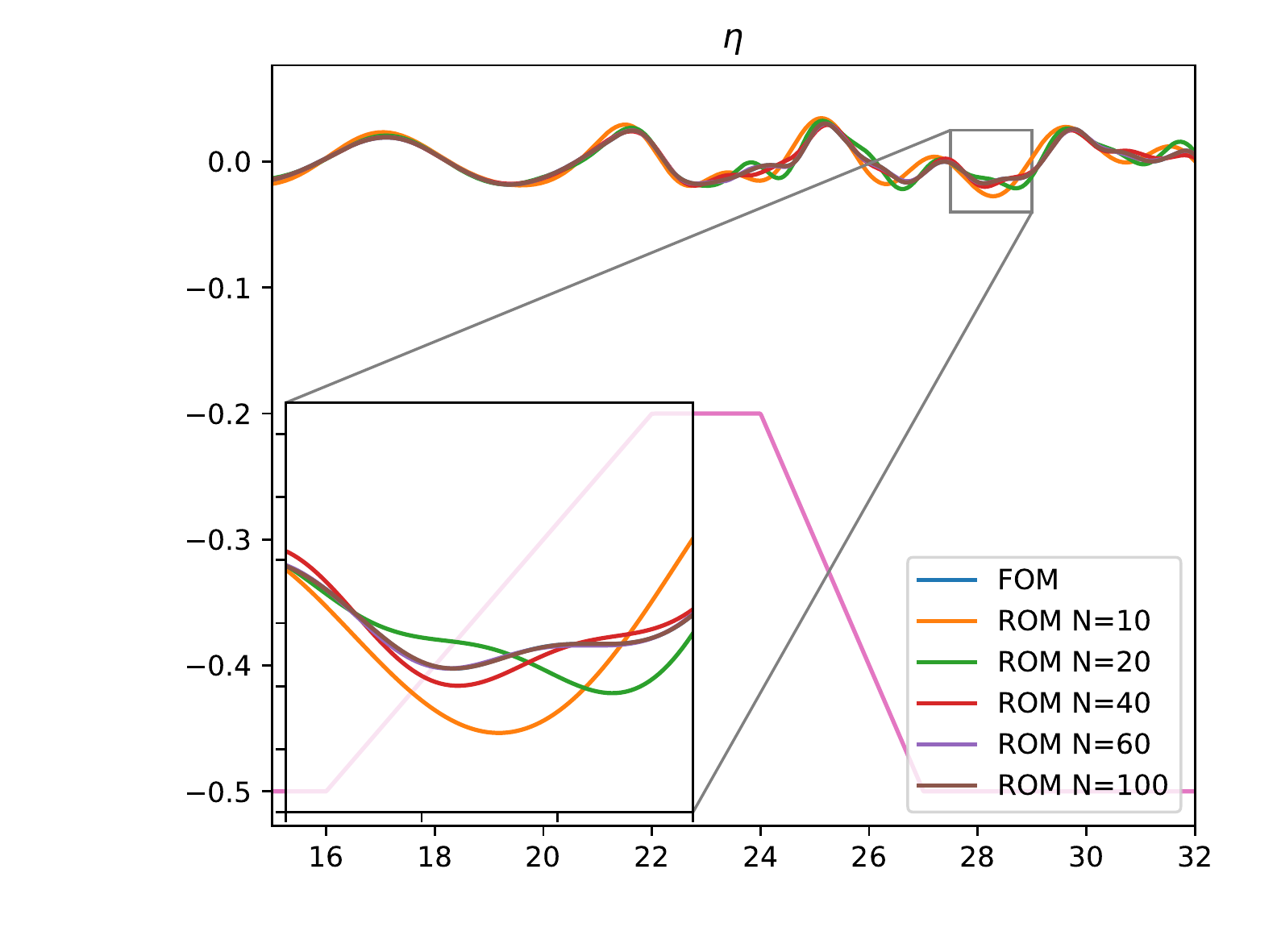}}\;
		\subfigure[EIM ROM errors time-only reduction \label{fig:ROMEIMerrSysBump}]{\includegraphics[width=0.32\textwidth]{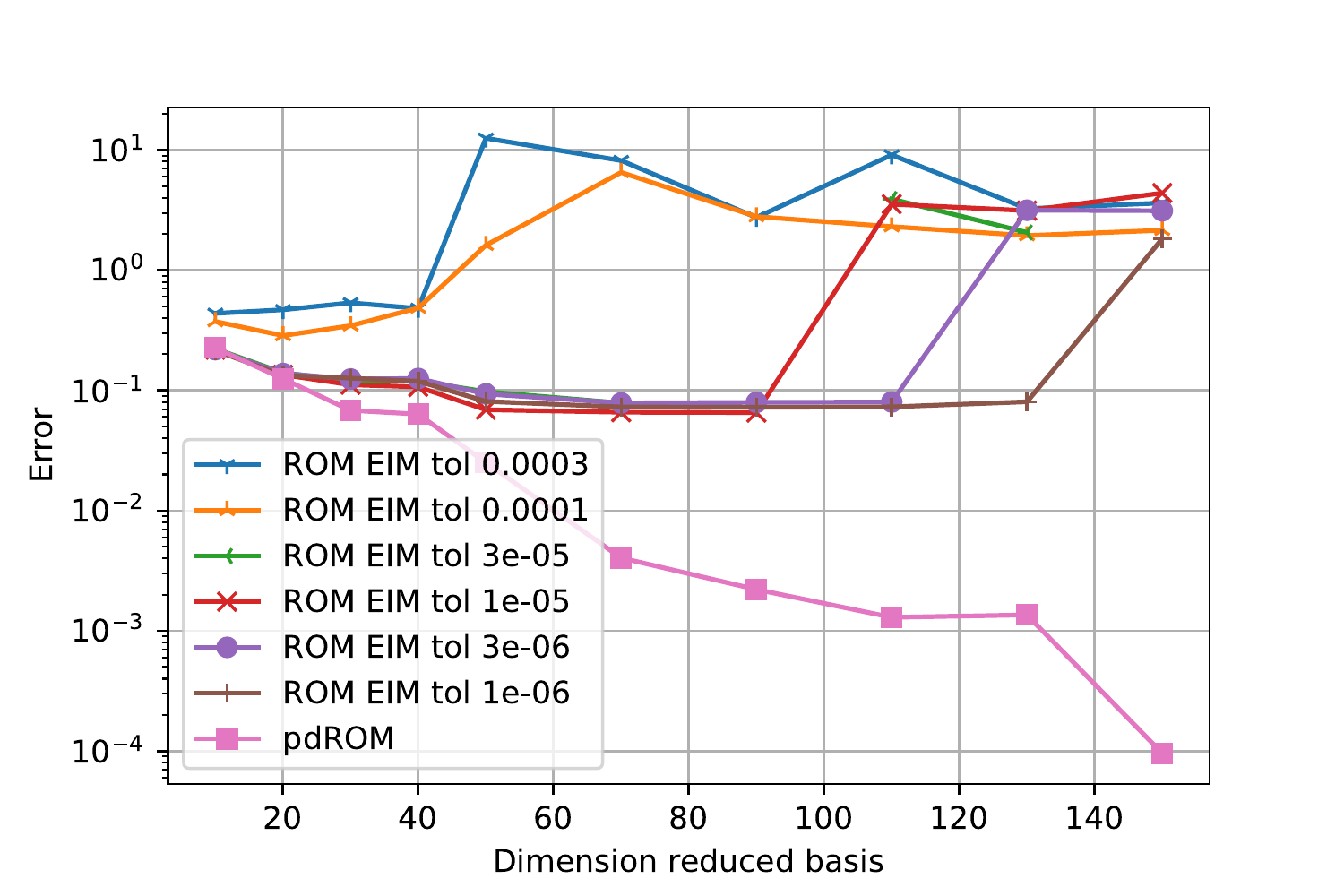}}\;
		\subfigure[EIM ROM computational times time-only reduction \label{fig:ROMEIMtimeSysBump}]{\includegraphics[width=0.32\textwidth]{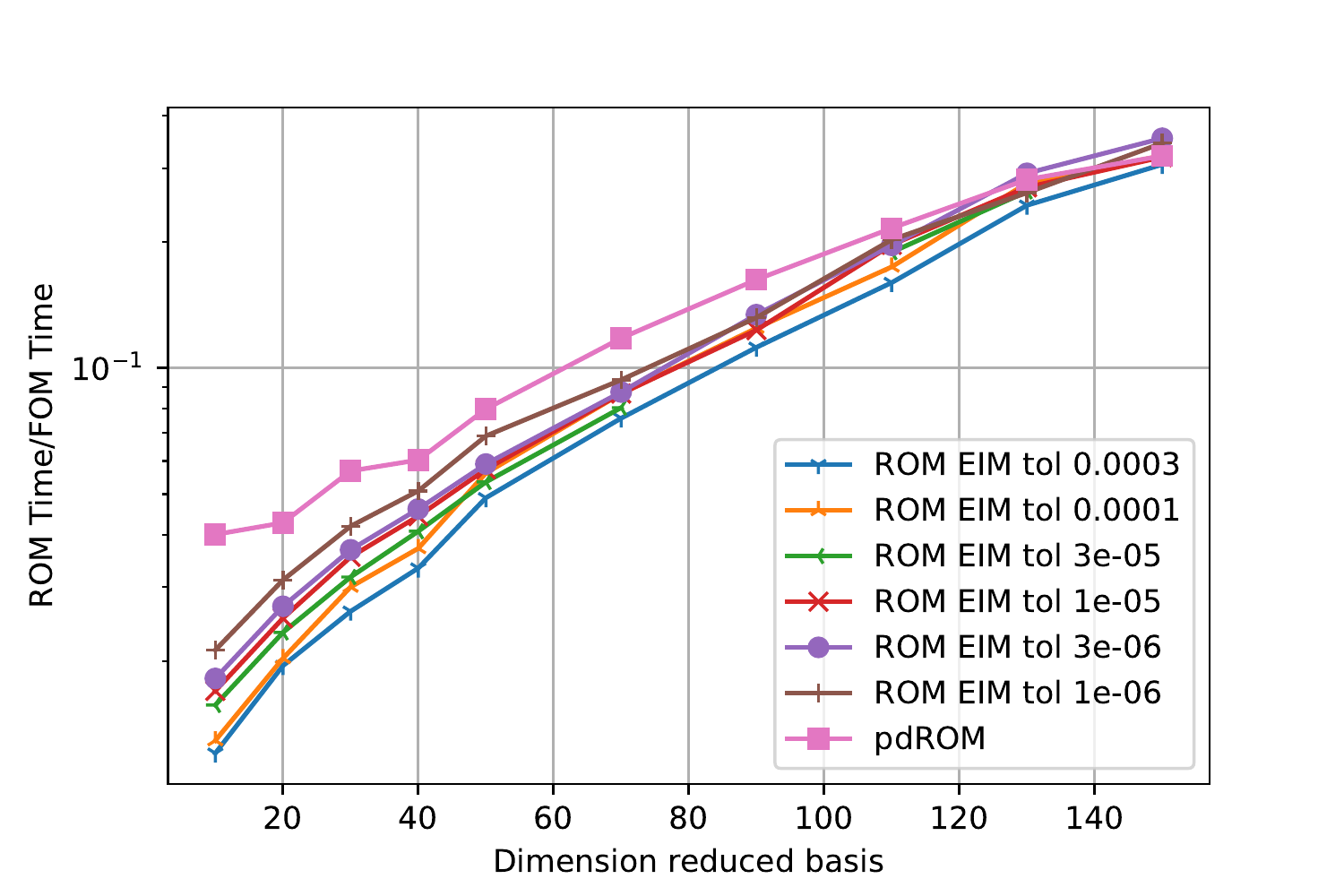}}\\
		\subfigure[EIM ROM simulation time-only reduction with $\NRB^\eta=\NRB^q = 50$ and $tol_{EIM} = 10^{-4},\, 3 \cdot 10^{-5},\,10^{-6}$ respectively from left to right \label{fig:ROMEIMsimSysBump}]{\includegraphics[width=0.32\textwidth]{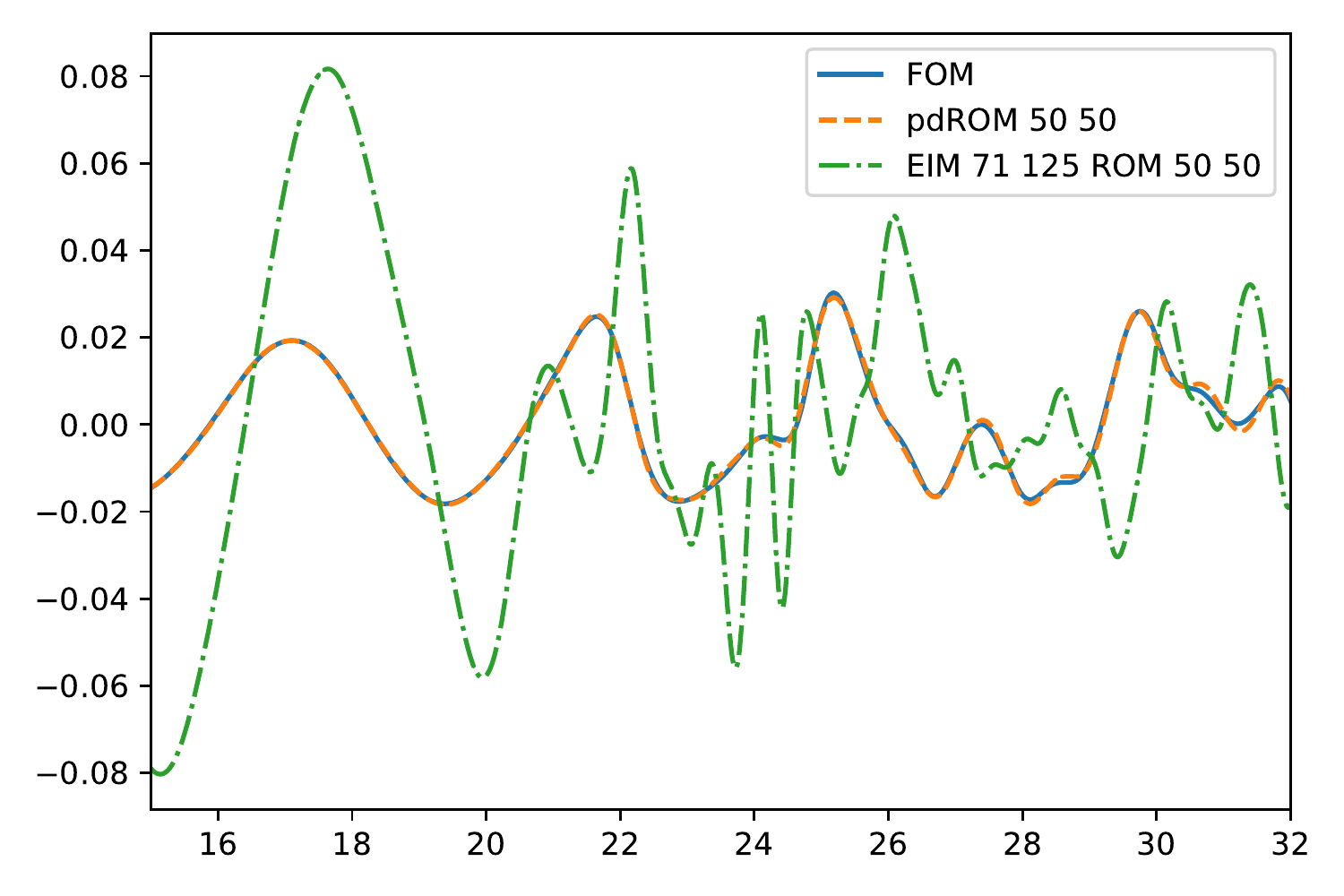}
			\includegraphics[width=0.32\textwidth]{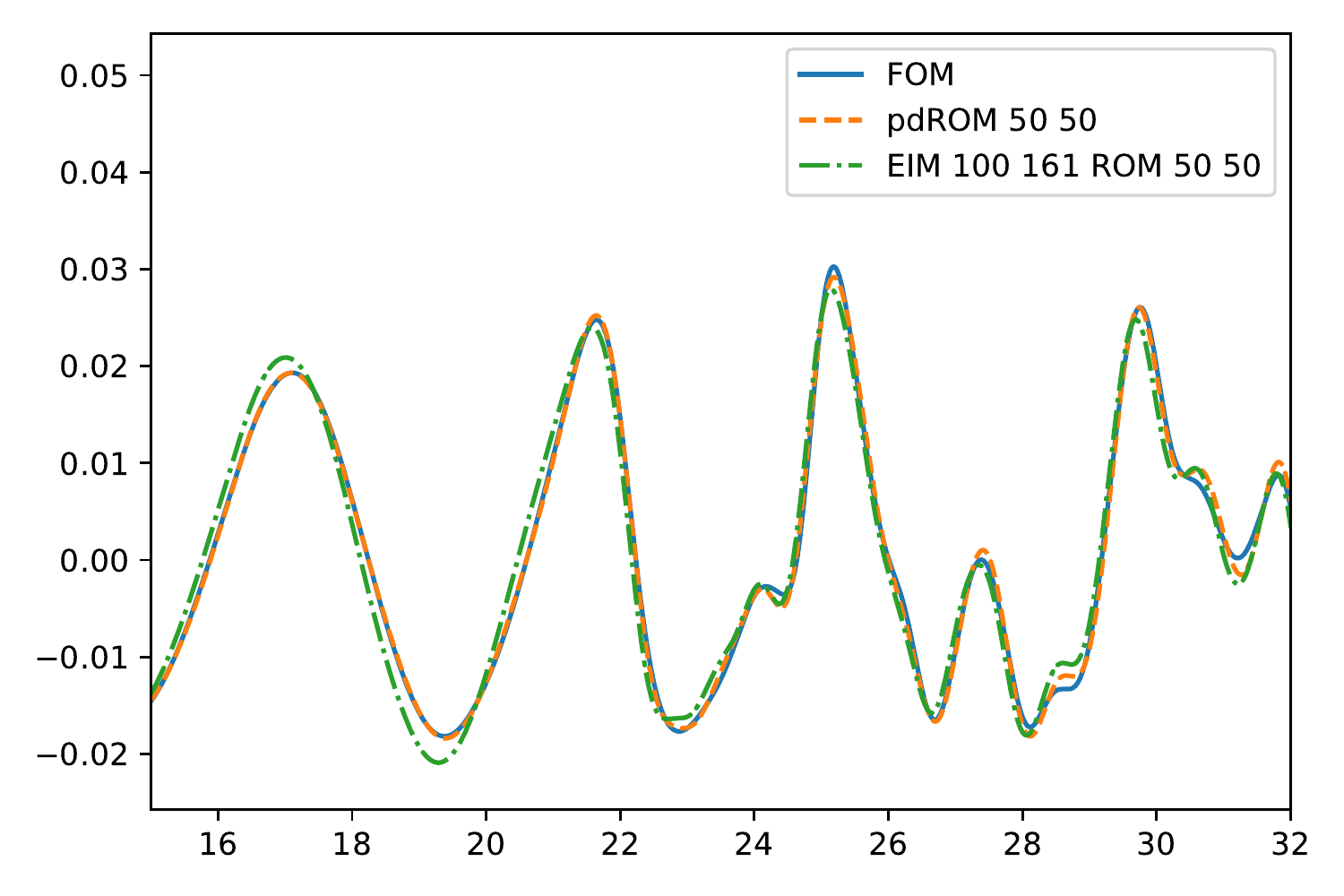}\;
			\includegraphics[width=0.32\textwidth]{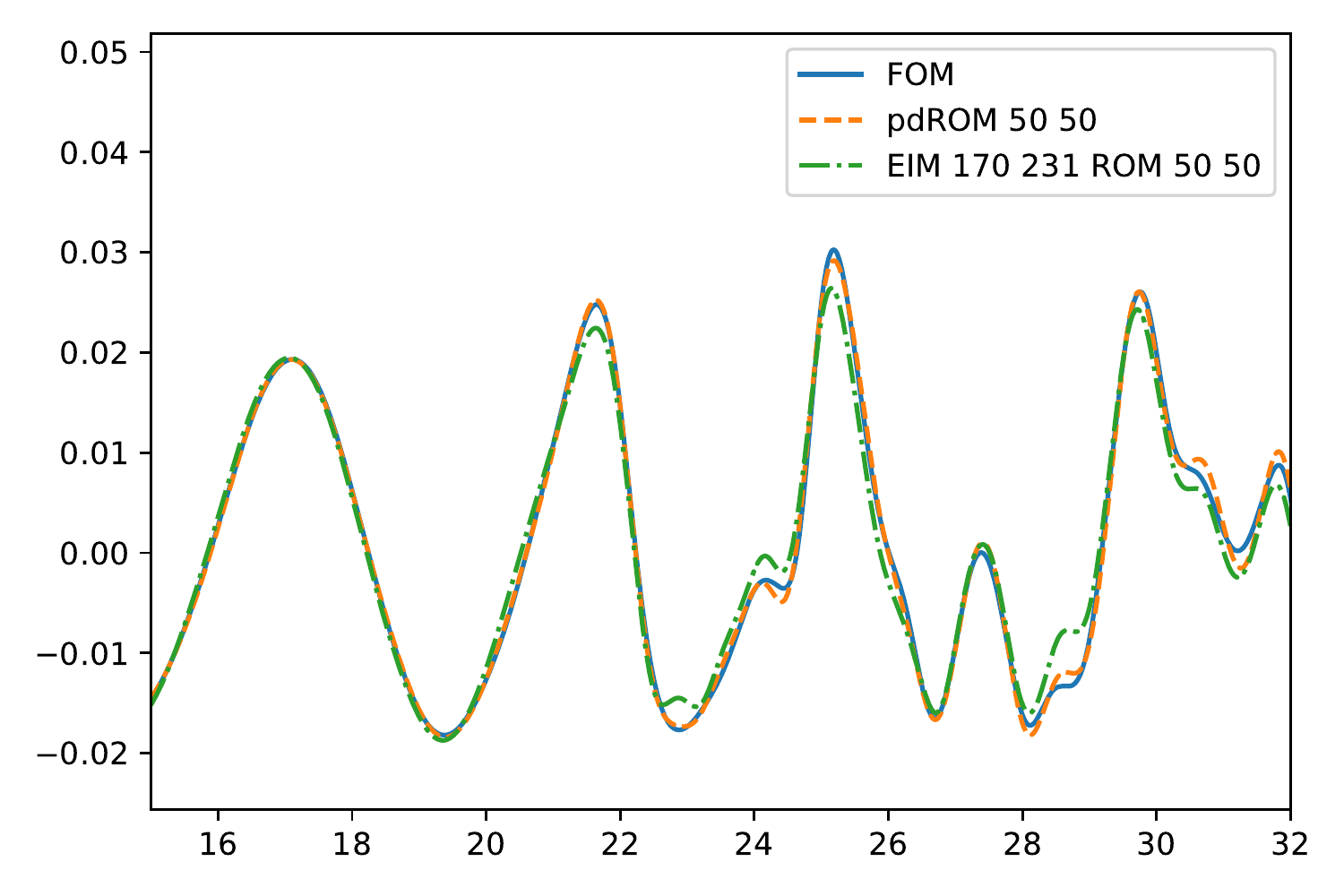}}
	\end{center}	\caption{ \textbf{Monochromatic waves on a submerged bar.}   pdROM and EIM (time-only) reduction: simulation, errors and computational time, varying POD and EIM dimensions}
\end{figure}

\begin{figure}
	\begin{center}
		\subfigure[ROM errors time-parameter reduction on parameter in the training set \label{fig:ROMEIMparamErrSysBumpin}]{\includegraphics[width=0.36\textwidth,trim={20 0 360 35},clip]{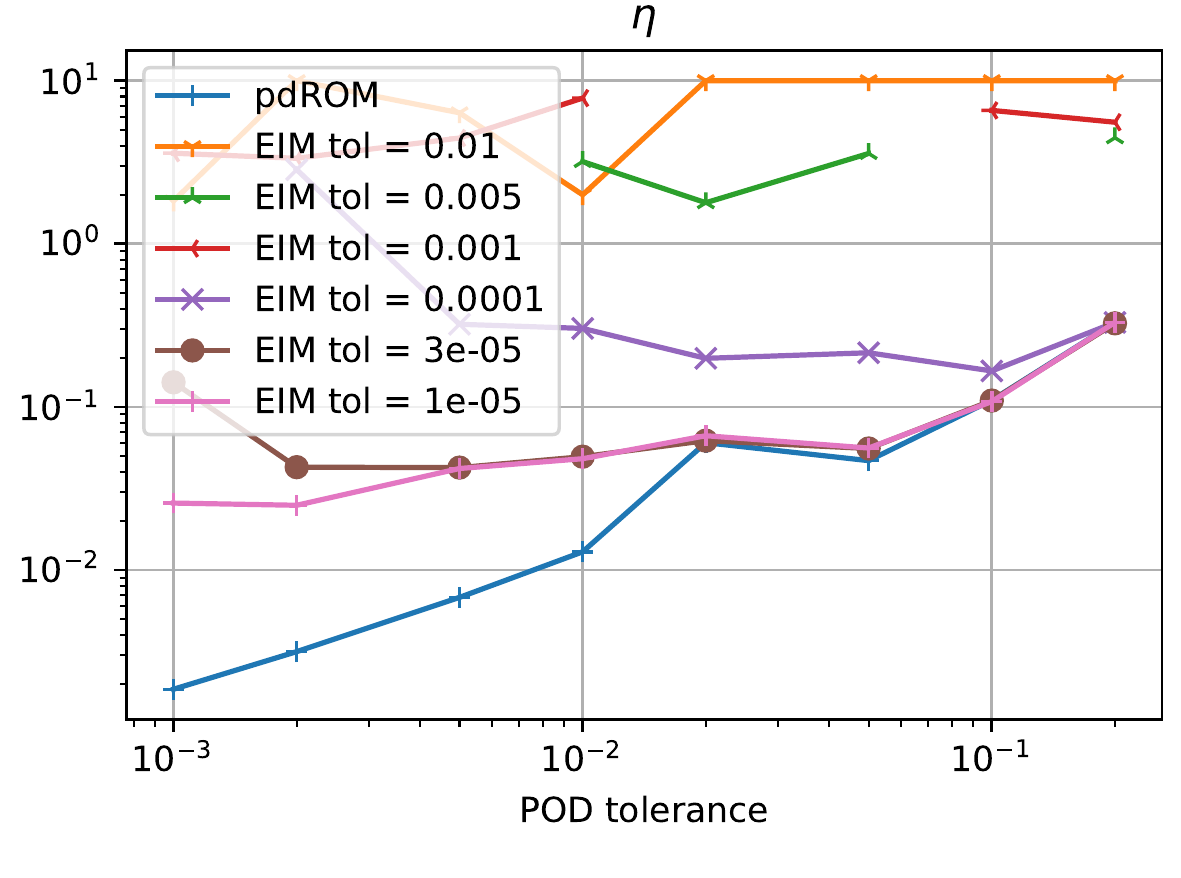}}\;
		\subfigure[ROM computational time time-parameter reduction on parameter in the training set \label{fig:ROMEIMparamTimeSysBumpin}]{\includegraphics[width=0.33\textwidth]{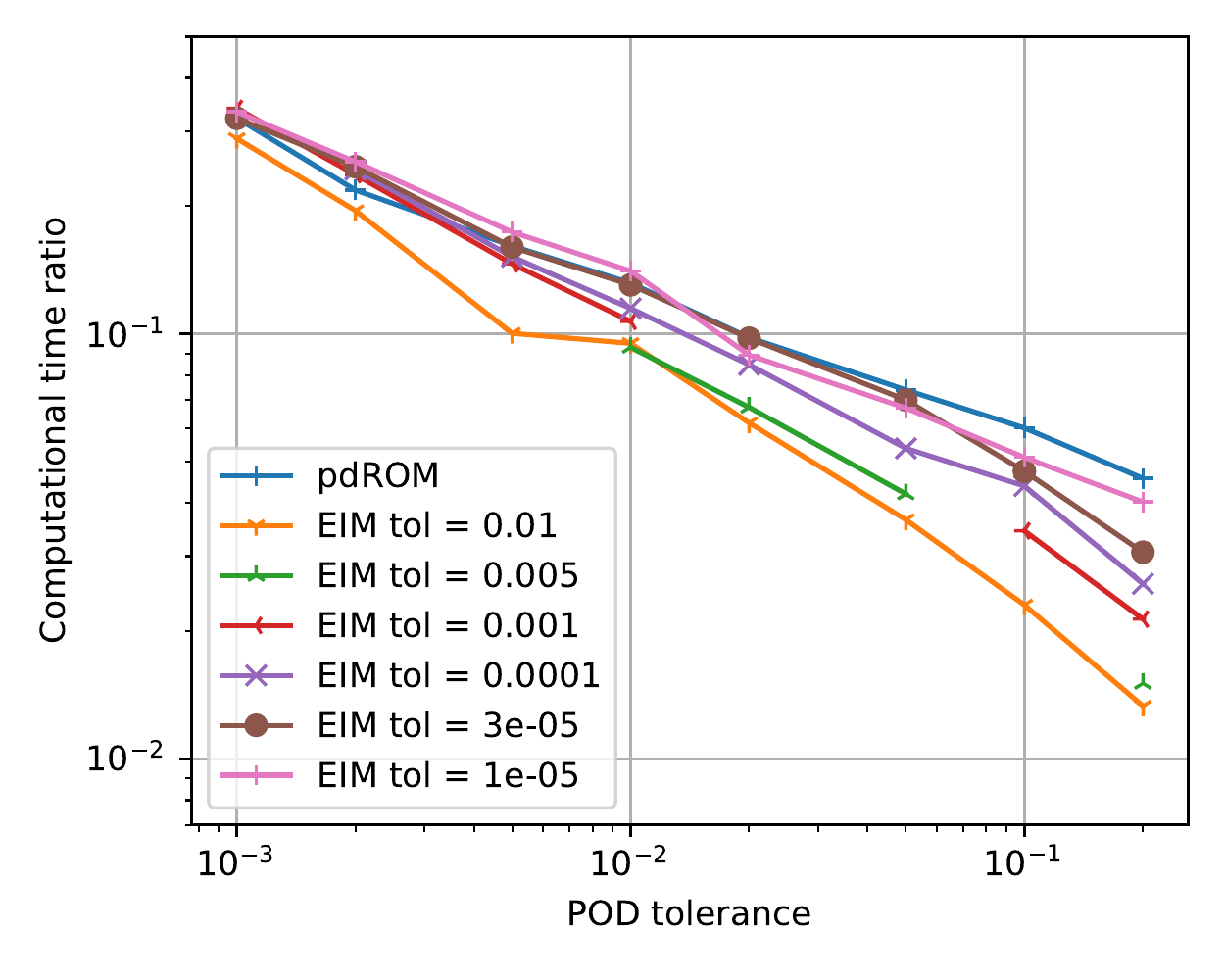}}\;
		\subfigure[\label{fig:ROMEIMparamSimSysBumpin}pdROM and EIMROM solutions time-parameter reduction]
		{\includegraphics[width=0.28\textwidth,trim={10 0 290 20},clip]{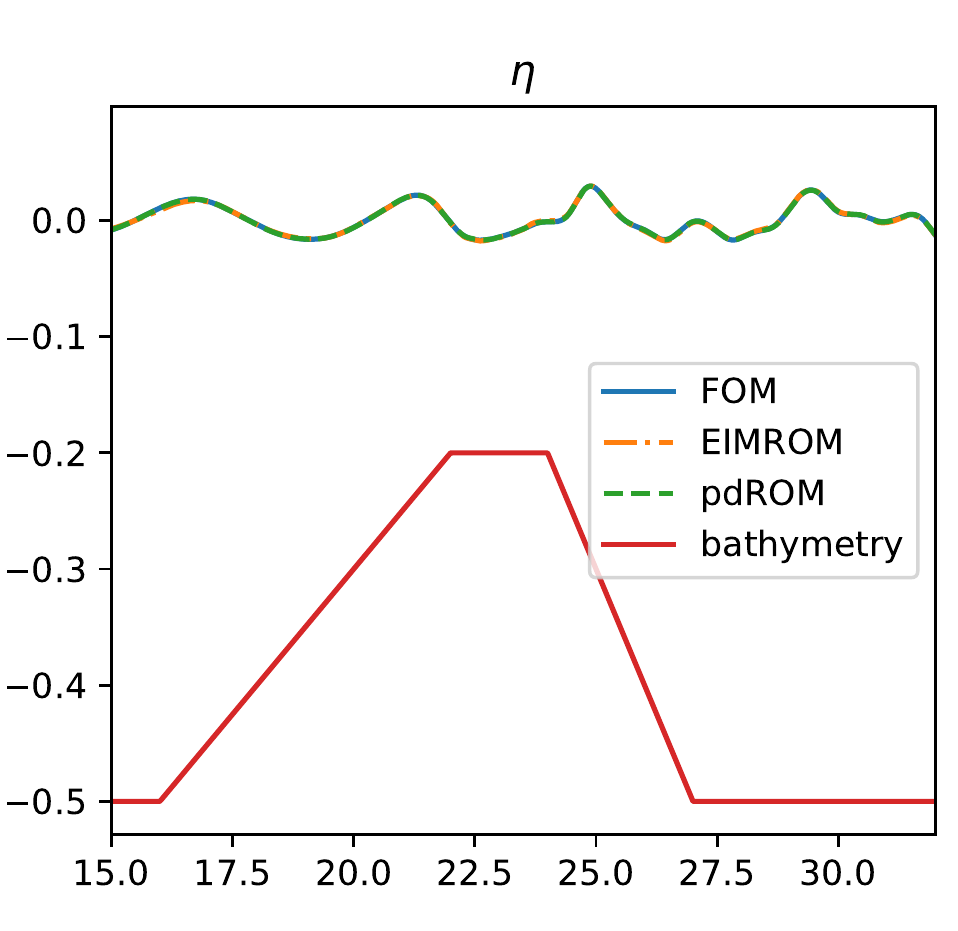}}
	\end{center}	\caption{\textbf{Monochromatic waves on a submerged bar.} pdROM and EIM (time-parameter) reduction: simulation, errors and computational time, varying POD and EIM dimensions. Parameter inside the training domain. }
\end{figure}
\begin{figure}
	\begin{center}
		\subfigure[ROM errors time-parameter reduction on parameter outside the training set \label{fig:ROMEIMparamErrSysBumpout}]{\includegraphics[width=0.36\textwidth,trim={20 0 360 35},clip]{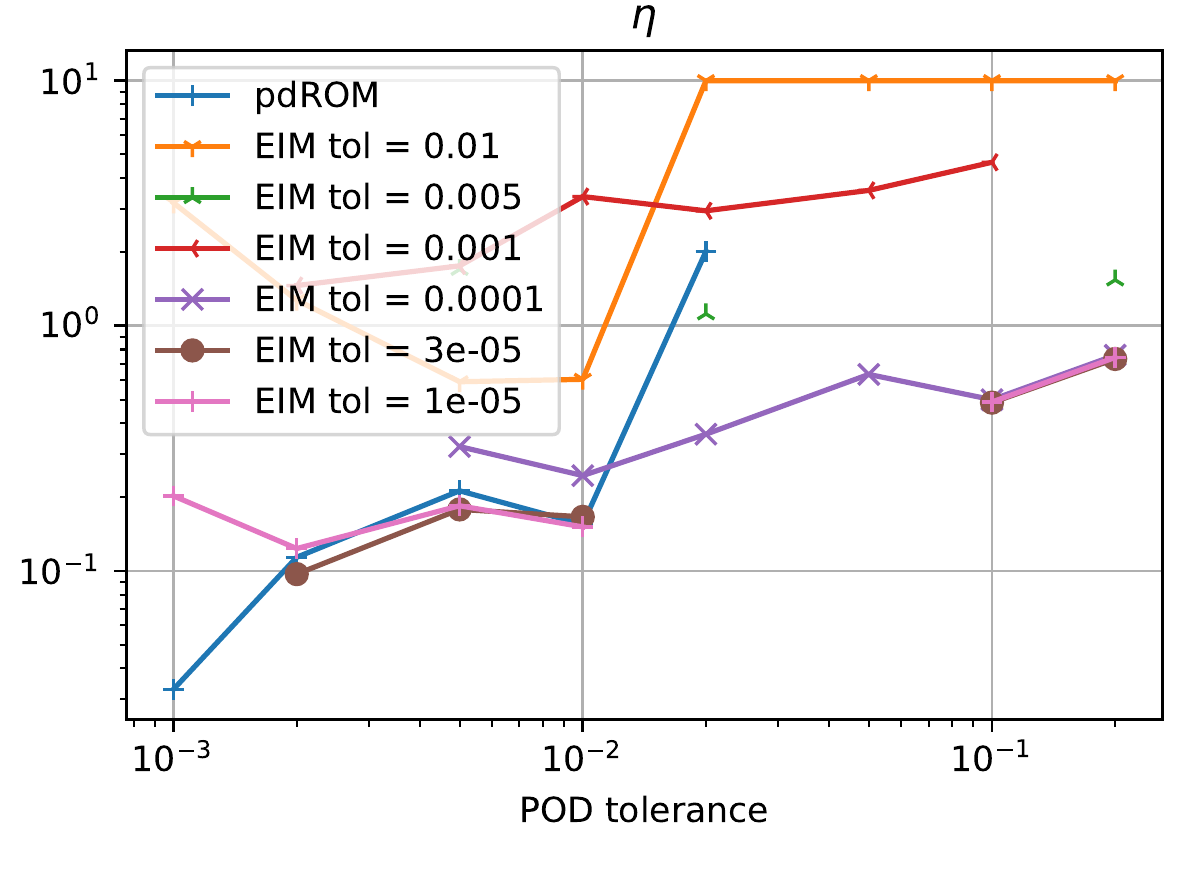}}\;
		\subfigure[ROM computational time time-parameter reduction on parameter and $u_0$ outside the training set \label{fig:ROMEIMparamTimeSysBumpout}]{\includegraphics[width=0.33\textwidth]{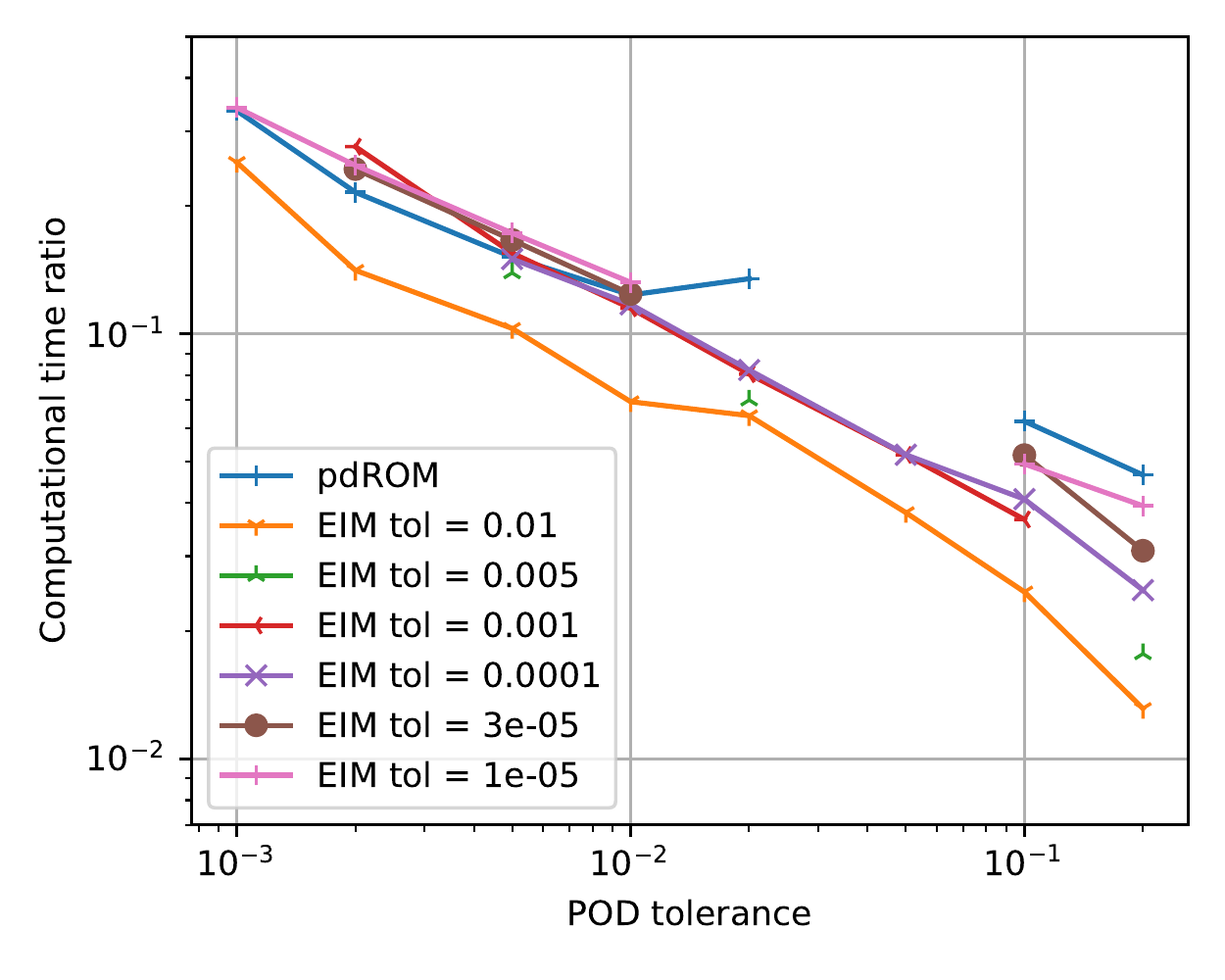}}\;
		\subfigure[\label{fig:ROMEIMparamSimSysBumpout}pdROM and EIMROM solutions time-parameter reduction]
		{\includegraphics[width=0.28\textwidth,trim={10 0 290 20},clip]{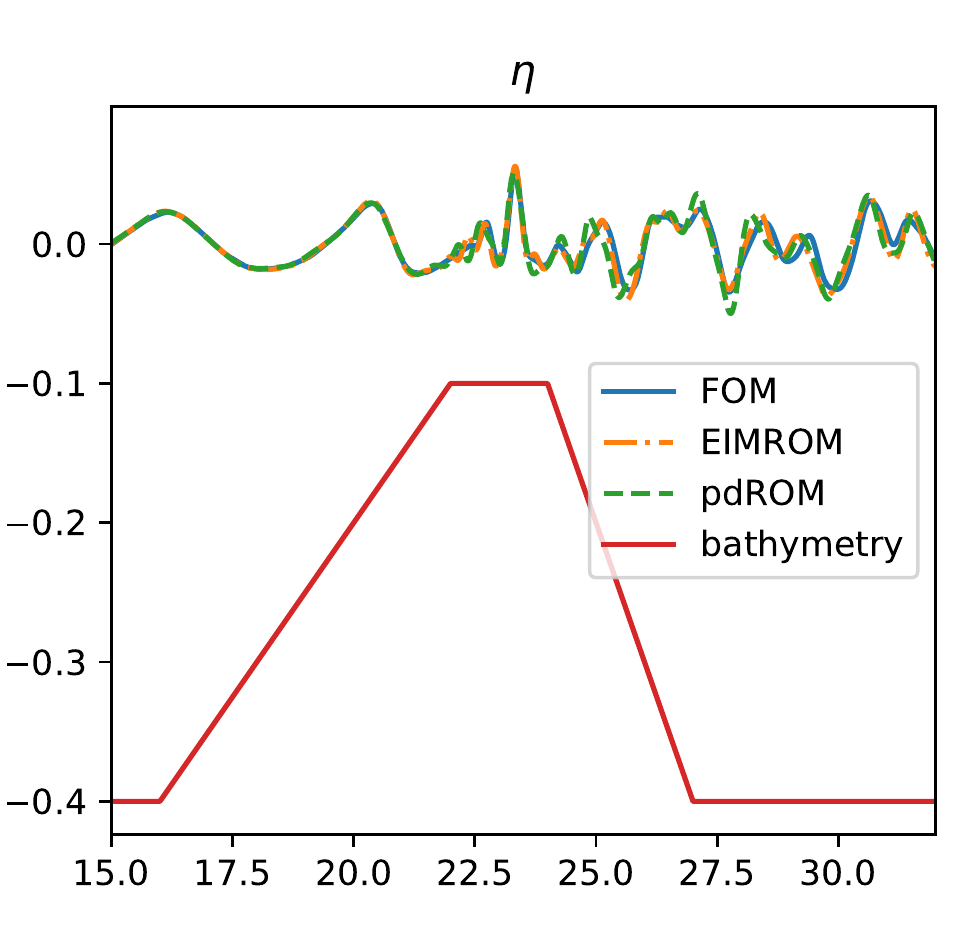}}
		\caption{\textbf{Monochromatic waves on a submerged bar.}  pdROM and EIM (time-parameter) reduction: simulation, errors and computational time, varying POD and EIM dimensions. Parameter outside the training domain.}		
	\end{center}
	
\end{figure}

\begin{figure}
	\centering
	\includegraphics[width=0.9\textwidth]{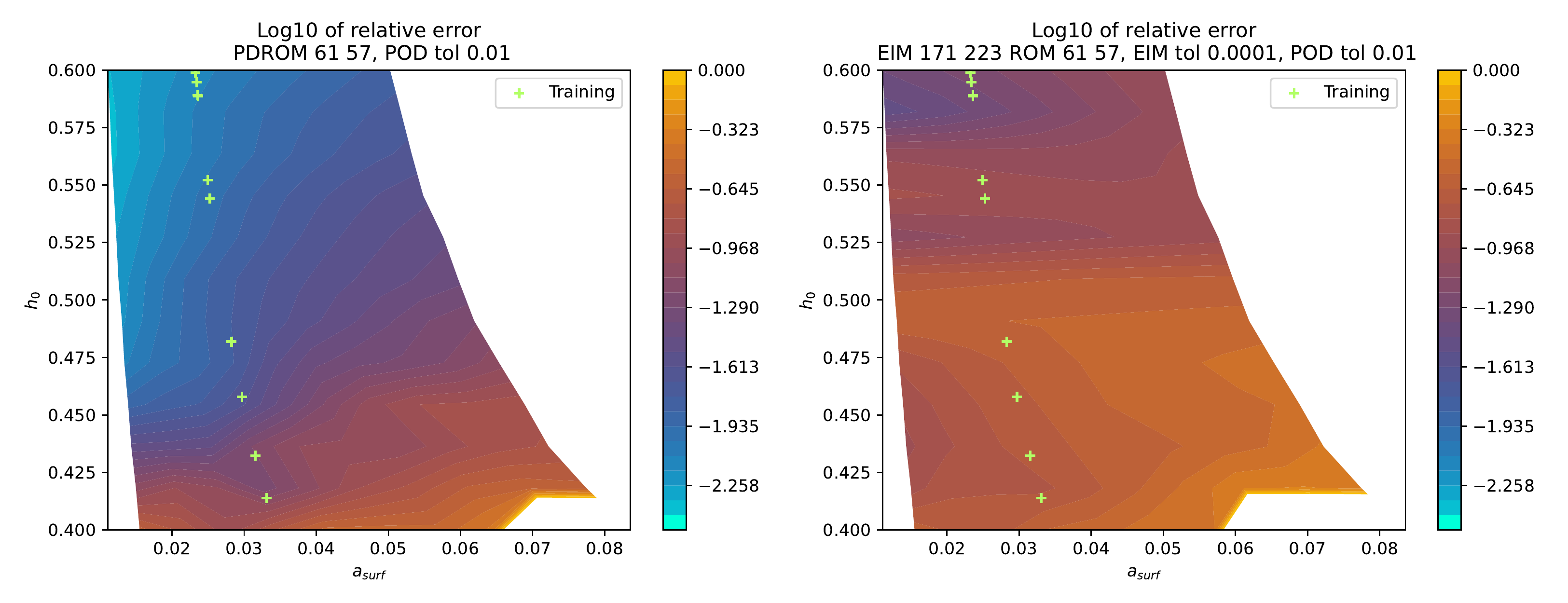}
	\includegraphics[width=0.9\textwidth]{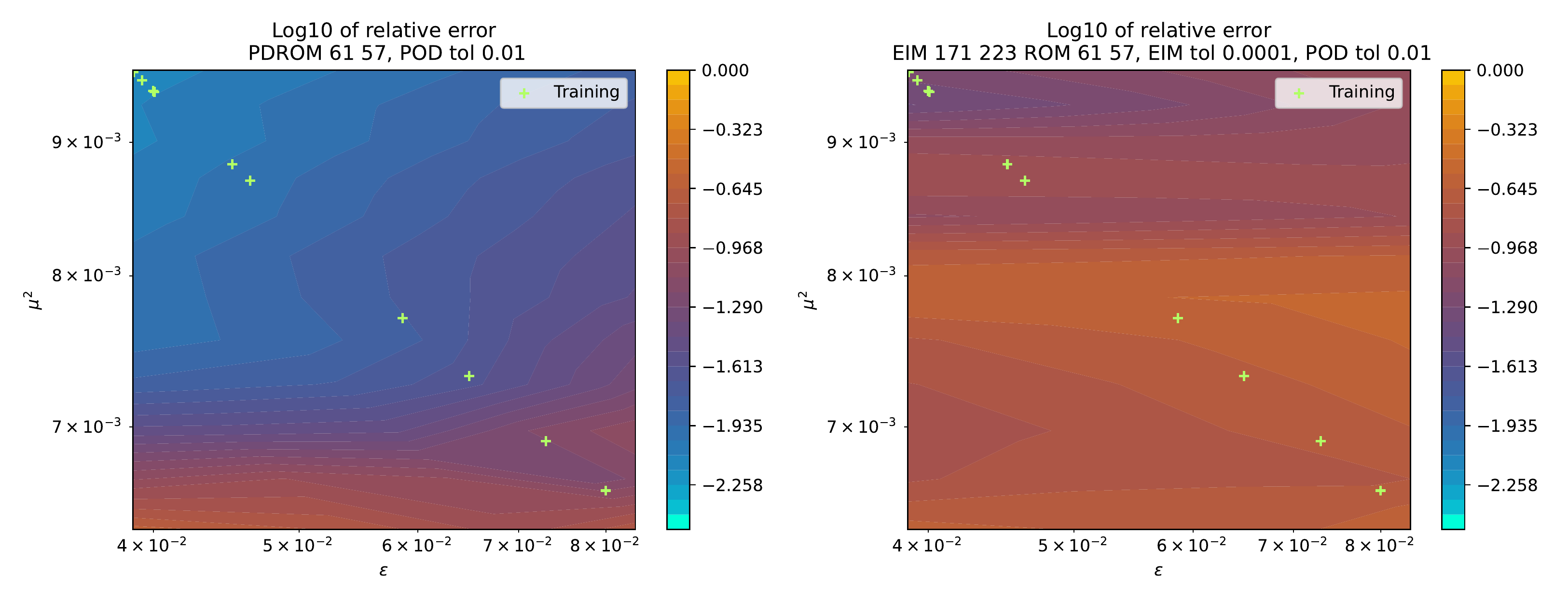}
	\caption{\label{fig:varyingSysBump}\textbf{Monochromatic waves on a submerged bar.} Error plot varying the parameters $a_0$ and $h_0$ for pdROM $tol_{POD}=10^{-2}$ $\NRB^\eta=61$ and $\NRB^q=57$ and EIMROM with $tol_{EIM}=10^{-4}$ and $\NEIM^\eta= 121$ and $\NEIM^q=235$}
\end{figure}
The last test is very challenging as its solutions are very close to show steep gradients and transport phenomena. Moreover, the water level is often close to the bathymetry and, if not well represented, the approximation can result into invalid values. Fortunately, the steep gradients are all close to the end of the trapezoid and with not so many basis functions it is possible to approximate them. We take $a_0=0.027$, $\bar{h}_0=0.5$ and final time $T=40$.
Taken the POD over the time evolution of this FOM, we approximate in \cref{fig:ROMsolSysBump} with pdROM the same evolution for different $\NRB^\eta =\NRB^q$ of the same problem. We see that with 40 basis functions the pdROM approximation is close to the FOM one, but only for 60 basis functions we do not observe any large deformation. In \cref{fig:ROMEIMerrSysBump} we see more precisely the error decay for different $\NRB$ and different EIM tolerances. 
The error is quite large, and the decay with the EIM approximation is very slow. 
It must be noticed that in many of these EIMROM simulations, the water height reaches negative values and stop at an earlier time, in those cases the error is not reported or it is larger than the scale. We see for an EIM tolerances of $10^-4$, $3\cdot10^{-5}$ and 
$10^{-6}$  in \cref{fig:ROMEIMsimSysBump} how the error for the first one is very large, but also for the other two cases the error does not decrease a lot in particular in the area after the trapezoid, where the steepest part of the solutions are located.  The computational costs for this problem range between 4\% and 30\% for pdROM and between 1\% and 30\% for EIMROM as shown in \cref{fig:ROMEIMtimeSysBump}. 

Introducing more parameters in the training set, with $h_0\in[ 0.4,0.6 ]$ and $a_0\in [ 0.0216, 0.0324]$, leads to a much more rich manifold of solutions. This implies that the eigenvalue decay of the POD gets slower. For the test parameter $\bar{h}_0$ and $a_0=0.027$, we observe in \cref{fig:ROMEIMparamTimeSysBumpin} that the computational costs for pdROM and EIMROM are very similar and they are around the 5\% of the FOM computational time to obtain an error of $10^{-1}$, around 15\% of the time for an error of $10^{-2}$ and around 30\% for an error of $10^{-3}$. For very accurate results with error lower than $10^{-4}$ a computational time of more of the 30\% of the FOM computational time may be required. The pdROM simulation in \cref{fig:ROMEIMparamErrSysBumpin} shows the accurate result for $tol_{POD}=0.001$, $\NRB^\eta=99$, $\NRB^\eta=92$ and $tol_{EIM}=0.00003$, $\NEIM^\eta=213$ $\NEIM^q=291$. The number of EIM basis functions is large in order to obtain a good approximation of the solution and almost not gaining anything in terms of computational time with respect to the pdROM. This number must be large in order to avoid the simulation to oscillate so much that it hits negative values and it crashes the simulation. Again, this behavior has been observed in many other works \cite{peherstorfer2020stability,zimmermann2018geometric,ghavamian2017pod,argaud2017stabilization,chen2021eim} and we do not aim at finding better hyper reduction algorithm in this work.

For the parameter outside the training set, we chose $a_0=0.027$ and $h_0=0.4$, which correspond to the test case in \cite{ricchiuto2014upwind} which was validated with experimental data. This simulation is quite challenging also for pdROM and when it is not enough resolved it can have oscillations in $h$ going below 0. 
This happens for example for  $tol_{POD}=5\cdot10^{-2}$, as shown from the incomplete curves in 
\cref{fig:ROMEIMparamErrSysBumpout}. Increasing the dimension of the reduced space we overcome this issue. 
Computational costs are similar to the ones discussed above.  
In \cref{fig:ROMEIMparamSimSysBumpout} we plot the simulation for $tol_{POD}=0.005$, $\NRB^\eta=99$, $\NRB^\eta=92$ and $tol_{EIM}=0.00003$, $\NEIM^\eta=213$ $\NEIM^q=291$. 
Here it is the EIMROM simulation with the largest tolerance that does not crashes along the simulation. 
The computational time reduction due to the EIM at this stage is almost negligible with respect to the pdROM and the approximation is very similar. 
Both simulations oscillate more than the original FOM after the trapezoid at this resolution, though using only 15\% of the computational time of the FOM. To reach errors of the order of $10^{-1}$ we need computational costs of around the 25\% of the FOM ones as one can observe comparing \cref{fig:ROMEIMparamErrSysBumpout} and \cref{fig:ROMEIMparamTimeSysBumpout} for $tol_{POD} = 2 \cdot 10^{-3}$.

In \cref{fig:varyingSysBump} we see the strong influence that $h_0$ and $a_0$ have on these simulations. In this case, it is crucial the level of $h_0$. When this level is too low, we fall in another regime and all the hypotheses made on the dispersive equations are not valid anymore. Moreover, the simulations risk to hit negative values and break in this regime. Above this regime, for pdROM we can simply say that as the oscillations and the nonlinearity increase, the error increases. On the other side, the dispersion term is actually helping in reducing the error in the reduced space. For EIMROM we observe more areas where the algorithm struggles with obtaining good results, also for not too small $h_0$ values.

\section{Conclusion}\label{sec:conclusion}
In this work we have proposed some strategies for 
model order reduction  for dispersive waves equations.
Applications have been shown to  a BBM-KdV type model, and to the enhanced Boussinesq equations of Madsen and S{\o}rensen. Both models contain hyperbolic and dispersive terms that can be decoupled in the numerical solver. The dispersive terms can be obtained solving an elliptic and linear problem, which is well suited to be reduced with standard projection based reduction techniques after the choice of a reduced basis space. We apply the proper orthogonal decomposition on some snapshots obtained for different times and parameters to obtain a reduced space. When the gradients of the solutions are not so steep and the nonlinear character is not the dominant one, we can obtain large reduction with the pdROM in computational times, up to 20 times less, and still obtain a reliable solution. Introducing a second level of approximation with the empirical interpolation method, this factor can increase up to 100, but stability issues may arise. Hence, it is not always obvious how to choose how many interpolation functions are needed. There are some clear bounds that cannot be overtaken when using these simple reduction techniques: the nonlinearity cannot be too pronounced and the water level must be far enough from zero.
Nevertheless, the results are encouraging in particular knowing that already for one dimensional problems we can can obtain good reduction.

This idea can be easily applied to a pre-existent hyperbolic algorithm, as a shallow water code, with a cheap additional term which takes into account of the dispersive effect. 
This can, in fact, be used without the burden of having to solve a large linear system. 
In the future, we aim to apply this algorithm to more complicated two dimensional problems, where we expect larger reduction in computational times (as the FOM will become more expensive), and to selectively switch on and off the dispersive term into a shallow water code, in order to allow to pass between different regimes, according to the solution shape.

To our knowledge  this is the first work in this direction, and the initial results are quite promising. Several future improvements and developments  are foreseen both involving  improved reduction strategies, and 
more complex wave dynamics and models. Notable the extension to multiple space dimensions, as well as breaking waves are under investigation.

\section*{Acknowledgments}
D.~T. has been funded by an Inria Postdoc in Team Cardamom and by a SISSA Mathematical Fellowship. M.~R. is a  permanent member of   Inria  Team Cardamom.

Declarations of interest: none.

\appendix
\section{$\mathbb{P}^1$ finite elements for the MS model: full expressions}\label{sec:FEM_MS} 

We report here the expressions of the finite element approximations  of some of the operators arising in the Madsen and S{\o}rensen model.
The first is  $\dispersion\in \R^{N\times N}$, representing  the matrix discretization of the terms $\mathcal{T}^t$. It  can be split into two matrices. The first one is simply defined as a tridiagonal matrix with entries for the $j$th row $B\bar{h}^2_j [1\;\; -2\; 1 ]$, while the second one has the following entries for the $j$th row
\begin{equation}
	-\frac{1}{18} \begin{bmatrix}
		(\bar{h}_j-\bar{h}_{j-1})(2\bar{h}_j+\bar{h}_{j-1}) \\  -(\bar{h}_j-\bar{h}_{j-1})(2\bar{h}_j+\bar{h}_{j-1})+ (\bar{h}_{j+1}-\bar{h}_{j})(2\bar{h}_j+\bar{h}_{j+1})\\ -(\bar{h}_{j+1}-\bar{h}_{j})(2\bar{h}_j+\bar{h}_{j+1})
	\end{bmatrix} ^T.
\end{equation}
The sum of the two defines the dispersion matrix $\dispersion\in \R^{N\times N}$. The term $\mathcal{T}^x[\eta]$ directly depends on $\partial_{xx} \eta$, is discretized using an auxiliary variable $w \approx \partial_{xx} \eta$ with the definition $w_j = \frac{\eta_{j-1} -2 \eta_j + \eta_{j+1}}{\Delta x^2}$. The operator can be approximated by a matrix $\mathbb{T}^x$ with the multiplication $\mathbb{T}^x\boldeta$. It can be conveniently written with two matrices $\mathbb{T}^{x,1}$ and $\mathbb{T}^{x,2}$ such that $\mathbb{T}^{x}\boldeta = \mathbb{T}^{x,1}\boldsymbol{w}+\mathbb{T}^{x,2}\boldsymbol{w}
$, where $\mathbb{T}^{x,1}$ has the following 3 entries on the $j$-th row in the 3 main diagonals
\begin{equation}
	-\frac{\beta g}{3}\begin{bmatrix}
		-\frac{(\bar{h}_j+\bar{h}_{j-1})^3}{4} -\bar{h}_j^3, & \frac{(\bar{h}_j+\bar{h}_{j-1})^3}{4} -\frac{(\bar{h}_j+\bar{h}_{j+1})^3}{4}, & \frac{(\bar{h}_j+\bar{h}_{j+1})^3}{4} +\bar{h}_j^3
	\end{bmatrix},
\end{equation}
and where $\mathbb{T}^{x,2}$ has the following 3 entries on the $j$-th row in the 3 main diagonals
\begin{equation}
	-\frac{\beta g}{6}\begin{bmatrix}
		(\bar{h}_{j}-\bar{h}_{j-1})(\bar{h}_j^2w_j+\frac{1}{4}(\bar{h}_j+\bar{h}_{j-1} )^2) \\ (\bar{h}_{j}-\bar{h}_{j-1})(\bar{h}_j^2w_j+\frac{1}{4}(\bar{h}_j+\bar{h}_{j-1} )^2) +(\bar{h}_{j+1}-\bar{h}_j)(\bar{h}_j^2w_j+\frac{1}{4}(\bar{h}_j+\bar{h}_{j+1} )^2)\\ (\bar{h}_{j+1}-\bar{h}_j)(\bar{h}_j^2w_j+\frac{1}{4}(\bar{h}_j+\bar{h}_{j+1} )^2)
	\end{bmatrix}^T.
\end{equation}
Finally, the discretization of this sponge layer term leads to an operator $\sponge$  which is  a tridiagonal matrix with entries for the $j$th row 
\begin{equation}
	\sponge_{j,j-1} = -\frac{\nu_{j}+\nu_{j-1}}{8\Delta x} ,\,  \sponge_{j,j} = \frac{\nu_{j-1} + 2\nu_{j}+\nu_{j+1}}{8\Delta x} ,\, \sponge_{j,j+1} = -\frac{\nu_{j}+\nu_{j+1}}{8\Delta x}.
\end{equation}

\section{Computational costs of different solvers}\label{app:cost}

\begin{figure}
	\centering
	\subfigure[Ratio between computational time of pdROM with respect to FOM computed either with Thomas solver in \texttt{numba} or with the sparse solver of \texttt{scipy}\label{fig:scipy_times}]{
		\begin{tikzpicture}
			\begin{axis}[
				ymode = log,
				ymin=0.03,ymax=1.4,
				grid=major,
				xlabel={$\NRB$},
				ylabel={ROM time/FOM time},
				legend pos=north west,
				legend style={nodes={scale=0.6, transform shape}},
				width=.45\textwidth
				]
				\newcommand{\test}{Test_periodicSoliton}				
				\addplot[mark=square, color=red] table [x=NROM, y expr=\thisrow{compTimeROM}] {figures/\test/small_edit_vs_ROM.txt};
				\addlegendentry{Test \ref{test1} FOM \texttt{numba}};		
				\addplot[mark=square*, color=red] table [x=NROM, y expr=\thisrow{compTimeROM}] {figures/\test/small_edit_vs_ROM_on_scipy.txt};
				\addlegendentry{Test \ref{test1} FOM \texttt{scipy}};
				\renewcommand{\test}{Test_damBreak}					
				\addplot[mark=o, color=blue] table [x=NROM, y expr=\thisrow{compTimeROM}] {figures/\test/small_edit_vs_ROM.txt};
				\addlegendentry{Test \ref{sec:undular} FOM \texttt{numba}};		
				\addplot[mark=*, color=blue] table [x=NROM, y expr=\thisrow{compTimeROM}] {figures/\test/small_edit_vs_ROM_on_scipy.txt};
				\addlegendentry{Test \ref{sec:undular} FOM \texttt{scipy}};	
				\renewcommand{\test}{Test_periodicSolitonOnBump}				
				\addplot[mark=diamond, color=green] table [x=NROM, y expr=\thisrow{compTimeROM}] {figures/\test/small_edit_vs_ROM.txt};
				\addlegendentry{Test \ref{sec:benchSolitaryTrapezoid} FOM \texttt{numba}};		
				\addplot[mark=diamond*, color=green] table [x=NROM, y expr=\thisrow{compTimeROM}, color=red] {figures/\test/small_edit_vs_ROM_on_scipy.txt};
				\addlegendentry{Test \ref{sec:benchSolitaryTrapezoid} FOM \texttt{scipy}};
			\end{axis}
		\end{tikzpicture}
	}\hfill
	\subfigure[Ratio between computational time of ROM with reduction only on $\Phi$ and pdROM \label{fig:ratio_small_edit}]{
		\begin{tikzpicture}
			\begin{axis}[
				grid=major,
				xlabel={$\NRB$},
				ylabel={Computational time ratios},
				legend pos=north east,
				legend style={nodes={scale=0.6, transform shape}},
				width=.45\textwidth,
				cycle list name=exotic
				]
				\newcommand{\test}{Test_periodicSoliton}				
				\addplot table [x=NROM, y expr=\thisrow{compTimeSmall}/\thisrow{compTimeROM}] {figures/\test/small_edit_vs_ROM.txt};
				\addlegendentry{Test \ref{test1}};
				\renewcommand{\test}{Test_damBreak}				
				\addplot table [x=NROM, y expr=\thisrow{compTimeSmall}/\thisrow{compTimeROM}] {figures/\test/small_edit_vs_ROM.txt};
				\addlegendentry{Test \ref{sec:undular}};
				\renewcommand{\test}{Test_periodicSolitonOnBump}				
				\addplot table [x=NROM, y expr=\thisrow{compTimeSmall}/\thisrow{compTimeROM}] {figures/\test/small_edit_vs_ROM.txt};
				\addlegendentry{Test \ref{sec:benchSolitaryTrapezoid}};
			\end{axis}
		\end{tikzpicture}
	}
	\caption{Computational costs for different algorithms: different FOMs (left), only projection of $\Phi$ for ROM (right)}
\end{figure}
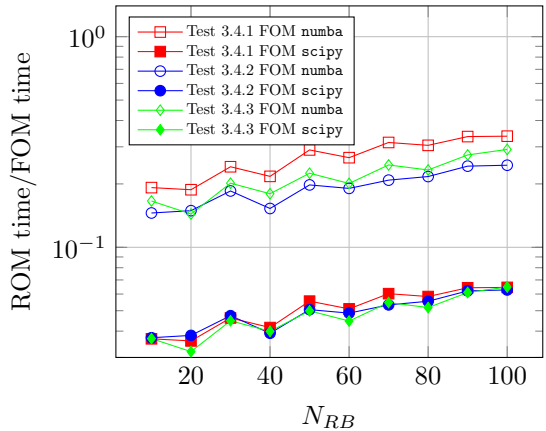
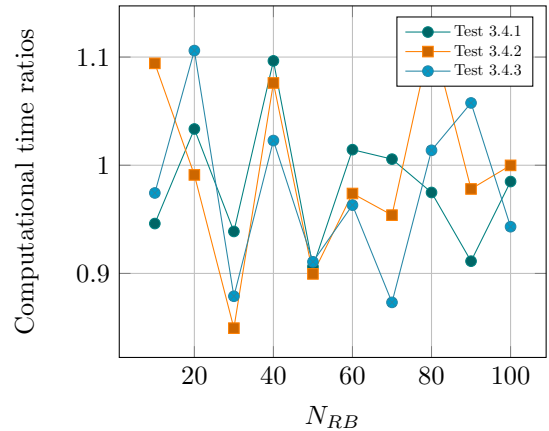
We compare here the computational costs of different sparse linear solvers and different strategies for the reduction. 
We will refer to the KdV-BBM problem \eqref{eq:BBMKdV}, similar results can be drawn from the enhanced Boussinesq system \eqref{eq:EBsystemFOM_SW}.

First of all, let us compare different sparse linear solvers for the FOM. We have used Thomas algorithm in the previously proposed examples. This linear solver very optimized and tailored towards the problem we are solving. As soon as the discretization gets more involved or as soon as we move to more dimensions, we would lose the tridiagonal structure of the system matrices and it would become impossible to apply Thomas algorithm. In order to have a fairer comparison, in \cref{fig:scipy_times} we show the ROM computational costs over the FOM computational costs when the FOM uses Thomas algorithm implemented in \texttt{numba} and when the FOM uses the sparse solver of \texttt{scipy}. We observe a huge difference between the two solvers of a factor between 4 and 6. This would make the pdROM computational costs even more appealing for  a general sparse solver implementation in the FOM,  passing from 20\% of the computational time to only 4\%!

Secondly, we want to prove again that the main computational cost of the FOM is given by the linear system and that the advantage obtained with the pdROM is focused on this step of the algorithm. Indeed, we can compare the performance of the pdROM with the scheme given by the reduction only of the $\Phi$ equation, as suggested in \cref{rem:modification_sw}. In that algorithm, $\Phi$ is computed by projecting the RHS of its equation onto the reduced basis space, then solving the reduced system and, finally, reconstructing the full $\Phi$. Its computational costs are comparable to the one of the pdROM. In \cref{fig:ratio_small_edit} we plot the ratio of these costs for some $\NRB$ for all the previous tests. We can see that they are very close to one and the fluctuations around that are probably due to the low precision in measuring the times.
We remind that the ROM where only $\Phi$ is reduced has larger errors than the pdROM. This is maybe due to a mismatch of the spaces along the computations that allows the propagation of spurious modes. We are still investigating the phenomenon. Hence, in general, we would recommend a fully projected approach.

\bibliographystyle{siam}
\bibliography{biblio}

\end{document}